\documentclass[12pt, draft]{amsart}
\usepackage[russian]{babel}
\usepackage{amsmath,amsfonts,amssymb,euscript,eufrak}

\theoremstyle{plain} \textwidth=170mm \textheight=230mm \hoffset=-20mm \voffset=-15mm

\newtheorem{Lemma}{╦хььр}
\newtheorem{Proposition}{╧Ёхфыюцхэшх}
\newtheorem{Theorem}{╥хюЁхьр}

\newtheorem{Corollary}{╤ыхфёЄтшх}

\newtheorem{Note}{╟рьхўрэшх}

\renewcommand\Re{\mathrm{Re}\,}
\renewcommand\Im{\mathrm{Im\,}}
\renewcommand{\le}{\leqslant}
\renewcommand{\ge}{\geqslant}

\newcommand\wt{\widetilde}

\def\Dom{{\mathfrak D}}

\def\al{\alpha}
\def\la{\lambda}

\def\eps{\varepsilon}

\def\cT{{\mathcal T}}
\def\cE{{\EuScript E}}
\def\cA{{\EuScript A}}
\def\cS{{\mathcal S}}

\def\tr{\operatorname{tr}}

\def\bR{{\mathbb R}}

\def\bC{{\mathbb C}}

\def\C{{\mathcal C}}
\newcommand\AC{\mathrm{AC}}

\def\z{{\mathbf z}}
\def\y{{\mathbf y}}

\def\D{{\mathrm D}}
\def\A{{\mathrm A}}
\def\C{{\mathrm C}}
\def\B{{\mathrm B}}
\def\F{{\mathrm F}}
\def\T{{\mathrm T}}
\def\V{{\mathrm V}}
\def\W{{\mathrm W}}
\def\Y{{\mathrm Y}}
\def\Z{{\mathrm Z}}
\def\M{{\mathrm M}}
\def\O{{\mathrm O}}
\def\I{{\mathrm I}}
\def\E{{\mathrm E}}
\def\Q{{\mathrm Q}}

\def\U{{\mathrm U}}
\def\f{{\mathbf f}}
\def\g{{\mathbf g}}
\def\R{{\mathrm R}}

\def\u{\mathbf{u}}
\def\v{{\mathbf v}}

\def\rv{{\mathrm v}}
\def\i{{\mathfrak i}}
\newlength{\lenun}
\newlength{\lendu}
\def\diag{\operatorname{diag}}

\begin{document}
{\Large \begin{center} └ёшьяЄюЄшўхёъшщ рэрышч  Ёх°хэшщ юс√ъэютхээ√ї фшЇЇхЁхэЎшры№э√ї
єЁртэхэшщ ё ъю¤ЇЇшЎшхэЄрьш-ЁрёяЁхфхыхэш ьш \end{center}}
\bigskip

\centerline{└.~╠.~╤ртўєъ, └.~└.~╪ърышъют\footnote{╚ёёыхфютрэш 
яюффхЁцрэ√ уЁрэЄюь ╨═╘  17-11-01215. This work is supported by Russian Science Foundation under  grant No 17-11-011215.}}
\bigskip
%Asymptotic formulas for  fundamental system of solutions of high order ordinary differential equations with
%coefficients --- distributions.
%Artem M.Savchuk, Andrei A.Shkalikov
%Moscow State University

\section{┬тхфхэшх}

╬ёэютэющ Ўхы№■ эрёЄю ∙хщ ёЄрЄ№ш ртЄюЁ√ ёЄртшыш яюыєўхэшх рёшьяЄюЄшўхёъшї ЇюЁьєы фы  Ёх°хэшщ юс√ъэютхээ√ї фшЇЇхЁхэЎшры№э√ї єЁртэхэшщ
тшфр
\begin{equation}\label{1.1}
\tau(y)- \lambda ^{2m} \sigma(x) y = 0, \quad \tau(y) =\sum_{k,\,s=0}^m(\tau_{k,\,s}(x)y^{(m-k)}(x))^{(m-s)},
\end{equation}
яЁш ьшэшьры№э√ї яЁхфяюыюцхэш ї эр уырфъюёЄ№   ъю¤ЇЇшЎшхэЄют $\tau_{k,\,s}$ ш $\sigma$.  ╟фхё№ $\lambda$ --- ъюьяыхъёэ√щ ярЁрьхЄЁ,
ъюЄюЁ√щ яЁш фюяюыэшЄхы№эющ яюёЄрэютъх ъЁрхт√ї єёыютшщ шуЁрхЄ Ёюы№ ёяхъЄЁры№эюую ярЁрьхЄЁр. └ёшьяЄюЄшўхёъшх ЇюЁьєы√ яЁхфяюырурхЄё  яюыєўшЄ№ яЁш
$\lambda \to \infty$  т  ёхъЄюЁрї ъюьяыхъёэющ яыюёъюёЄш, ъюЄюЁ√х т ёютюъєяэюёЄш яюъЁ√тр■Є тё■ ъюьяыхъёэє■ яыюёъюёЄ№.
┬ ¤Єющ ЁрсюЄх ь√ юуЁрэшўштрхьё  шчєўхэшхь єЁртэхэш  ўхЄэюую яюЁ фър $n= 2m$,   ёыєўрщ эхўхЄэюую яюЁ фър ЄЁхсєхЄ юЄфхы№эюую ЁрёёьюЄЁхэш .

▌Єр ЄхьрЄшър шьххЄ сюыхх ўхь ёЄюыхЄэ■■ шёЄюЁш■, эрўрыю ъюЄюЁющ шэшЎшшЁютрэю т ъырёёшўхёъшї ЁрсюЄрї ─ц.~┴шЁъуюЇр \cite{Bir1}, \cite{Bir2}.  ╚ьх■∙шхё  т эрёЄю ∙хх тЁхь  Ёхчєы№ЄрЄ√ яючтюы ■Є єёЄрэютшЄ№ сшЁъуюЇютёъшх ршёьяЄюЄшъш яЁш  ёыхфє■∙шї єёыютш ї эр ъю¤ЇЇшЎшхэЄ√ фшЇЇхЁхэЎшры№эюую т√Ёрцхэш  $\tau$: \ $\tau_{0,0} \equiv 1, $
р $\tau_{k,s}$  шьх■Є $m-s-1$  рсёюы■Єэю эхяЁхЁ√тэ√ї яЁюшчтюфэ√ї,  ўЄю ¤ътштрыхэЄэю яЁшэрфыхцэюёЄш ЇєэъЎшщ  $\tau_{k,s}$ яЁюёЄЁрэёЄтє ╤юсюыхтр $W^{m-s}_1 [0,1]$
(чфхё№ $W^0_1[0,1] = L_1[0,1]$). ╥ръшх єёыютш  ёыш°ъюь юуЁрэшўшЄхы№э√.
┬ ¤Єющ ёЄрЄ№х ь√ яюыєўшь эєцэ√х ршьяЄюЄшъш эх Єюы№ъю, ъюуфр ъю¤ЇЇшЎшхэЄ√ $\tau_{k,s}$  ты ■Єё  ъырёёшўхёъшьш (ёєььшЁєхь√ьш) ЇєэъЎш ьш, эю ш фы  ёыєўр , ъюуфр
юэш  ты ■Єё  юсюс∙хээ√ьш ЇєэъЎш ьш ъюэхўэюую яюЁ фър ёшэуєы ЁэюёЄш. └ шьхээю, ь√ яюърцхь, ўЄю эєцэ√х рёшьяЄюЄшўхёъшх ЇюЁьєы√  ёяЁртхфышт√ Єюы№ъю яЁш т√яюыэхэшш єёыютшщ
\begin{equation}\label{eq:cond2}
\frac1{\sqrt{|\tau_0|}},\ \frac1{\sqrt{|\tau_0|}}\tau_{k,s}^{(-l)},\in
L_2[0,1], \quad  0\leqslant k,s \leqslant m,  \quad l = \min\{k,s\},
\end{equation}
уфх $\tau_0 : =\tau_{0,0}$, р  $f^{(-k)}$  юсючэрўрхЄ $k$-■ яхЁтююсЁрчэє■ ЇєэъЎшш $f$, яюэшьрхьє■ т ёь√ёых ЄхюЁшш ЁрёяЁхфхыхэшщ. ─Ёєушьш ёыютрьш,
юЄ ъю¤ЇЇшЎшхэЄют $\tau_{k,s}$ яЁш фюяюыэшЄхы№эюь єёыютшш ЁртэюьхЁэющ яюыюцшЄхы№эюёЄш ЇєэъЎшш  $\tau_0$ фюёЄрЄюўэю ЄЁхсютрЄ№, ўЄюс√ юэш с√ыш юсюс∙хээ√ьш яЁюшчтюфэ√ьш $l$-ую яюЁ фър юЄ ътрфЁрЄшўэю шэЄхуЁшЁєхь√ї ЇєэъЎшщ,
уфх $l= \min \{k,s\}$.  ╩юэхўэю, яюыєўхэшх Єръюую Ёхчєы№ЄрЄр яюЄЁхсєхЄ эх Єюы№ъю ёхЁ№хчэющ Єхїэшўхёъющ ЁрсюЄ√ (ъюЄюЁє■ ь√ я√Єрышё№ ьръёшьры№эю єяЁюёЄшЄ№),
эю ш яЁшьхэхэш  эют√ї яЁшхьют ш ьхЄюфют. ┬ ўрёЄэюёЄш, юфэшь шч эршсюыхх ёє∙хёЄтхээ√ї ьюьхэЄют яЁш ЁхрышчрЎшш яюёЄртыхээющ Ўхыш  ты хЄё  шёяюы№чютрэшх
 \textit{ьхЄюфр Ёхуєы ЁшчрЎшш}. ─ы  єЁртэхэш  ╪ЄєЁьр-╦шєтшыы   ¤ЄюЄ ьхЄюф с√ы шэшЎшшЁютрэ т ЁрсюЄрї ртЄюЁют \cite{SavSh99}, \cite{SavSh03}.  ┬ яюёыхфє■∙шї ЁрсюЄрї \cite{Vlad} -- \cite{Mir1} ьхЄюф с√ы ЁрёяЁюёЄЁрэхэ фы  юс∙шї єЁртэхэшщ тЄюЁюую яюЁ фър ш фы  эхъюЄюЁ√ї ёяхЎшры№э√ї ъырёёют єЁртэхэшщ ўхЄтхЁЄюую ш т√ёюъшї яюЁ фъют.  ╟фхё№ сєфхЄ шёяюы№чютрэ эхфртэшщ Ёхчєы№ЄрЄ ╠шЁчюхтр ш ╪ърышъютр
\cite{MirzSh}, уфх яЁютхфхэр Ёхуєы ЁшчрЎш  юс∙шї єЁртэхэшщ \eqref{1.1}.

 ╚ёяюы№чютрэшх ьхЄюфр Ёхуєы ЁшчрЎшш тхфхЄ ъ шфхх  яюыєўхэш  рёшьяЄюЄшўхёъшї ЇюЁьєы яєЄхь ётхфхэш  єЁртэхэш  \eqref{1.1} ъ ёшёЄхьх $2m$  єЁртэхэшщ яхЁтюую яюЁ фър. ═ю шьх■∙шхё   Ёхчєы№ЄрЄ√ фы  ёшёЄхь эх тяюыэх эрь яюфїюф Є, ш ь√ яЁшїюфшь ъ тЄюЁющ Ўхыш --- яюыєўхэш■ рёшьяЄюЄшўхёъшї ЇюЁьєы фы  ьрЄЁшЎ√ ЇєэфрьхэЄры№э√ї Ёх°хэшщ  ёшёЄхь√ яхЁтюую яюЁ фър тшфр
\begin{equation}\label{eq:main}
\y'=\la\V(x)\y+\A(x)\y+\C(x,\la)\y, \quad x\in[0,1],
\end{equation}
ё ъюьяыхъёэ√ь  сюы№°шь ярЁрьхЄЁюь $\la$. ╟фхё№
$$
\y=\y(x)=(y_1(x),y_2(x),\dots,y_n(x))^t
$$
--- тхъЄюЁ-ёЄюысхЎ, ёюёЄртыхээ√щ шч рсёюы■Єэю эхяЁхЁ√тэ√ї эр $[0,1]$ ЇєэъЎшщ
$y_j(x)$.
%┬ ърўхёЄтх ьэюцхёЄтр чэрўхэшщ ярЁрьхЄЁр $\la$ ь√ сєфхь ЁрёёьрЄЁштрЄ№ юсырёЄ№ тшфр $\{\la\in\bC:|\la|>\la_0\}$.
╤ЇююЁьєышЁєхь єёыютш , яЁш ъюЄюЁ√ї ь√ сєфхь шчєўрЄ№ ¤Єє ёшёЄхьє.
%┴єфхь ёўшЄрЄ№, ўЄю ъю¤ЇЇшЎшхэЄ√ ёшёЄхь√ \eqref{eq:main} єфютыхЄтюЁ ■Є ёыхфє■∙шь єёыютш ь.

\hfill\break (i)\quad ╠рЄЁшЎр $\V(x)$  шьххЄ тшф  $\V(x) =\rho(x)\B$, уфх $\B=\diag\{b_1,\dots,b_n\}$ ---   фшруюэры№эр , яюёЄю ээр 
ьрЄЁшЎр, р ъюьяыхъёэ√х ўшёыр
$b_j$ юЄышўэ√ юЄ эєы . ╘єэъЎш■ $\rho$  ёўшЄрхь ёєььшЁєхьющ эр $[0,1]$ ш яюыюцшЄхы№эющ яюўЄш тё■фє. ┼х яхЁтююсЁрчэє■ юсючэрўрхь
$p(x):=\int_0^x\rho(t)\,dt$.
\hfill\break
(ii)\quad ┬ёх ъю¤ЇЇшЎшхэЄ√ ёшёЄхь√
--- ЇєэъЎшш $a_{jk}(x)$ ш $c_{jk}(x,\la)$ яхЁхьхээющ $x$ (яЁш ърцфюь ЇшъёшЁютрээюь $\la$)  ёўшЄрхь
ёєььшЁєхь√ьш яю ╦хсхує эр тёхь юЄЁхчъх $[0,1]$.
%\hfill\break (iii)\quad ┴єфхь ёўшЄрЄ№, ўЄю ЇєэъЎш  $\rho(x)$ тх∙хёЄтхээр
%ш яюыюцшЄхы№эр яюўЄш тё■фє. ╧хЁтююсЁрчэє■ ЇєэъЎшш $\rho(x)$ сєфхь юсючэрўрЄ№ ўхЁхч
%$p(x):=\int_0^x\rho(t)\,dt$.
\hfill\break (iii)\quad ╧юырурхь, ўЄю т $L_1$- эюЁьх ЇєэъЎшш
$c_{jk}(x,\la)=o(1)$  яЁш $|\la|\to\infty$, р шьхээю, $\int_0^1|c_{jk}(x,\la)|\,dx\to 0$ фы  тёхї $1\le j,\,k\le n$.

┬ фхщёЄтшЄхы№эюёЄш,  ъ ёшёЄхьх \eqref{eq:main} ётюф Єё  ёшёЄхь√ ё ьрЄЁшЎхщ $\V(x)$  сюыхх юс∙хую тшфр,  эхцхыш фшруюэры№э√х ьрЄЁшЎ√,
єърчрээ√х т єcыютшш (i). ╧юфЁюсэхх юс ¤Єюь ёърцхь т эрўрых ёыхф■∙хую ярЁруЁрЇр.% ═ю ЁрсюЄрЄ№ єфюсэхх ё ьрЄЁшЎрьш, яюфўшэхээ√ьш єёыютш■ (i).

╠рЄЁшЎхщ ЇєэфрьхэЄры№эющ ёшёЄхь√ Ёх°хэшщ шыш ьрЄЁшЎхщ ьюэюфЁюьшш ь√ эрч√трхь ьрЄЁшЎє $\Y(x,\la)$ Ёрэур $n$,
юяЁхфхыхээє■ яЁш $x\in[0,1]$, $\la\in\bC$, $|\la|>\la_0$, ърцф√щ ёЄюысхЎ ъюЄюЁющ $\Y_k(x,\la)$  ты хЄё  Ёх°хэшхь
ёшёЄхь√ \eqref{eq:main}. ╥ръшь юсЁрчюь, ьрЄЁшЎр $\Y$ єфютыхЄтюЁ хЄ ьрЄЁшўэюьє фшЇЇхЁхэЎшры№эюьє єЁртэхэш■
\begin{equation}\label{eq:matrmain}
\Y'(x,\la)=\left(\la\rho(x) \B+\A(x)+\C(x,\la)\right)\Y(x,\la)
\end{equation}
яЁш ърцфюь ЇшъёшЁютрээюь $\la$. ╒юЁю°ю шчтхёЄэю (ёь., эряЁшьхЁ, \cite{CL}), ўЄю фы  юяЁхфхышЄхы  $\det \Y(x,\la)$
т√яюыэхэю ЁртхэёЄтю
$$
(\det \Y(x,\la))'=\tr\left(\la\rho(x)\B+\A(x)+\C(x,\la)\right)\det \Y(x,\la),
$$
Єръ ўЄю шч єёыютш  $\det \Y(\xi,\la)\ne0$ фы  эхъюЄюЁюую $\xi\in[0,1]$ ёыхфєхЄ, ўЄю $\det \Y(x,\la)\ne0$ яЁш тёхї
$x\in[0,1]$.

╩ръ єцх ёърчрэю, эр°хщ тЄюЁющ Ўхы№■  ты хЄё   яюыєўхэшх рёшьяЄюЄшўхёъшї яЁхфёЄртыхэшщ фы  ьрЄЁшЎ√ $\Y(x,\la)$ яЁш $|\la|\to\infty$.
▌Єр Єхьр ъырёёшўхёър  ш шьххЄ фртэ■■ шёЄюЁш■. ╧юёых ЁрсюЄ ┴шЁъуюЇр \cite{Bir1}, \cite{Bir2} ш ╧хЁЁюэр \cite{Per}
  яю тшырё№ ьюэюуЁрЇш  ╥рьрЁъшэр  \cite{Ta1} ш ёЄрЄ№   ┴шЁъуюЇр ш ╦рэухЁр
\cite{BirLa}, т ъюЄюЁ√ї с√ыю яЁхфяЁшэ Єю шчєўхэшх рёшьяЄюЄшъ ЇєэфрьхэЄры№э√ї ьрЄЁшЎ Ёх°хэшщ фы  ёшёЄхь ш шчєўхэшх
ёяхъЄЁры№э√ї чрфрў, яюЁюцфрхь√ї  ёшёЄхьрьш ё фюсртыхээ√ьш ъЁрхт√ьш єёыютш ьш. ▌Єш фтр шёЄюўэшър эх яюЄхЁ ыш ръЄєры№эюёЄ№
тяыюЄ№ фю эрёЄю ∙хую тЁхьхэш.  ╨хчєы№ЄрЄ√ \cite{Ta1} с√ыш тяюёыхфёЄтшш фюяюыэхэ√ т \cite{Ta2}.
  ╒юЁю°ю шчтхёЄэю, ўЄю  рёшьяЄюЄшўхёъшх Ёхчєы№ЄрЄ√
фы  ёшёЄхь  ты ■Єё  ъы■ўхт√ьш фы  шёёыхфютрэш  ёяхъЄЁры№э√ї
ётющёЄт юс√ъэютхээ√ї фшЇЇхЁхэЎшры№э√ї юяхЁрЄюЁют яюЁ фър $n\ge2$ (ёь., эряЁшьхЁ, ьюэюуЁрЇш■ ═рщьрЁър \cite{Na}).
═ю шёёыхфютрэшх ёшёЄхь, эхёюьэхээю шьххЄ ёрьюёЄю Єхы№эюх чэрўхэшх.  ┬ ўрёЄэюёЄш, юфшэ шч трцэхщ°шї юяхЁрЄюЁют ьрЄхьрЄшўхёъющ Їшчшъш
--- юяхЁрЄюЁ ─шЁрър яюЁюцфрхЄё  ёшёЄхьющ тЄюЁюую яюЁ фър (ўрёЄю  т ёыєўрх $n=2m$, $b_1=\dots=b_m=\i$,
$b_{m+1}=\dots=b_n=-\i$, ёююЄтхЄёЄтє■∙є■ ёшёЄхьє  т ышЄхЁрЄєЁх  Єръцх   эрч√тр■Є \textit{ёшёЄхьющ Єшяр ─шЁрър}).
╩ ёшёЄхьх ─шЁрър ётюфшЄё  єЁртэхэшх  ╪ЄєЁьр-╦шєтшыы 
$-y''+qy=\la^2\varrho y$. ╚ёёыхфютрэш  ёяхъЄЁры№э√ї чрфрў фы  ¤Єюую єЁртэхэш  т ёшэуєы Ёэюь ёыєўрх фы   яюЄхэЎшрыют $q\in W^{-1}_2[0,1]$
с√ыю яЁютхфхэю т ЁрсюЄх ртЄюЁют \cite{SavSh03} ё яюью∙№■ рёшьяЄюЄшўхёъшї ЇюЁьєы фы  ЇєэфрьхэЄры№э√ї Ёх°хэшщ, ъюЄюЁ√х с√ыш єёЄрэютыхэ√ ёэрўрыр Єюы№ъю т ъЁшЄшўхёъшї яюыюёрї. ╧ючфэхх  рёшьяЄюЄшўхёъшх ЇюЁьєы√ фы  юс∙шї єЁртэхэшщ тЄюЁюую яюЁ фър ё ёшэуєы Ёэ√ьш ъю¤ЇЇшЎшхэЄрьш
т яюыєяыюёъюёЄ ї ъюьяыхъёэющ яыюёъюёЄш с√ыш яюыєўхэ√ ┬ырф√ъшэющ ш ╪ърышъют√ь \cite{VlaSh}. ╧юьшью рёшьяЄюЄшўхёъшї
ьхЄюфют фы  юяхЁрЄюЁют ╪ЄєЁьр-╦шєтшыы  ¤ЇЇхъЄштэю ЁрсюЄрхЄ ьхЄюф \textit{юяхЁрЄюЁр яЁхюсЁрчютрэш } (ёь. ъырёёшўхёъє■ ьюэюуЁрЇш■ \cite{Ma} ш эют√х ЁрсюЄ√
└ы№схтхЁшю, ├Ёшэштр ш ╠шъшЄ■ър \cite{HM}, \cite{AHM}, уфх ¤ЄюЄ ьхЄюф ЁрчтшЄ т ёшэуєы Ёэюь ёыєўрх).

%╧Ёшўшэє ьюцэю юс· ёэшЄ№ Єхь, ўЄю т ¤Єюь ёыєўрх
%тэхфшруюэры№эр  ўрёЄ№ ьрЄЁшЎ√ $\A$ рэЄшъюььєЄшЁєхЄ ё ьрЄЁшЎхщ $\B$
%$$
%\begin{pmatrix}\i&0\\0&-\i\end{pmatrix}\begin{pmatrix}0&a_{12}\\ a_{21}&0\end{pmatrix}=
%-\begin{pmatrix}0&a_{12}\\ a_{21}&0\end{pmatrix}\begin{pmatrix}\i&0\\0&-\i\end{pmatrix}.
%$$
%═р°ш яюёЄЁюхэш  т фрээющ ЁрсюЄх сєфєЄ юс∙шьш, ъЁюьх єёыютшщ (i)--(iv), ь√ эх сєфхь фхырЄ№ эшъръшї яЁхфяюыюцхэшщ ю
%ъю¤ЇЇшЎшхэЄрї ёшёЄхь√.

 ┬  ётюхь шёёыхфютрэшш \cite{Ta1} ш хую фюяюыэхэшш \cite{Ta2} ╥рьрЁъшэ яЁхфяюыруры, ўЄю яюёых
фшруюэрышчрЎшш ьрЄЁшЎр $\V$ ёшёЄхь√ \eqref{eq:main} шьххЄ тшф $\B=\B(x)=\diag\{\varphi_1(x),\dots,\varphi_n(x)\}$  ё фюяюыэшЄхы№э√ьш єёыютш ьш
эр чэрўхэш  ЇєэъЎшщ $\varphi_j$,  ъюЄюЁ√х ё єўхЄюь ╧Ёхфыюцхэш  \ref{pr:gen}  ¤ътштрыхэЄэ√ эр°хьє єёыютш■ (i).  ╧Ёш ¤Єюь ¤ыхьхэЄ√ ьрЄЁшЎ√ $\V(x)$,
 $\A(x)$ ш $\C(x,\la)$ яЁхфяюырурышё№ фтрцф√ эхяЁхЁ√тэю фшЇЇхЁхэЎшЁєхь√ьш, эхяЁхЁ√тэю фшЇЇхЁхэЎшЁєхь√ьш ш эхяЁхЁ√тэ√ьш ёююЄтхЄёЄтхээю. ╧Ёхфяюырурыюё№
 Єръцх єс√трэшх эр схёъюэхўэюёЄш ьрЄЁшўэ√ї ¤ыхьхэЄют  $\C(x,\la)$ ё яюЁ фъюь
$O(\lambda^{-1})$. ┬яюёыхфёЄтшш єёыютш  уырфъюёЄш
эхюфэюъЁрЄэю юёырсы ышё№.  ┬ ьюэюуЁрЇшш ╨ряюяюЁЄр \cite{Rap}  рёшьяЄюЄшўхёъшх
ётющёЄтр Ёх°хэшщ с√ыш шчєўхэ√ фы   ёшёЄхь тшфр \eqref{eq:main} ё тх∙хёЄтхээ√ь ярЁрьхЄЁюь $\la$ т яЁхфяюыюцхэшш  яЁшэрфыхцэюёЄш ¤ыхьхэЄют ьрЄЁшЎ $\V(x)$,
$\A(x)$ ш $\C(x,\la)$ яЁюёЄЁрэёЄтрь  $W_1^2[0,1]$,\ $\AC[0,1] = W^1_1[0,1]$ ш $\C[0,1]$  ёююЄтхЄёЄтхээю ё фюяюыэшЄхы№э√ь єёыютшхь єс√трэш  $\C(x,\la)=O(\la^{-1})$.
 ┬ ёЄрЄ№х  ┬рурсютр \cite{Vag} яЁютхфхэю фры№эхщ°хх юёырсыхэшх єёыютшщ уырфъюёЄш: $\V(x)\in\C^1[0,1]$, $\A(x)\in\C[0,1]$.
% ╥ръюую Ёюфр єёыютш  эр уырфъюёЄ№  ты ■Єё  ёЄрэфрЁЄэ√ьш фы 
%ъырёёшўхёъющ ЄхюЁшш (ёь. Єръцх \cite[уыртр II, \S 1]{Na}).
╧ю-тшфшьюьє, эршсюыхх юс∙шщ эр Єхъє∙шщ ьюьхэЄ ёыєўрщ  фы  ёшёЄхь тшфр \eqref{eq:main}
ЁрчюсЁрэ ╨√їыют√ь т ёЄрЄ№х \cite{Ryhl}, уфх ¤ыхьхэЄ√ ьрЄЁшЎ√ $\V$ яЁхфяюырур■Єё  рсёюы■Єэю эхяЁхЁ√тэ√ьш, р
¤ыхьхэЄ√ ьрЄЁшЎ $\A(x)$ ш $\C(x)$ --- ёєььшЁєхь√ьш яю ╦хсхує. ╬ЄьхЄшь, ўЄю ёЄрЄ№  \cite{Ryhl} юёЄрырё№ эхчрьхўхээющ
ёяхЎшрышёЄрьш, ┬ ўрёЄэюёЄш, т ЁрсюЄх \cite{Vag1} хх Ёхчєы№ЄрЄ√ с√ыш яхЁхфюърчрэ√, яЁшўхь яЁш сюыхх ёшы№э√ї єёыютш ї
эр ъю¤ЇЇшЎшхэЄ√:  $\V\in W_q^1[0,1]$, $\A(x),\,\C(x)\in L_q[0,1]$, $q>1$, р юЎхэър юёЄрЄър с√ыр яюыєўхэр
т сюыхх ёырсющ эюЁьх $\|\cdot\|_{L_{q'}}$, $1/q+1/q'=1$, тьхёЄю
$\|\cdot\|_{L_\infty}$.

%─ы  эрё юёэютэ√ь яюсєфшЄхы№э√ь ьюЄштюь ъ шёёыхфютрэш■ ёшёЄхь тшфр \eqref{eq:main}  ты хЄё  Єю,
%ўЄю ъ эшь ётюф Єё  ышэхщэ√х фшЇЇхЁхэЎшры№э√х єЁртэхэш  т√ёюъюую яюЁ фър ё ъю¤ЇЇшЎшхэЄрьш
%--- ЁрёяЁхфхыхэш ьш тшфр \eqref{1.1}. ╧Ёш ¤Єюь ¤ыхьхэЄ√ ьрЄЁшЎ $\A$ ш $\C$ юърч√тр■Єё 
%ыш°№ ёєььшЁєхь√ьш эр $[0,1]$, ўЄю ш т√ч√трхЄ эхюсїюфшьюёЄ№ юёырсыхэш  ёЄрэфрЁЄэ√ї ЄЁхсютрэшщ уырфъюёЄш.

╟фхё№ ь√ юуЁрэшўшышё№ ъюЁюЄъющ ёяЁртъющ яю Єхьх эр°хщ ЁрсюЄ√.  ╧юфЁюсэхх ю Ёхчєы№ЄрЄрї рёьшьяЄюЄшўхёъющ ЄхюЁшш фы  юс√ъэютхээ√ї фшЇЇхЁхэЎшры№э√ї єЁртэхэшщ ш ёшёЄхь ш сюыхх яюыэ√ь яхЁхўэхь  ышЄхЁрЄєЁ√ яю Єхьх ьюцэю ючэръюьшЄ№ё  т ьюэюуЁрЇш ї ╨ряюяюЁЄр \cite{Rap}, ┬рчютр \cite{Vaz}, ═рщьрЁър \cite{Na}, ёЄрЄ№ ї ╪ърышъютр \cite{Sh1}, ╨√їыютр \cite{Ryhl}, ╠рырьєфр ш ╬ЁшфюЁюуш \cite{MaO}. ─рыхх яюыхчэю юЄьхЄшЄ№ юёэютэ√х ¤ыхьхэЄ√ эютшчэ√ ¤Єющ эр°хщ ЁрсюЄ√ т ёЁртэхэшш ё яЁхф°хёЄтє■∙шьш. ╤шЄєрЎш  ё рёшьяЄюЄшўхёъшь рэрышчюь Ёх°хэшщ фы  єЁртэхэш  \eqref{1.1}  ёэр : Ёрэхх Єръющ рэрышч яЁш эрышўшш ъю¤ЇЇшЎшхэЄют-ЁрёяЁхфхыхэшщ яЁютюфшыё  Єюы№ъю фы  єЁртэхэшщ тЄюЁюую яюЁ фър, яЁшўхь фЁєушьш ьхЄюфрьш, ъюЄюЁ√х фы  т√ёюъшї яюЁ фъют ЁхрышчютрЄ№ эрь эх єфрыюё№. ╧юфЁюсэхх ёърцхь ю эютшчэх
яЁш рэрышчх ёшёЄхь тшфр \eqref{eq:main}.
%└ёшьяЄюЄшўхёъшщ рэрышч Ёх°хэшщ ¤Єшї ёшёЄхь \eqref{eq:main} шьххЄ Ё ф юўхэ№ трцэ√ї Єюэъшї ьюьхэЄют.
╒юЁю°ю шчтхёЄэю, ўЄю рёшьяЄюЄшўхёъюх
яютхфхэшх ьрЄЁшЎ√ $\Y(x,\la)$ яЁш $|\la|\to+\infty$ ърЁфшэры№э√ь юсЁрчюь чртшёшЄ юЄ ьэюцхёЄтр $\Omega$ т ъюьяыхъёэющ
яыюёъюёЄш, ёюфхЁцр∙хую ярЁрьхЄЁ $\la$. ╙цх схуы√щ тчуы ф эр ёшёЄхьє \eqref{eq:main} яючтюы хЄ ёфхырЄ№ т√тюф, ўЄю
уыртэ√ьш ўыхэрьш рёшьяЄюЄшъш фюыцэ√ с√Є№ ЇєэъЎшш тшфр $\exp\{\la b_j p(x)\}, \ p(x) = \int \rho(t)\, dt, $ рёшьяЄюЄшўхёъюх яютхфхэшх ъюЄюЁ√ї (ЁюёЄ,
єс√трэшх шыш юёЎшыы Ўш ) юяЁхфхы хЄё  чэръюь тхышўшэ√ $\Re(\la b_j p(x))$. ╙ўшЄ√тр  эхюЄЁшЎрЄхы№эюёЄ№ ЇєэъЎшш $p$,
тшфшь, ўЄю рёшьяЄюЄшъє ЇєэъЎшш $\Y(x,\la)$ эрфю шёърЄ№ т ёхъЄюЁрї $\arg\la\in(\al,\beta)$ т ъюьяыхъёэющ яыюёъюёЄш. ╧Ёш
¤Єюь тэєЄЁш ёхъЄюЁр сєфхЄ эрсы■фрЄ№ё  ышсю ¤ъёяюэхэЎшры№э√щ ЁюёЄ, ышсю ¤ъёяюэхэЎшры№эюх єс√трэшх Ёх°хэшщ. ╬ёЎшыы Ўш 
тючэшърхЄ ыш°№ т яюыєяюыюёрї ъюьяыхъёэющ яыюёъюёЄш, ъюЄюЁ√х ёюфхЁцрЄ уЁрэшЎ√ ёхъЄюЁют. ▌Єш уЁрэшЎ√  ь√ фрыхх сєфхь эрч√трЄ№
\textit{ъЁшЄшўхёъшьш  ыєўрьш}, р яЁюшчтюы№э√х яюыєяюыюё√, ёюфхЁцр∙шх Єръшх ыєўш, \textit{ъЁшЄшўхёъшьш  яюыєяюыюёрьш}.
 ╬фэръю чрфрўш юс рёшьяЄюЄшўхёъюь яютхфхэшш ёюсёЄтхээ√ї чэрўхэшщ фшЇЇхЁхэЎшры№э√ї
юяхЁрЄюЁют т√ёюъюую яюЁ фър шыш юяхЁрЄюЁют, эхяюёЁхфёЄтхээю яюЁюцфрхь√ї фшЇЇхЁхэЎшры№э√ь т√Ёрцхэшхь \eqref{eq:main},
ЄЁхсє■Є чэрэш  рёшьяЄюЄшъш $\Y(x,\la)$ шьхээю т ъЁшЄшўхёъшї яюыєяюыюёрї. ┬ ёшыє ¤Єшї яЁшўшэ ь√ эх ьюцхь шёяюы№чютрЄ№
Ёхчєы№ЄрЄ√ эхфртэхщ ЁрсюЄ√ ╠рырьєфр ш ╬ЁшфюЁюуш \cite{MaO}, т ъюЄюЁющ шчєўрышё№ ёшёЄхь√ тшфр \eqref{eq:main} c
$\C(x,\la)=0$, эю рёшьяЄюЄшўхёъшх яЁхфёЄртыхэш  с√ыш яюыєўхэ√ т "<ёєцхээ√ї"> ёхъЄюЁрї
$\arg\la\in(\al+\eps,\beta-\eps), \ \eps>0$. ╥ръшх ёхъЄюЁ√ эх яюъЁ√тр■Є тё■ ъюьяыхъёэє■ яыюёъюёЄ№, ш яюыєўшЄ№ шэЇюЁьрЎш■
ю яютхфхэшш Ёх°хэшщ т ъЁшЄшўхёъшї яюыєяюыюёрї шч Ёхчєы№ЄрЄют \cite{MaO} эх єфрхЄё .

╨хчєы№ЄрЄ√ ЁрсюЄ√ \cite{Ryhl}, т ъюЄюЁющ эрщфхэю рёшьяЄюЄшўхёъюх яютхфхэшх ЇєэъЎшш $\Y(x,\la)$ т чрьъэєЄ√ї ёхъЄюЁрї
$\arg\la\in[\al,\beta]$, Єръцх юърчрышё№ эхфюёЄрЄюўэ√ьш фы  эр°шї эєцф. ╧Ёшўшэ ъ ¤Єюьє ЄЁш, яЁшўхь юфэр шч эшї
ёє∙хёЄтхээр. ╧хЁтр  ёюёЄюшЄ т Єюь, ўЄю Єръюую Ёюфр рёшьяЄюЄшўхёъшх ЇюЁьєы√  эхюсїюфшью шьхЄ№ т сюыхх °шЁюъшї юсырёЄ ї.
 ╟рЇшъёшЁєхь фы  юяЁхфхыхээюёЄш ыєў $\arg\la=\al$ ш яЁхфяюыюцшь, ўЄю юэ  ты хЄё  уЁрэшЎхщ фтєї ёюёхфэшї ёхъЄюЁют. ╥юуфр, эрщфхээ√х т
\cite{Ryhl} рёшьяЄюЄшъш фы  ьрЄЁшЎ√ Ёх°хэшщ $\Y_+(x,\la)$, ёяЁртхфышт√ тэєЄЁш ёхъЄюЁр $\arg\la\in[\al,\beta]$, Є.х. т
яюыютшэх ъЁшЄшўхёъющ яюыєяюыюё√, эряЁртыхээющ тфюы№ ыєўр $\arg\la=\al$. ╩юэхўэю, т фЁєующ яюыютшэх ¤Єющ яюыєяюыюё√
рёшьяЄюЄшър Ёх°хэш  Єръцх юсхёяхўштрхЄё  ЄхюЁхьющ ЁрсюЄ√ \cite{Ryhl}, эю ¤Єє рёшьяЄюЄшъє сєфхЄ шьхЄ№ \textit{єцх
фЁєур } ьрЄЁшЎр Ёх°хэшщ $\Y_-(x,\la)$. ╥ръшь юсЁрчюь, фы  яюыєўхэш  рёшьяЄюЄшъш ьрЄЁшЎ√ Ёх°хэшщ тэєЄЁш тёхщ ъЁшЄшўхёъющ
яюыєяюыюё√ эрфю яЁютюфшЄ№ "<ёъыхщъє"> Ёх°хэшщ, шчєўрЄ№ рёшьяЄюЄшъє ьрЄЁшЎ яхЁхїюфр ш Є.ф. ┬ эр°хщ ЁрсюЄх ь√ яюыєўшь
рёшьяЄюЄшўхёъшх ЇюЁьєы√ фы  $\Y(x,\la)$ т "<Ёрё°шЁхээ√ї"> ёхъЄюЁрї, яюыєўхээ√ї ёфтшуюь шёїюфэ√ї ёхъЄюЁют тфюы№
яЁюфюыцхэш  сшёёхъЄЁшё√. ╧Ёш ¤Єюь ь√ фюяєёърхь ёфтшу эр яЁюшчтюы№эюх ЁрёёЄю эшх, ўЄю яючтюы хЄ юїтрЄшЄ№ ъЁшЄшўхёъшх
яюыєяюыюё√ яЁюшчтюы№эющ °шЁшэ√.

┬ЄюЁр  (ёє∙хёЄтхээр ) яЁюсыхьр ёюёЄюшЄ т т√сюЁх эюЁь√ фы  юёЄрЄюўэ√ї ўыхэют. ├Ёєсю уютюЁ , ¤Єш ўыхэ√ шьх■Є тшф
$\int_0^xe^{\i rt}a(t)\,dt$, уфх $r\to+\infty$, ЇєэъЎш  $a(t)$ ёєььшЁєхьр эр юЄЁхчъх $[0,1]$, р Єюўъє $x\in[0,1]$, т
ъюЄюЁющ ш∙хЄё  чэрўхэшх Ёх°хэш  $\Y(x,\la)$, ьюцэю ёўшЄрЄ№ ЇшъёшЁютрээющ. ╬фэръю ёяхЎшЇшър ьхЄюфр яюёыхфютрЄхы№э√ї
яЁшсышцхэшщ, ъюЄюЁ√щ хёЄхёЄтхээю шёяюы№чєхЄё  яЁш яюшёъх Ёх°хэш  ёшёЄхь√ \eqref{eq:main}, ЄЁхсєхЄ юЎхэъш тхышўшэ тшфр
$\int_0^\xi e^{\i rt}a(t)\,dt$, уфх $\xi\in[0,x]$. ╚ чфхё№ єцх Єюўъє $\xi$ ЇшъёшЁютрээющ ёўшЄрЄ№ эхы№ч , эхюсїюфшью
юЎхэштрЄ№ ъръє■-Єю ЇєэъЎшюэры№эє■ эюЁьє юёЄрЄър, чртшё ∙хую юЄ яхЁхьхээющ $\xi$. ┬√сюЁ ¤Єющ эюЁь√ юърч√трхЄё 
ўЁхчт√ўрщэю трцэ√ь. ─ы  ёыєўр  $a\in L_1[0,1]$ хёЄхёЄтхээющ  ты хЄё  эюЁьр $\|\cdot\|_{L_\infty}$ фы  юёЄрЄъют (шьхээю т
¤Єющ эюЁьх яюыєўхэ√ юЎхэъш т ЁрсюЄх \cite{Ryhl}). ╦хуъю тшфхЄ№, ўЄю яюЄЁхсютрт фюяюыэшЄхы№эє■ уырфъюёЄ№ юЄ ЇєэъЎшш
$a(x)$ (т√Ёрцхээє■, эряЁшьхЁ, т ЄЁхсютрэш ї эр шэЄхуЁры№э√щ ьюфєы№ эхяЁхЁ√тэюёЄш), ь√ яюыєўшь єыєў°хэшх юЎхэюъ юёЄрЄъют
т рёшьяЄюЄшўхёъшї ЇюЁьєырї. ╬фэръю яюЄЁхсютрт юЄ ЇєэъЎшш $a$ сюы№°хщ ёЄхяхэш ёєььшЁєхьюёЄш, Є.х. яюыюцшт $a\in
L_q[0,1]$, $q>1$, ь√ эх фюс№хьё  єыєў°хэш  юЎхэюъ юёЄрЄъют. ─ы  эрё ¤ЄюЄ тюяЁюё шьххЄ сюы№°юх чэрўхэшх, яюёъюы№ъє яЁш
ётхфхэшш фшЇЇхЁхэЎшры№эюую єЁртэхэш  т√ёюъюую яюЁ фър ё ъю¤ЇЇшЎшхэЄрьш--ЁрёяЁхфхыхэш ьш ъ ёшёЄхьх, ¤ыхьхэЄ√ ьрЄЁшЎ√
$\A$ юърч√тр■Єё  ътрфЁрЄшўэю ёєььшЁєхь√ьш. ╚ёяюы№чютрэшх чфхё№ эюЁь√ $\|\cdot\|_{L_\infty}$ фы  юЎхэштрэш  юёЄрЄъют
яЁштхыю с√ ъ ёхЁ№хчэюьє єїєф°хэш■ юЎхэюъ т рёшьяЄюЄшўхёъшї ЇюЁьєырї фы  ёюсёЄтхээ√ї чэрўхэшщ ш ёюсёЄтхээ√ї ЇєэъЎшщ.
═ряЁшьхЁ, фы  юяхЁрЄюЁр ╪ЄєЁьр--╦шєтшыы  тьхёЄю фюърчрээющ т \cite{SavSh99} ш \cite{SavSh03} рёшьяЄюЄшъш
$\sqrt{\la_n}=n+s_n$, уфх $\{s_n\}_1^\infty\in l_2$, %(т ёырсю Ёхуєы Ёэюь ёыєўрх $\{s_n^2\}_1^\infty\in l_2$),
ь√ яюыєўшыш с√ ыш°№ юЎхэъє $s_n=o(1)$.

═ръюэхЎ, ЄЁхЄ№  яЁшўшэр ёюёЄюшЄ т Єюь, ўЄю т \cite{Ryhl} эрщфхэ ыш°№ уыртэ√щ ўыхэ рёшьяЄюЄшъш. ╩юэхўэю, фхщёЄтє 
ьхЄюфюь яюёыхфютрЄхы№э√ї яЁшсышцхэшщ, ьюцэю т√яшёрЄ№ ы■сюх эхюсїюфшьюх ўшёыю ёырурхь√ї т рёшьяЄюЄшўхёъшї яЁхфёЄртыхэш ї
(т ЄхюЁхьх \ref{tm:app} фрээющ ЁрсюЄ√ ь√ фюїюфшь фю тЄюЁюую яюЁ фър ьрыюёЄш т юёЄрЄърї). ═ю яЁюсыхьр ёюёЄюшЄ т Єюь, ўЄю
яЁш ¤Єюь ёрьш ёырурхь√х рёшьяЄюЄшўхёъюую Ё фр ёЄрэют Єё  эхюсючЁшь√ьш (ЇюЁьры№эю, фы  чряшёш тЄюЁюую ўыхэр рёшьяЄюЄшъш
ЄЁхсєхЄё  $\thicksim n^3$ ёырурхь√ї). ═р ёрьюь фхых, эх тёх ¤Єш ёырурхь√х шьх■Є юфшэръют√щ яюЁ фюъ. ╥ръ, т
\cite{SavSh03}, фы  ёыєўр  юяхЁрЄюЁр ╪ЄєЁьр--╦шєтшыы  $-y''+qy$, ъюуфр $n=2$, ртЄюЁ√ яюърчрыш, ўЄю юёэютэющ тъырф
тэюёшЄ ыш°№ юфшэ рёшьяЄюЄшўхёъшщ ўыхэ, ёыхфє■∙шх ЄЁш шьх■Є яюфўшэхээ√щ їрЁръЄхЁ, р юёЄрт°шхё  шьх■Є сюы№°шщ яюЁ фюъ ьрыюёЄш.
 ▌Єю, т ўрёЄэюёЄш, яючтюышыю ртЄюЁрь  \cite{SavSh10} фюърчрЄ№ ъы■ўхтє■ фы  Ёх°хэш  юсЁрЄэющ чрфрўш ЄхюЁхьє ю
 яЁюшчтюфэющ ╘Ёх°х %т эєых
 эхышэхщэюую юЄюсЁрцхэш  $q\mapsto\{\la_n\}_1^\infty$ (чфхё№ $\la_n$ --- ёюсёЄтхээ√х чэрўхэш ). ╚Єръ, трцэ√щ ¤ыхьхэЄ эютшчэ√
 ЁрсюЄ√ ёюёЄюшЄ т Єюь, ўЄю ь√ т√фхы хь т  тэюь тшфх юёэютэє■ ъюьяюэхэЄє юёЄрЄър т рёшьяЄюЄшўхёъшї ЇюЁьєырї ш яюърч√трхь,
 ўЄю юёЄрЄюъ ьюцэю ¤ЇЇхъЄштэю юЎхэштрЄ№ яЁш эрышўшш фюяюыэшЄхы№эющ шэЇюЁьрЎшш ю ьрЄЁшўэ√ї ¤ыхьхэЄрї ёшёЄхь√.

═р°ш яЁхфяюыюцхэш  ю ЇєэъЎшш $\rho$ эр эрёЄю ∙шщ ьюьхэЄ  ты ■Єё  эршсюыхх юс∙шьш. ┬ ъырёёшўхёъющ ЄхюЁшш ¤Єє ЇєэъЎш■ ёўшЄр■Є
рсёюы■Єэю эхяЁхЁ√тэющ (эрь эхшчтхёЄэ√ ЁрсюЄ√, уфх ЄЁхсютрышё№ сюыхх ёырс√х єёыютш ). ╥ръюх яЁхфяюыюцхэшх  ты хЄё 
хёЄхёЄтхээ√ь фы  ёшёЄхь, тючэшър■∙шї яЁш шчєўхэшш фшЇЇхЁхэЎшры№э√ї юяхЁрЄюЁют т√ёюъюую яюЁ фър. ─хыю т Єюь, ўЄю ¤Єш
юяхЁрЄюЁ√ ётюф Єё  ъ ёшёЄхьрь тшфр \eqref{eq:main}, р ъръ тшфэю шч фюърчрЄхы№ёЄтр ╧Ёхфыюцхэш  \ref{pr:gen},
фшруюэрышчрЎш  ьрЄЁшЎ√ $\V(x)$ тючьюцэр Єюы№ъю яЁш єёыютшш рсёюы■Єэющ эхяЁхЁ√тэюёЄш ¤ыхьхэЄют ¤Єющ ьрЄЁшЎ√.
 ╥хь эх ьхэхх, ёфхырээюх эрьш юёырсыхэшх ЄЁхсютрэшщ эр ЇєэъЎш■ $\rho(x)$ (Єюы№ъю  яюыюцшЄхы№эюёЄ№ ш ёєььшЁєхьюёЄ№), ёє∙хёЄтхээю яЁш шчєўхэшш ёшёЄхь
\eqref{eq:main} ъръ Єръют√ї (эряЁшьхЁ, яЁш шчєўхэшш ёшёЄхь√ Єшяр ─шЁрър).  ╟фхё№ цх юЄьхЄшь, ўЄю эр°х єёыютшх (i) эр ЇєэъЎш■ $\rho$
єцх  эхы№ч 
юёырсшЄ№ схч ёє∙хёЄтхээюую шчьхэхэш  шЄюуютюую Ёхчєы№ЄрЄр. ╬сэєыхэшх ЇєэъЎшш $\rho$ эр ьэюцхёЄтх яюыюцшЄхы№эющ ьхЁ√
тхфхЄ ъ яюЄхЁх уыртэюую ёырурхьюую т яЁртющ ўрёЄш \eqref{eq:main}. ╤ьхэр чэрър $\rho(x)$ (Єюўъш ёьхэ√ чэрър эрч√тр■Є
Єюўърьш яютюЁюЄр) Ёртэюёшы№эр ёьхэх чэрър ёяхъЄЁры№эюую ярЁрьхЄЁр, Є.х. яхЁхїюфє т фЁєующ ёхъЄюЁ ъюьяыхъёэющ яыюёъюёЄш.
╚чєўхэшх рёшьяЄюЄшъ  Ёх°хэшщ єЁртэхэшщ тЄюЁюую яюЁ фър ёю  чэръюяхЁхьхээющ ЇєэъЎшхщ $\rho(x)$ юсЁрчє■Є юЄфхы№эє■ ЄхьрЄшъє.
╤ ¤Єшьш шёёыхфютрэш ьш ьюцэю ючэръюьшЄ№ё  т ъэшурї ш ёЄрЄ№ ї  т \cite{Ben, Fleige, Zettl, BPT}, уфх шьх■Єё  фюяюыэшЄхы№э√х ёё√ыъш.

╩ЁрЄъю яЁюъюььхэЄшЁєхь  яЁшэ Є√х эрьш єёыютш  (i)--(iii).
╠√ ёючэрЄхы№эю юуЁрэшўштрхьё  тшфюь $\rho(x)\B$ уыртэюую яю яюЁ фъє ўыхэр т яЁртющ ўрёЄш \eqref{eq:main}, їюЄ  эр°
юёэютэющ Ёхчєы№ЄрЄ юс рёшьяЄюЄшъх ЇєэфрьхэЄры№эющ ьрЄЁшЎ√ $\Y(x,\la)$ хёЄхёЄтхээ√ь юсЁрчюь  ьюцэю яхЁхэхёЄш эр сюыхх юс∙шщ
ёыєўрщ $\B=\B(x)=\diag\{\varphi_1(x),\dots,\varphi_n(x)\}$ яЁш эхъюЄюЁ√ї фюяюыэшЄхы№э√ї єёыютш ї эр ЇєэъЎшш $\varphi_j$.
═ряЁшьхЁ, эр° ьхЄюф яючтюы хЄ яюыєўшЄ№  рэрыюушўэ√х Ёхчєы№ЄрЄ√, хёыш ЇєэъЎшш $\varphi_j$ яЁшэшьр■Є чэрўхэш  эр Ёрчышўэ√ї ыєўрї
$\{\gamma_j\}_1^n$  т ъюьяыхъёэющ яыюёъюёЄш. └ёшьяЄюЄшўхёъшх яЁхфёЄртыхэш  ьюцэю яюыєўшЄ№ ш яЁш сюыхх юс∙шї яЁхфяюыюцхэш ї эр ЇєэъЎшш
$\varphi_j$, эю юэш сєфєЄ ёяЁртхфышт√ эх т ёхъЄюЁрї, р т сюыхх ёыюцэ√ї юсырёЄ ї $\Omega_j\subset \mathbb C$,  ъюЄюЁ√х єцх ьюуєЄ эх яюъЁ√трЄ№ тё■
 ъюьяыхъёэє■ яыюёъюёЄ№.
 ╥хь эх ьхэхх, Єрър  яюёЄрэютър чрфрўш тяюыэх
ёюфхЁцрЄхы№эр. ╬эр тючэшърхЄ, эряЁшьхЁ, яЁш шчєўхэшш яюышэюьшры№э√ї юяхЁрЄюЁэ√ї яєўъют --- фшЇЇхЁхэЎшры№э√ї єЁртэхэшщ
т√ёюъюую яюЁ фър, ъю¤ЇЇшЎшхэЄ√ ъюЄюЁ√ї  ты ■Єё  ьэюуюўыхэрьш ёяхъЄЁры№эюую ярЁрьхЄЁр $\la$. ╧Ёш Ёх°хэшш эхъюЄюЁ√ї чрфрў
(эряЁшьхЁ, яЁш фюърчрЄхы№ёЄтх яюыэюЄ√ ёюсёЄтхээ√ї ш яЁшёюфшэхээ√ї ЇєэъЎшщ яєўъют) фюёЄрЄюўэю чэрЄ№ рёшьяЄюЄшъє Ёх°хэш 
$\Y(x,\la)$ эр эхёъюы№ъшї ыєўрї т $\bC$ ш юс∙є■ юЎхэъє эр ЁюёЄ $\Y(x,\la)$ яЁш $|\la|\to\infty$ --- фры№°х ЁрсюЄрхЄ
ЄхюЁхьр ╘Ёруьхэр--╦шэфхыхЇр (ёь. \cite{Sh1}, \cite{ShPliev}). ═ю т ¤Єющ ЁрсюЄх ь√ юёЄрты хь ¤Єш чрфрўш схч ЁрёёьюЄЁхэш ,
шэрўх ¤Єю єтхыю с√ эрё т ёЄюЁюэє Єхїэшўхёъшї єёыюцэхэшщ.

─рыхх,  ь√ фюяєёърхь ЁртхэёЄтр ёЁхфш ўшёхы $b_j$, Є.х. эх ЄЁхсєхь, ўЄюс√ $b_j$
с√ыш яюярЁэю Ёрчышўэ√ьш. ╘ръЄшўхёъш ¤Єю ючэрўрхЄ, ўЄю ёшёЄхьє \eqref{eq:main} єфюсэю ЁрёёьрЄЁштрЄ№ сыюърьш, юс·хфшэ   т
юфэє уЁєяяє єЁртэхэш  ё ёютярфр■∙шьш ўшёырьш $b_j$. ╥ръюх юсюс∙хэшх эх чрЄЁруштрхЄ яюърчрЄхыш ¤ъёяюэхэЄ√ т
рёшьяЄюЄшўхёъюь яЁхфёЄртыхэшш Ёх°хэш  $\Y(x,\la)$, эю ьхэ хЄ ьэюцшЄхыш яхЁхф ¤Єшьш ¤ъёяюэхэЄрьш, ъюЄюЁ√х  ты ■Єё  ЇєэъЎш ьш,
чртшё ∙шьш Єюы№ъю юЄ  $x$. ┬ ёыєўрх, ъюуфр ўшёыр $b_j$ яюярЁэю Ёрчышўэ√ ¤Єш ЇєэъЎшш т√яшё√тр■Єё   тэю, р т юс∙хщ ёшЄєрЎшш
юяЁхфхы ■Єё  ъръ Ёх°хэш  ёшёЄхь фшЇЇхЁхэЎшры№э√ї єЁртэхэшщ, эх чртшё ∙шї юЄ $\la$. ┼ёыш цх  фюяєёЄшЄ№ эрышўшх
цюЁфрэют√ї ъыхЄюъ т ёЄЁєъЄєЁх ьрЄЁшЎ√ $\B$, ь√ яюыєўшь фЁєує■ ёЄЁєъЄєЁє уыртэюую рёшьяЄюЄшўхёъюую ўыхэр --- єяюь эєЄ√х
т√°х ьэюцшЄхыш ёЄрэєЄ яюышэюьрьш ёяхъЄЁры№эюую ярЁрьхЄЁр. ▌ЄюЄ ёыєўрщ  Єхїэшўхёъш ёыюцхэ, тЁ ф ыш ЎхыхёююсЁрчэю яЁютюфшЄ№
 хую шчєўхэшх т юс∙хь тшфх, р эх фы  ъюэъЁхЄэющ ьюфхы№эющ чрфрўш.  ╙ёыютш 
ёєььшЁєхьюёЄш ъю¤ЇЇшЎшхэЄют ёшёЄхь√ \eqref{eq:main}  ты ■Єё  тяюыэх хёЄхёЄтхээ√ьш, їюЄ  эрь шчтхёЄхэ Ёхчєы№ЄрЄ
└ьшЁютр ш ├єёхщэютр \cite{AmGus}, ъюуфр яют√°хэшх уырфъюёЄш  ЇєэъЎшш $a_{12}(x)$  т ёшёЄхьх ─шЁрър  фрхЄ тючьюцэюёЄ№ ЁрёёьрЄЁштрЄ№ ёыєўрщ юсюс∙хээющ ЇєэъЎшш
$a_{21}(x)$. ╥Ёхсютрэшх єс√трэш  ЇєэъЎшщ $c_{jk}(x,\la)$ яЁш $|\la|\to\infty$  ты ■Єё , юўхтшфэю,
ьръёшьры№эю юс∙шьш (яЁш ЁрёёьюЄЁхэшш ёшёЄхь, тючэшър■∙шї яЁш шчєўхэшш юяхЁрЄюЁют ё ъю¤ЇЇшЎшхэЄрьш-ЁрёяЁхфхыхэш ьш ш
яЁш шчєўхэшш яюышэюьшры№э√ї яєўъют, шьххь $c_{jk}(x,\la)=O(|\la|^{-1})$). ╩юэхўэю, яют√°хэшх єёыютшщ эр уырфъюёЄ№
ЇєэъЎшщ $a_{ij}(x)$, юЄъЁ√тр■∙хх яєЄ№ ъ єыєў°хэш■ юЎхэюъ юёЄрЄъют т рёшьяЄюЄшъх ьрЄЁшЎ√ $\Y(x,\la)$, фюыцэю
ёюяЁютюцфрЄ№ё  єёшыхэшхь єёыютшщ эр єс√трэшх ЇєэъЎшщ $c_{jk}(x,\la)$.

┬ чртхЁ°хэшх --- ъЁрЄъю ю ёЄЁєъЄєЁх ЁрсюЄ√.  ┬ю тЄюЁюь ярЁуЁрЇх ь√ юяЁхфхы хь ёхъЄюЁ√, т ъюЄюЁ√ї сєфхь ЁрсюЄрЄ№ фы  яюыєўхэш  рёшьяЄюЄшъ
   ш яЁютюфшь яюфуюЄютшЄхы№эє■ ЁрсюЄє   фы  фюърчрЄхы№ёЄтр  т ЄЁхЄ№хь ярЁруЁрЇх яхЁтющ юёэютэющ ЄхюЁхь√ ю ЇєэфрьхэЄры№эющ ьрЄЁшЎх
   Ёх°хэшщ ёшёЄхь√ \eqref{eq:main}.  ┬ ўхЄтхЁЄюь ярЁруЁрЇх, шёяюы№чє   ьрЄЁшЎє Ёхуєы ЁшчрЎшш  ╠шЁчюхтр-╪ърышъютр,
     ь√ яЁютюфшь  ётхфхэшх єЁртэхэш  \eqref{1.1}  ъ ёшёЄхьх \eqref{eq:main} ш яюыєўрхь тЄюЁющ юёэютэющ Ёхчєы№ЄрЄ
      юс рёшьяЄюЄшърї ЇєэфрьхэЄры№эющ ёшёЄхь√ Ёх°хэшщ єЁртэхэш  \eqref{1.1}  ё ъю¤ЇЇшЎшхэЄрьш-ЁрёяЁхфхыхэш ьш,
      яюфўшэхээ√ьш єёыютш ь \eqref{eq:cond2}. ┬ ъюэЎх ЁрсюЄ√ ь√ юЄфхы№эю т√яшё√трхь яюыєўхээ√х юс∙шх Ёхчєы№ЄрЄ√  фы  эршсюыхх трцэюую ёыєўр  $n=2$. ╬ЄьхЄшь, ўЄю  Ёхчєы№ЄрЄ√ фрээющ
ЁрсюЄ√ ёюёЄрты ■Є срчє фы  шёёыхфютрэш  ёяхъЄЁры№э√ї ётющёЄт фшЇЇхЁхэЎшры№э√ї юяхЁрЄюЁют т√ёюъюую яюЁ фър ё
ъю¤ЇЇшЎшхэЄрьш--ЁрёяЁхфхыхэш ьш (яюф ёяхъЄЁры№э√ьш ётющёЄтрьш ь√ шьххь т тшфє рёшьяЄюЄшъє ёюсёЄтхээ√ї чэрўхэшщ ш
ёюсёЄтхээ√ї ЇєэъЎшщ, юЎхэъш  Ёхчюы№тхэЄ√ юяхЁрЄюЁют, ъюЄюЁ√х яюыєўр■Єё   яюёых фюяюыэшЄхы№эющ яюёЄрэютъш ъЁрхт√ї єёыютшщ, ЄхюЁхь√ ю схчєёыютэющ срчшёэюёЄш ёюсёЄтхээ√ї ш яЁшёюхфшэхээ√ї ЇєэъЎшщ ш фЁ.) ▌Єш шёёыхфютрэш  ртЄюЁ√ яырэшЁє■Є яЁютхёЄш т яюёыхфє■∙шї ЁрсюЄрї. ╬ЄьхЄшь, ўЄю яхЁтр  тхЁёш  ¤Єющ ЁрсюЄ√ с√ыр  яЁхфёЄртыхэр ртЄюЁрьш т яЁхяЁшэЄх \cite{SavSh17}.

\section{┬ёяюьюурЄхы№э√х Ёхчєы№ЄрЄ√}
\setcounter{equation}{0}

┬эрўрых юЄьхЄшь яюыхчэ√щ Ёхчєы№ЄрЄ ю тючьюцэюёЄш фшруюэрышчрЎшш ьрЄЁшЎ√ $\V(x)$ т ёшёЄхьх \eqref{eq:main}.

\begin{Proposition}\label{pr:gen} ╧єёЄ№ ьрЄЁшЎр $\V(x)$ ёшёЄхь√
\begin{equation}\label{eq:gen}
\u'=\la \V(x)\u+\A_0(x)\u+\C_0(x,\la)\u.
\end{equation}
фюяєёърхЄ яЁхфёЄртыхэшх $\V(x)=\W(x)(\rho(x)\B)\W^{-1}(x)$, уфх ьрЄЁшЎ√ $\W(x)$ ш $\W^{-1}(x)$ рсёюы■Єэю эхяЁхЁ√тэ√,
р ьрЄЁшЎр $\B$  ш ЇєэъЎш  $\rho(x)$ яюфўшэхэ√ єёыютш■ (i).\footnote{╥ръшь юсЁрчюь ьрЄЁшЎр $\V(x)$ яЁш ърцфюь $x\in[0,1]$ шьххЄ $n$ ёюсёЄтхээ√ї чэрўхэшщ, Ёртэ√ї
$\rho(x)b_j$, $1\le j\le n$. ╤Ёхфш ¤Єшї ёюсёЄтхээ√ї чэрўхэшщ фюяєёър■Єё  ъЁрЄэ√х, эю эх фюяєёър■Єё  цюЁфрэют√ ъыхЄъш,
Є.х. ь√ ЄЁхсєхь ёютярфхэш  ухюьхЄЁшўхёъющ ш рыухсЁршўхёъющ ъЁрЄэюёЄш фы  ърцфюую ёюсёЄтхээюую чэрўхэш  ьрЄЁшЎ√
$\V(x)$.}
─рыхх, яєёЄ№ ьрЄЁшЎ√ $\A_0(x)$ ш $\C_0(x,\la)$  яюфўшэхэ√ єёыютш ь (ii) ш (iii). ╥юуфр ёшёЄхьр
\eqref{eq:gen} ётюфшЄё  ъ тшфє \eqref{eq:main} чрьхэющ $\u(x)=\W(x)\y(x)$.
\end{Proposition}
\begin{proof}
╤фхырт чрьхэє $\u(x)=\W(x)\y(x)$, яюыєўшь
\begin{gather*}
\W'(x)\y(x)+\W(x)\y'(x)=\la \W(x)(\rho(x)\B)\y(x)+\A_0(x)\W(x)\y(x)+\C_0(x,\la)\W(x)\y(x)\\ \Longleftrightarrow\\
\y'=\la\rho(x)\B\y+\big(\W^{-1}(x)\A_0(x)\W(x)-\W^{-1}(x)\W'(x)\big)\y+\W^{-1}(x)\C_0(x,\la)\y.
\end{gather*}
╧юыюцшь $\A(x)=\W^{-1}(x)\A_0(x)\W(x)-\W^{-1}(x)\W'(x)$ ш $\C(x,\la)=\W^{-1}(x)\C_0(x,\la)$. ┬√яюыэхэшх єёыютшщ
(ii)--(iii) фы  ¤Єшї ьрЄЁшЎ юўхтшфэю.
\end{proof}

─рыхх эрь яюЄЁхсєхЄё   ЄхюЁхьр ёє∙хёЄтютрэш  ш хфшэёЄтхээюёЄш фы  ёшёЄхь√ $n\times n$ $\y'=\T(x)\y+\f(x)$ ё эрўры№э√ь
єёыютшхь $\y(\xi)=\y^0$ ш ёєььшЁєхь√ьш тхъЄюЁ--ЇєэъЎшхщ $\f(x)$ ш ьрЄЁшЎ--ЇєэъЎшхщ $\T(x)$. ▌Єр ЄхюЁхьр, ъюэхўэю,
їюЁю°ю шчтхёЄэр (ёь., эряЁшьхЁ, \cite{CL}). ═ю эрь  сєфєЄ эєцэ√, т яхЁтє■ юўхЁхф№, юЎхэъш
эр эюЁьє Ёх°хэш  ш эр эюЁьє ЁрчэюёЄш Ёх°хэшщ яЁш трЁ№шЁютрэшш $\xi$, $\f$ ш $\T$.

 ╧юыюцшь
\begin{gather*}
\|\y\|_{\AC}:=\sum_{j=1}^n\int_0^1|y_j(x)|+|y_j'(x)|\,dx,\qquad \|\y\|_\xi:=\sum_{j=1}^n|y_j(\xi)|+\int_0^1|y_j'(x)|\,dx,\\
\|\y\|_\C=\|\y\|_{L_\infty}:=\max_{\substack{1\le j\le
n\\x\in[0,1]}}|y_j(x)|,\qquad\|\y\|_{L_p}=\left(\sum_{j=1}^n\int_0^1|y_j(x)|^p\,dx\right)^{1/p}
\end{gather*}
ш ёЁрчє цх чрьхЄшь, ўЄю $\|\y\|_{\AC}\le2\|\y\|_\xi$ фы  ы■сюую $\xi\in[0,1]$ ш $\|\y\|_\C\le\|\y\|_{\AC}$. ─ы 
ЇшъёшЁютрээюую тхъЄюЁр $\y$ яюыюцшь $|\y|:=\max_{1\le j\le n}|y_j|$. ╬сючэрўшь ўхЁхч $|\f(x)|:=\max_{1\le j\le
n}|f_j(x)|$ фы  яЁюшчтюы№эющ тхъЄюЁ--ЇєэъЎшш $\f(x)$ ш ўхЁхч $|\T(x)|=\max_{1\le j\le n}\sum_{k=1}^n|t_{jk}(x)|$ фы 
яЁюшчтюы№эющ ьрЄЁшЎ--ЇєэъЎшш $\T(x)$. ─ы  ъЁрЄъюёЄш сєфхь яшёрЄ№ $\f\in L_1[0,1]$ шыш $\T\in L_1[0,1]$ т Єюь ёь√ёых,
ўЄю тёх $f_j(x)\in L_1[0,1]$ шыш $t_{jk}(x)\in L_1[0,1]$. ╬яЁхфхышь Єръцх ьрЄЁшўэ√х эюЁь√
$$
\|\T(x)\|_{L_p}=\sum_{jk=1}^{n}\left(\sum_{k=1}^n\int_0^1|t_{jk}(x)|^pdx\right)^{\frac1p},\qquad\|\T(x)\|_\C=\max_{1\le
j,\,k\le n}|t_{jk}(x)|,
$$
$$
\|\T(x)\|_{\AC}=\|\T(x)\|_{L_1}+\|\T'(x)\|_{L_1}.
$$
\begin{Proposition}\label{pr:1}
╧єёЄ№ ¤ыхьхэЄ√ ьрЄЁшЎ√ $\T$ ёшёЄхь√
\begin{equation}\label{eq:syst}
\z'=\T(x)\z+\f(x),\qquad x\in[0,1],\quad\z(\xi)=\z^0,\quad\f\in L_1[0,1],
\end{equation}
ЁрчьхЁр $n\times n$ яЁшэрфыхцрЄ яЁюёЄЁрэёЄтєя $L_1[0,1]$. ╥юуфр ёшёЄхьр \eqref{eq:syst} шьххЄ хфшэёЄтхээюх Ёх°хэшх $\z(x)\in \AC[0,1]$, яЁшўхь
\begin{equation}\label{eq:zgrowth}
\|\z\|_{\AC}\le(1+2\tau e^\tau)\|\g\|_{\AC},\quad\text{ш}\quad|\z(x)|\le e^{\tau(x)}\|\g\|_{\C},
\end{equation}
уфх
\begin{gather*}
\tau(x)=\left|\int_\xi^x|\T(t)|\,dt\right|,\qquad\tau=\tau(1),\qquad\text{р}\quad \g(x)=\z^0+\int_\xi^x\f(t)\,dt.
\end{gather*}
\end{Proposition}
\begin{proof}
╟ряш°хь ёшёЄхьє \eqref{eq:syst} т шэЄхуЁры№эюь тшфх
$$
\z(x)=\g(x)+(\T\z)(x),\quad\text{уфх}\ \ (\T\z)(x)=\int_\xi^x\T(t)\z(t)\,dt
$$
ш, яЁхцфх тёхую, фюърцхь юЎхэъє
\begin{equation}\label{eq:Tnpoint}
|(\T^l\z)(x)|\le\frac{\tau^l(x)}{l!}\|\z\|_C,\quad l\ge1.
\end{equation}
╦хуъю тшфхЄ№, ўЄю
$$
|(\T\z)(x)|\le\tau(x)\|\z\|_C.
$$
─рыхх яЁютхфхь фюърчрЄхы№ёЄтю яю шэфєъЎшш. ╧єёЄ№ юЎхэър єцх яюыєўхэр фы  $l-1$, Єюуфр
\begin{multline}
|(\T^l\z)(x)|\le\left|\int_\xi^x|\T(t)|\cdot|(\T^{l-1}\z)(t)|\,dt\right|\le\|\z\|_{\C}\left|\int_\xi^x|\T(t)|\frac{\tau^{l-1}(t)}{(l-1)!}\,dt\right|=\\
=\|\z\|_{\C}\left|\int_\xi^x\frac{\tau^{l-1}(t)}{(l-1)!}\,d\tau(t)\right|=\frac{\tau^l(x)}{l!}\|\z\|_{\C}
\end{multline}
ш юЎхэър \eqref{eq:Tnpoint} фюърчрэр. ╬Єё■фр ёыхфєхЄ фЁєур  юЎхэър
$$
\|\T^l\z\|_{\AC}\le2\|\T^l\z\|_\xi\le2\!\int_0^1\!\!|\T(x)|\cdot|(\T^{l-1}\z)(x)|\,dx\le2\|\z\|_{\C}\!\int_0^1\!\!|\T(x)|\frac{\tau^{l-1}(x)}{(l-1)!}\,dx
\le\frac{2\tau^l}{(l-1)!}\|\z\|_{\C}.
$$
┬ ўрёЄэюёЄш,
\begin{equation}\label{eq:Tnnorm}
\|\T^l\|_{\AC}\le\frac{2\tau^l}{(l-1)!}.
\end{equation}
╥хяхЁ№ Ёх°хэшх єЁртэхэш  \eqref{eq:syst} ь√ ьюцхь чряшёрЄ№ т тшфх Ё фр
\begin{equation}\label{eq:zseries}
\z(x)=\sum_{l=0}^\infty (\T^l\g)(x),
\end{equation}
яЁшўхь т ёшыє \eqref{eq:Tnpoint} ш \eqref{eq:Tnnorm} ¤ЄюЄ Ё ф ёїюфшЄё  ш т ърцфющ Єюўъх $x\in[0,1]$, ш яю эюЁьх
яЁюёЄЁрэёЄтр $\AC[0,1]$. ╧Ёш ¤Єюь
\begin{gather*}
|\z(x)|\le\sum_{l=0}^\infty|(\T^l\g)(x)|\le\|\g\|_\C\sum_{l=0}^\infty\frac{\tau^l(x)}{l!}\le\|\g\|_{\C}e^{\tau(x)},\\
\|\z\|_{\AC}\le\sum_{l=0}^\infty\|\T^l\g\|_{\AC}\le\|\g\|_{\AC}\left(1+2\sum_{l=1}^\infty\frac{\tau^l}{(l-1)!}\right)=\|\g\|_{\AC}(1+2\tau
e^\tau).
\end{gather*}
┼фшэёЄтхээюёЄ№ Ёх°хэш  фюърч√трхЄё  ёЄрэфрЁЄэю: хёыш $\z$ ш $\wt\z$ --- фтр Ёх°хэш  ёшёЄхь√ \eqref{eq:syst}, Єю
$$
\z-\wt\z=\T(\z-\wt\z)=\dots=\T^l(\z-\wt\z),
$$
юЄъєфр $\z-\wt\z=0$, яюёъюы№ъє юяхЁрЄюЁ $\T^l$  ты хЄё  ёцрЄшхь т яЁюёЄЁрэёЄтх $\AC[0,1]$ яЁш фюёЄрЄюўэю сюы№°юь $l$.
\end{proof}
\begin{Note}
╧Ёш ЇшъёшЁютрээ√ї $\xi$, $\z^0$ ш $\f\in L_1[0,1]$, Ёх°хэшх єЁртэхэш  \eqref{eq:syst} эхяЁхЁ√тэю чртшёшЄ юЄ ьрЄЁшЎ√
$\T\in L_1[0,1]$ т ёыхфє■∙хь ёь√ёых. ┼ёыш $\wt\z$ --- Ёх°хэшх эрўры№эющ чрфрўш \eqref{eq:syst} ё Єхьш цх $\xi$, $\z^0$
ш $\f$, эю ё фЁєующ ьрЄЁшЎхщ $\wt \T$, Єю
\begin{equation}
\|\z-\wt\z\|_{\AC}\le2\|\T-\wt \T\|_{L_1}(1+2\tau e^\tau)e^{\wt\tau}\|\g\|_{\AC}.
\end{equation}
\end{Note}
\begin{proof}
╧юыюцшь $\u(x)=\z(x)-\wt\z(x)$, ¤Єр ЇєэъЎш  рсёюы■Єэю эхяЁхЁ√тэр, Ёртэр эєы■ т Єюўъх $\xi$ ш єфютыхЄтюЁ хЄ єЁртэхэш■
$$
\u'=\T(x)\u+(\T(x)-\wt \T(x))\wt\z(x).
$$
╤юуырёэю \eqref{eq:zgrowth},
\begin{gather*}
\|\u\|_{\AC}\le(1+2\tau e^\tau)\left\|\int_\xi^x(\T(t)-\wt \T(t))\wt\z(t)\,dt\right\|_{\AC}\le2(1+2\tau
e^\tau)\left\|\int_\xi^x(\T(t)-\wt \T(t))\wt\z(t)\,dt\right\|_\xi=\\
=2(1+2\tau e^\tau)\int_0^1\left|(\T(x)-\wt \T(x))\wt\z(x)\right|\,dx\le2(1+2\tau e^\tau)\|\T-\wt \T\|_{L_1}\|\wt\z\|_\C\le\\
\le2\|\T-\wt \T\|_{L_1}(1+2\tau e^\tau)e^{\wt\tau}\|\g\|_{\AC}.
\end{gather*}
\end{proof}
\begin{Note}\label{note:2}
╧єёЄ№ ьрЄЁшЎр $\T$, ЇєэъЎш  $\f$ ш эрўры№эюх чэрўхэшх $\z^0$ чртшё Є юЄ ярЁрьхЄЁр $\la$, ьхэ ■∙хуюё  т юсырёЄш
$D\subset\bC$. ╧Ёхфяюыюцшь, ўЄю ЇєэъЎшш $z^0_j(\la)$ ш юЄюсЁрцхэш  $\la\mapsto f_j(\cdot,\la)$, $\la\mapsto
t_{jk}(\cdot,\la)$ шч $D$ т $L_1[0,1]$, $1\le j,\,k\le n$, уюыюьюЁЇэ√ яю $\la$. ╥юуфр ш Ёх°хэшх \eqref{eq:syst},
юЄюсЁрцхэшх $\la\mapsto \z(\cdot,\la)$ шч $D$ т $(\AC[0,1])^n$, уюыюьюЁЇэю яю $\la\in D$.
\end{Note}
\begin{proof}
├юыюьюЁЇэюёЄ№ ЇєэъЎшш $\g(\cdot,\la)$ т $(\AC[0,1])^n$ ёыхфєхЄ шч яЁхфёЄртыхэш  $\g(\cdot,\la)=\z^0(\la)+\int
\f(\cdot,\la)$. ├юыюьюЁЇэюёЄ№ ЇєэъЎшш $(\T\g)(\cdot,\la)$ ёыхфєхЄ шч юяЁхфхыхэш  юяхЁрЄюЁр $T$. ├юыюьюЁЇэюёЄ№ ЇєэъЎшщ
$(\T^m\g)(\cdot,\la)$ єёЄрэртыштрхЄё  яю шэфєъЎшш. ╧юёъюы№ъє ЇєэъЎш  $\tau(\la)$ эхяЁхЁ√тэр яю $\la\in D$, Єю юэр
юуЁрэшўхэр эр ърцфюь ъюьяръЄх. ╥юуфр юЎхэъш \eqref{eq:Tnnorm} урЁрэЄшЁє■Є ЁртэюьхЁэє■ эр ¤Єюь ъюьяръЄх ёїюфшьюёЄ№ Ё фр
\eqref{eq:zseries}. ╥хяхЁ№ уюыюьюЁЇэюёЄ№ ЇєэъЎшш $\z(\cdot,\la)$ ёыхфєхЄ шч ЄхюЁхь√ ┬хщхЁ°ЄЁрёёр\footnote{─ы  єфюсёЄтр ўшЄрЄхы  ь√ фрхь чфхё№
ъЁрЄъюх фюърчрЄхы№ёЄтю, яюфЁюсэюёЄш ёь., эряЁшьхЁ, т \cite[├ыртр 12]{BogSmol}.}.
\end{proof}
\begin{Corollary}
╧Ёш ърцфюь ЇшъёшЁютрээюь $\la\in\bC$, $|\la|>\la_0$, ёшёЄхьр \eqref{eq:main} яюёых фюсртыхэш  эрўры№эюую єёыютш 
$\y(\xi,\la)=\y^0(\la)$ шьххЄ хфшэёЄтхээюх Ёх°хэшх т ъырёёх рсёюы■Єэю эхяЁхЁ√тэ√ї эр $[0,1]$ ЇєэъЎшщ. ╧Ёш ¤Єюь
$\AC$--эюЁьр ¤Єюую Ёх°хэш  фюяєёърхЄ юЎхэъє
$$
\|\y\|_{\AC}\le (1+2\tau e^{\tau})|\y^0(\la)|,\qquad\text{уфх}\ \tau=|\la|\|\rho\|_{L_1}\max_{1\le j\le
n}|b_j|+a+c,\quad |\y^0(\la)|=\max_{1\le j\le n}|y^0_j(\la)|,
$$
яЁш тёхї $ |\la|>\max\{a+c,\la_0\}$. ╧єёЄ№ тёх ЇєэъЎшш $y^0_j(\la)$, $1\le j\le n$, уюыюьюЁЇэ√ яю
$\la\in\{|\la|>\la_0\}$, Єръ цх, ъръ ш юЄюсЁрцхэш  $\la\mapsto c_{jk}(\cdot,\la)$, $1\le j,\,k\le n$, шч юсырёЄш
$\la\in\{|\la|>\la_0\}$ т $L_1[0,1]$. ╥юуфр уюыюьюЁЇэю ш Ёх°хэшх, юЄюсЁрцхэшх $\la\mapsto\y(\cdot,\la)$ шч
$\{|\la|>\la_0\}$ т $(\AC[0,1])^n$.
\end{Corollary}
╬яЁхфхышь ёхьхщёЄтю ёхъЄюЁют $\Gamma_\kappa=\{\la:\arg\la\in(\al_{\kappa-1},\al_{\kappa})\}$, $\kappa=1,\dots,J$. ┼ёыш
тёх ўшёыр $b_j$ ёютярфр■Є, Єю эр°х ёхьхщёЄтю сєфхЄ ёюёЄю Є№ Ёютэю шч юфэюую ёхъЄюЁр $\Gamma_1=\bC$. ┬ яЁюЄштэюь ёыєўрх
чрЇшъёшЁєхь фтр яЁюшчтюы№э√ї шэфхъёр $1\le k<l\le n$ Єръшї, ўЄю $b_k\ne b_l$, ш ЁрёёьюЄЁшь єЁртэхэшх
\begin{equation}\label{eq:sect1}
\Re(b_k\la)=\Re(b_l\la)\ \Longleftrightarrow\ \Re((b_k-b_l)\la)=0.
\end{equation}
╦хуъю тшфхЄ№, ўЄю Ёх°хэшхь ¤Єюую єЁртэхэш   ты хЄё  эхъюЄюЁр  яЁ ьр , яЁюїюф ∙р  ўхЁхч эрўрыю ъююЁфшэрЄ. ╬с∙хх ўшёыю
єЁртэхэшщ тшфр \eqref{eq:sect1} Ёртэю $n(n-1)/2$, Єръ ўЄю т Ёхчєы№ЄрЄх яюыєўрхь Ёрчсшхэшх ъюьяыхъёэющ яыюёъюёЄш эр
$1\le J\le(n^2-n)$ ёхъЄюЁют. ╧Ёш $J>1$ ърцф√щ ёхъЄюЁ юуЁрэшўхэ фтєь  яюы Ёэ√ьш ыєўрьш, ъюЄюЁ√х ь√ чрэєьхЁєхь яю
тючЁрёЄрэш■ рЁуєьхэЄр $\al_0\le0<\al_1<\al_2<\dots<\al_{J-1}<\al_J=\al_0+2\pi$.

╟рЇшъёшЁєхь эхъюЄюЁ√щ эюьхЁ $\kappa\le J$, ЁрёёьюЄЁшь ёхъЄюЁ $\Gamma_\kappa$ ш эрщфхь т ¤Єюь ёхъЄюЁх Ёх°хэшх єЁртэхэш 
\eqref{eq:matrmain}, шьх■∙хх юяЁхфхыхээюх рёшьяЄюЄшўхёъюх яЁхфёЄртыхэшх. ┴юыхх Єюую, ь√ яюърцхь, ўЄю Ёх°хэшх ё Єръшь
рёшьяЄюЄшўхёъшь ётющёЄтюь юяЁхфхыхэю эх Єюы№ъю т $\Gamma_\kappa$, эю ш т сюыхх °шЁюъюь ёхъЄюЁх $\wt \Gamma_\kappa$. ╠√
юяЁхфхышь ёхъЄюЁ $\wt \Gamma_\kappa$ ъръ Ёхчєы№ЄрЄ ярЁрыыхы№эюую яхЁхэюёр ёхъЄюЁр $\Gamma_\kappa$ тфюы№ хую сшёёхъЄЁшё√
(Єюўэхх, тфюы№ яЁюфюыцхэш  ¤Єющ сшёёхъЄЁшё√ чр яЁхфхы√ $\Gamma_\kappa$):
$$
\wt \Gamma_\kappa=\wt \Gamma_\kappa(r):=\{\la\in\bC:\la+re^{\frac{i}2(\al_{\kappa-1}+\al_\kappa)}\in \Gamma_\kappa\},
$$
уфх $r>0$ ЇшъёшЁютрэю. ┬ ёыєўрх, ъюуфр $J=1$ ш $\Gamma_1=\bC$ яюыюцшь $\wt\Gamma_1=\bC$.

╟рьхЄшь, ўЄю фы  ы■сющ ярЁ√ шэфхъёют $1\le k,\,l\le n$ Єръшї, ўЄю $b_k\ne b_l$, т ёхъЄюЁх $\Gamma_\kappa$ т√яюыэхэю
ёыхфє■∙хх ётющёЄтю: ышсю $\Re(b_k\la)>\Re(b_l\la)$ фы  тёхї $\la\in \Gamma_\kappa$, ышсю $\Re(b_k\la)<\Re(b_l\la)$ фы 
тёхї $\la\in \Gamma_\kappa$. ╧єёЄ№ $\nu$ --- ўшёыю Ёрчышўэ√ї ўшёхы $b_j$. ╧хЁхэєьхЁєхь єЁртэхэш  ёшёЄхь√
\eqref{eq:main} Єръ, ўЄю
\begin{multline}\label{eq:ordnung}
\Re(b_1\la)=\Re(b_2\la)=\dots=\Re(b_{n_1}\la)>\Re(b_{n_1+1}\la)=\dots=\Re(b_{n_1+n_2}\la)>\dots\\
\ldots>\Re(b_{n_1+\dots+n_{\nu-1}+1}\la)=\dots=\Re(b_n\la)\quad\forall\la\in \Gamma_\kappa,\qquad
n=n_1+n_2+\dots+n_\nu.
\end{multline}
╙фюсэю сєфхЄ юсючэрўшЄ№ $b_1=\dots=b_{n_1}=:\beta_1$, $b_{n_1+1}=\dots=b_{n_1+n_2}=:\beta_2$ ш Є.ф., яхЁхщф  Єръшь
юсЁрчюь, ъ эрсюЁє $\beta_1,\dots,\beta_\nu$ яюярЁэю Ёрчышўэ√ї ўшёхы. ┬ ёхъЄюЁх $\wt \Gamma_\kappa$ эхЁртхэёЄтр
\eqref{eq:ordnung}, хёЄхёЄтхээю эх т√яюыэхэ√ (чр шёъы■ўхэшхь ёыєўр  $\wt\Gamma_1=\Gamma_1=\bC$).
\begin{Lemma}\label{lem:shift}
═рщфєЄё  Єръшх ўшёыр $h>0$ ш $\la_0>0$ (юэш чртшё Є юЄ тхышўшэ√ ёфтшур $r$), ўЄю т юсырёЄш
$\Dom_{\kappa,\la_0}:=\{\la\in\wt \Gamma_\kappa,\,|\la|>\la_0\}$ ёяЁртхфышт√ эхЁртхэёЄтр
\begin{equation}\label{eq:d}
\Re((b_k-b_l)\la)>-h, \quad 1\le k<l\le n;\qquad \Re((\beta_k-\beta_{l})\la)>-h, \quad 1\le k<l\le \nu.
\end{equation}
\end{Lemma}
\begin{proof}
╬ўхтшфэю, фюёЄрЄюўэю фюърчрЄ№ эхЁртхэёЄтр фы  ўшёхы $b_j$. ╧Ёхцфх тёхую, чрьхЄшь, ўЄю ёхъЄюЁ $\wt \Gamma_\kappa$
юуЁрэшўхэ ыєўюь $\ell_1$ тшфр $\la=-re^{\frac{\i}2(\al_{\kappa-1}+\al_\kappa)}+te^{\i\al_{\kappa-1}}$, $t>0$,
ярЁрыыхы№э√ь ыєўє $\arg\la=\al_{\kappa-1}$, ш ыєўюь $\ell_2$ тшфр
$\la=-re^{\frac{\i}2(\al_{\kappa-1}+\al_\kappa)}+te^{\i\al_{\kappa}}$, $t>0$, ярЁрыыхы№э√ь ыєўє $\arg\la=\al_\kappa$.
╬ўхтшфэю, ўЄю ыєў $\ell$ тшфр $\la=a+te^{\i\varphi}$, $t>0$, эрўшэр  ё эхъюЄюЁюую $t$, эх яюярфрхЄ т ёхъЄюЁ
$\Gamma=\{\la:\arg\la\in(\varphi_0,\varphi_1)\}$, хёыш $\varphi\notin[\varphi_0,\varphi_1]$. ╥ръшь юсЁрчюь, ыєў
$\ell_1$ ыхцшЄ, эрўшэр  ё эхъюЄюЁюую $t=t_1$, ышсю т ёхъЄюЁх $\Gamma_\kappa$, ышсю т ёюёхфэхь ёхъЄюЁх
$\Gamma_{\kappa-1}$ (ь√ ёўшЄрхь, ўЄю $\Gamma_0:=\Gamma_J$, р $\Gamma_{J+1}:=\Gamma_1$). ╧хЁтюх эхтючьюцэю, яюёъюы№ъє яю
юяЁхфхыхэш■ $\wt \Gamma_\kappa\supset \Gamma_\kappa$. └эрыюушўэю, ыєў $\ell_2$ ыхцшЄ, эрўшэр  ё эхъюЄюЁюую $t=t_2$, т
ёхъЄюЁх $\Gamma_{\kappa+1}$. ┬√сЁрт ўшёыю $\la_0$ Ёртэ√ь ьръёшьєьє шч $|\ell_1(t_1)|$ ш $|\ell_2(t_2)|$, яюыєўшь, ўЄю
юсырёЄ№ $D_{\kappa,\la_0}$ тыюцхэр т юс·хфшэхэшх ёхъЄюЁют $\Gamma_{\kappa-1}\cup\overline{\Gamma_\kappa}\cup
\Gamma_{\kappa+1}$.

╟рЇшъёшЁєхь ЄхяхЁ№ яЁюшчтюы№эє■ ярЁє шэфхъёют $1\le k<l\le n$ ш эрщфхь Єръюх ўшёыю $h_{kl}$, ўЄю эхЁртхэёЄтю
$\Re((b_k-b_l)\la)>-h_{kl}$ т√яюыэхэю тё■фє т $D_{\kappa,\la_0}$. ─ы  ¤Єюую ЁрёёьюЄЁшь $\bR$--ышэхщэє■ ЇєэъЎш■
$w_{kl}(\la)=\Re((b_k-b_l)\la)$. ┴єфхь ёўшЄрЄ№, ўЄю $b_k\ne b_l$, шэрўх ь√ яюыюцшь $h_{kl}=0$. ┬ ёшыє
\eqref{eq:ordnung} ¤Єр ЇєэъЎш  яюыюцшЄхы№эр яЁш $\la\in \Gamma_\kappa$. ┬ючьюцэ√ ўхЄ√Ёх ёыєўр : ¤Єр ЇєэъЎш  ьюцхЄ
юърчрЄ№ё  яюыюцшЄхы№эющ ш т $\Gamma_{\kappa-1}$, ш т $\Gamma_{\kappa+1}$. ┬ ¤Єюь ёыєўрх ь√ яюыюцшь $h_{kl}=0$ ш яЁшфхь
ъ эхЁртхэёЄтє $w_{kl}(\la)>-h_{kl}$ яЁш $\la\in D_{\kappa,\la_0}$. ╘єэъЎш  $w_{kl}$ ьюцхЄ юърчрЄ№ё  яюыюцшЄхы№эющ т
$\Gamma_{\kappa-1}$ ш т $\Gamma_\kappa$, эю Ёртэющ эєы■ эр ыєўх $\arg\la=\al_{\kappa}$. ╥юуфр $w_{kl}$ юЄЁшЎрЄхы№эр т
$\Gamma_{\kappa+1}$, р Єръ ъръ ыєў $\ell_2$ ярЁрыыхыхэ ыєўє $\arg\la=\al_\kappa$, Єю $w_{kl}$ яюёЄю ээр эр $\ell_2$. ┬
¤Єюь ёыєўрх ь√ яюыюцшь $h_{kl}=-w_{kl}\vert_{\ell_2}$. ╥юуфр эр ърцфюь ярЁрыыхы№эюь $\ell_2$ ыєўх, ыхцр∙хь т яюыюёх
ьхцфє ыєўрьш $\ell_2$ ш $\arg\la=\al_\kappa$, ЇєэъЎш  $w_{kl}$ Єръцх яюёЄю ээр ш яЁшэшьрхЄ чэрўхэш  т яЁюьхцєЄъх
$(-h_{kl},0)$. ╬Єё■фр ёыхфєхЄ эхЁртхэёЄтю $w_{kl}(\la)>-h_{kl}$ т $D_{\kappa,\la_0}$. ╘єэъЎш  $w_{kl}$ ьюцхЄ юърчрЄ№ё 
яюыюцшЄхы№эющ т $\Gamma_\kappa$ ш т $\Gamma_{\kappa+1}$, эю юЄЁшЎрЄхы№эющ т $\Gamma_{\kappa-1}$
--- ¤Єр ёшЄєрЎш  рэрыюушўэр яЁхф√фє∙хщ, ь√ яюыюцшь $h_{kl}=-w_{kl}\vert_{\ell_1}$ ш тэют№ яЁшфхь ъ эхЁртхэёЄтє
$w_{kl}(\la)>-h_{kl}$ тё■фє т $D_{\kappa,\la_0}$. ═ръюэхЎ, тючьюцхэ ёыєўрщ, ъюуфр $w_{kl}$ яюыюцшЄхы№эр Єюы№ъю т
ёхъЄюЁх $\Gamma_\kappa$, р т ёхъЄюЁрї $\Gamma_{\kappa-1}$ ш $\Gamma_{\kappa+1}$ юЄЁшЎрЄхы№эр. ╥юуфр $w_{kl}=0$ ш эр
ыєўх $\arg\la=\al_\kappa$, ш эр ыєўх $\arg\la=\al_{\kappa-1}$. ┬ ёшыє ышэхщэюёЄш $w_{kl}$ ¤Єю тючьюцэю Єюы№ъю Єюуфр,
ъюуфр ¤Єш ыєўш ярЁрыыхы№э√. ╧юёъюы№ъє $\al_\kappa\ne\al_{\kappa-1}$, Єю $\al_\kappa=\al_{\kappa-1}+\pi$. ┬ ¤Єюь ёыєўрх
ёхъЄюЁ $\Gamma_{\kappa}$ яЁхтЁр∙рхЄё  т яюыєяыюёъюёЄ№, р чэрўшЄ ыєўш $\ell_1$ ш $\ell_2$ ёюёЄрты ■Є юфэє яЁ ьє■.
╧юыюцшь $h_{kl}=-w_{kl}\vert_{\ell_1}=-w_{kl}\vert_{\ell_2}$ ш тэют№ яюыєўшь $w_{kl}(\la)>-h_{kl}$ т
$D_{\kappa,\la_0}$. ╬ёЄрхЄё  юсючэрўшЄ№ $h=\max_{1\le k<l\le n}h_{kl}$.
\end{proof}
╟рЇшъёшЁєхь эхъюЄюЁ√щ ёхъЄюЁ $\wt\Gamma_\kappa$ ш яхЁхэєьхЁєхь єЁртэхэш  ёшёЄхь√ \eqref{eq:main} Єръ, ўЄюс√ т
$\Gamma_\kappa$ с√ыш т√яюыэхэ√ эхЁртхэёЄтр \eqref{eq:ordnung}. ┴єфхь фрыхх ёўшЄрЄ№, ўЄю $\la\in\Dom_{\kappa,\la_0}$.
─ы  ъЁрЄъюёЄш сєфхь юсючэрўрЄ№ ¤Єє юсырёЄ№ $\Dom$. ╫хЁхч $h$ сєфхь юсючэрўрЄ№ ўшёыю, юяЁхфхыхээюх т ыхььх
\ref{lem:shift}. ╥ръшь юсЁрчюь, фы  ы■сюую $\la\in\Dom$ ёяЁртхфышт√ эхЁртхэёЄтр \eqref{eq:d}.

═рўэхь ё яюшёър уыртэюую ўыхэр т рёшьяЄюЄшўхёъюь яЁхфёЄртыхэшш фы  $\Y(x,\la)$ яЁш $\Dom\ni\la\to\infty$. ─ы  ¤Єюую ь√
юЄсЁюёшь т яЁртющ ўрёЄш \eqref{eq:matrmain} ёырурхьюх $\C(x,\la)\Y(x,\la)$. ─рыхх, ьрЄЁшЎє $\A(x)$ яЁхфёЄртшь т тшфх
$$
\A(x)=\D(x)+(\A(x)-\D(x)),\quad \text{уфх}\ \  \D(x)=\begin{cases}a_{jk}(x),\ \ &\text{хёыш}\
b_j=b_k,\\0,&\text{шэрўх},
\end{cases}
$$
ш ёырурхьюх $(\A(x)-\D(x))\Y(x,\la)$ т яЁртющ ўрёЄш \eqref{eq:matrmain} Єръцх юЄсЁюёшь. ┬ Ёхчєы№ЄрЄх ь√ яЁшфхь ъ
ьрЄЁшўэюьє фшЇЇхЁхэЎшры№эюьє єЁртэхэш■
\begin{gather}\label{eq:matrmajor}
\Y^0(x,\la)'=\la\rho(x)\B\Y^0(x,\la)+\D(x)\Y^0(x,\la),\\
\B=\begin{pmatrix}b_1&0&0&\ldots&0\\0&b_2&0&\ldots&0\\\hdotsfor{5}\\0&\ldots&0&0&b_n\end{pmatrix},\qquad
\D(x)=\begin{pmatrix}\D_1(x)&\O&\O&\hdotsfor{2}&\O\\\O&\D_2(x)&\O&\hdotsfor{2}&\O\\
\hdotsfor{6}\\\O&\hdotsfor{2}&\O&\O&\D_\nu(x)\end{pmatrix}\notag
\end{gather}
(ьрЄЁшЎ√ $\D_j(x)$ шьх■Є ЁрчьхЁ $n_j\times n_j$), ъюЄюЁюх ь√ сєфхь эрч√трЄ№ \textit{ьюфхы№э√ь єЁртэхэшхь} фы  юёэютэюую
єЁртэхэш  \eqref{eq:matrmain}.
\begin{Proposition}\label{pr:model}
╨х°хэшх $\Y^0(x,\la)$ єЁртэхэш  \eqref{eq:matrmajor} ё эрўры№э√ь єёыютшхь $\Y^0(0,\la)=\I$ шьххЄ тшф
\begin{equation}\label{eq:Y0}
\Y^0(x,\la)=\cE(x,\la)\cdot\M(x),\qquad\cE(x,\la)=\diag\{e^{b_1\la p(x)},\dots,e^{b_n\la p(x)}\},
\end{equation}
уфх $p(x)=\int_0^x\rho(t)\,dt$, ьрЄЁшЎ√ $\M_j(x)$, $1\le j\le \nu$, шьх■Є ЁрчьхЁ $n_j\times n_j$ ърцфр , р
$$
\M(x)=\begin{pmatrix}\M_1(x)&\O&\O&\hdotsfor{2}&\O\\\O&\M_2(x)&\O&\hdotsfor{2}&\O\\
\hdotsfor{6}\\\O&\hdotsfor{2}&\O&\O&\M_\nu(x)\end{pmatrix}.
$$
╩рцфр  шч ьрЄЁшЎ $\M_j(x)$  ты хЄё  Ёх°хэшхь єЁртэхэш  $\M_j'(x)=\D_j(x)\M_j(x)$ ё эрўры№э√ь єёыютшхь $\M_j(0)=\I$.
\end{Proposition}
\begin{proof}
╧ю юяЁхфхыхэш■, $\D(x)$ --- сыюъ--фшруюэры№эр  ўрёЄ№ ьрЄЁшЎ√ $\A(x)$, ёюёЄю ∙р  шч $\nu$ ътрфЁрЄэ√ї сыюъют ЁрчьхЁр
$n_j\times n_j$, $1\le j\le \nu$, ърцф√щ. ╧Ёш ¤Єюь т ёшёЄхьх \eqref{eq:Y0} ърцфюьє сыюъє ьрЄЁшЎ√ $\D$ юЄтхўрхЄ юфэю ш
Єю цх ўшёыю $b_j$, Ёрёяюыюцхээюх эр фшруюэрыш ьрЄЁшЎ√ $\B$. ╥ръшь юсЁрчюь, ёшёЄхьр \eqref{eq:matrmajor} ЁрёярфрхЄё  эр
$\nu$ эх ёт чрээ√ї ьхцфє ёюсющ ёшёЄхь
$$
\Y_j'(x,\la)=\la\beta_j\rho(x)\Y_j(x,\la)+\D_j(x)\Y_j(x,\la),\quad \Y_j(0,\la)=\I.
$$
╧Ё ьющ яюфёЄрэютъющ єсхцфрхьё , ўЄю $\Y_j(x,\la)=\exp\{\la\beta_jp(x)\}\M_j(x)$.
\end{proof}
\begin{Note}\label{nt:diagM}
┬ юс∙хь ёыєўрх ьрЄЁшЎє $\M(x)$ эхы№ч  т√ЁрчшЄ№  тэю ўхЁхч ¤ыхьхэЄ√ ьрЄЁшЎ√ $\A(x)$, юфэръю хёыш ьрЄЁшЎр $\D_j(x)$
 ты хЄё  фшруюэры№эющ, Є.х. $\D_j(x)=\diag\{a_{k_j,k_j}(x),\dots,\,a_{k_j+n_j,k_j+n_j}(x)\}$, уфх
$k_j=n_1+\dots+n_{j-1}+1$, Єю $\M_j(x)$ Єръцх фшруюэры№эр ш Ёртэр
$$
\M_j(x)=\diag\left\{\exp\left\{\int_0^x a_{k_j,k_j}(t)\,dt\right\},\ldots,\,\exp\left\{\int_0^x
a_{k_j+n_j,k_j+n_j}(t)\,dt\right\}\right\}.
$$
┬ ўрёЄэюёЄш, ¤Єю Єръ, хёыш тёх $b_j$, $1\le j\le n$, яюярЁэю Ёрчышўэ√.
\end{Note}
\begin{Note}
╠рЄЁшЎр $M(x)$ ъюььєЄшЁєхЄ ё ьрЄЁшЎрьш $B$ ш $\cE(x,\la)$.
%┬ ЁртхэёЄтю \eqref{eq:Y0} ьрЄЁшЎ√ ьюцэю ьхэ Є№ ьхёЄрьш
%\begin{multline*}
%\Y^0(x,\la)=\begin{pmatrix}e^{\la\beta_1p(x)}\M_1(x)&\O&\O&\hdotsfor{2}&\O\\\O&e^{\la\beta_2p(x)}\M_2(x)&\O&\hdotsfor{2}&\O\\
%\hdotsfor{6}\\\O&\hdotsfor{2}&\O&\O&e^{\la\beta_pp(x)}\M_p(x)\end{pmatrix}=\\=\M(x)\cdot\diag\{e^{b_1\la
%p(x)},\dots,e^{b_n\la p(x)}\}.
%\end{multline*}
\end{Note}
═ряюьэшь, ўЄю $a=\|\A(x)\|_{L_1}$, $\gamma(\la)=\|C(\cdot,\la)\|_{L_1}$, $c=\max_{|\la|>\la_0}\gamma(\la)$.
\begin{Lemma}
╠рЄЁшЎр $\M(x)$ фюяєёърхЄ юЎхэъє $\|\M(x)\|_\C\le e^{a}$. ─ы  ы■сюую $x\in[0,1]$ юэр юсЁрЄшьр, юсЁрЄэр  ьрЄЁшЎр
єфютыхЄтюЁ хЄ єЁртэхэш■ $(\M^{-1}(x))'=-\M^{-1}(x)\D(x)$ ё эрўры№э√ь єёыютшхь $\M(0)=\I$ ш Єръ цх фюяєёърхЄ юЎхэъє
$\|\M^{-1}(x)\|_\C\le e^{a}$.
\end{Lemma}
\begin{proof}
╬Ўхэър фы  ьрЄЁшЎ√ $\M$ ёыхфєхЄ шч \eqref{eq:zgrowth}. ═хт√ЁюцфхээюёЄ№ ьрЄЁшЎ√ $\M(x)$ фюърч√трхЄё  юс√ўэ√ь юсЁрчюь: хх
юяЁхфхышЄхы№ єфютыхЄтюЁ хЄ єЁртэхэш■ $(\det\M(x))'=\tr D(x)\cdot\det\M(x)$ ш эрўры№эюьє єёыютш■ $\det\M(0)=1$, р чэрўшЄ
юЄышўхэ юЄ эєы  эр тёхь юЄЁхчъх $[0,1]$. ─шЇЇхЁхэЎшЁє  ЄюцфхёЄтю $\I=\M^{-1}(x)\M(x)$, єўшЄ√тр  $\M'(x)=\D(x)\M(x)$ ш
фюьэюцр  ёяЁртр эр $\M^{-1}(x)$, яюыєўшь $(\M^{-1}(x))'=-\M^{-1}(x)\D(x)$. ╧хЁхїюф  ъ ёюяЁ цхээ√ь ьрЄЁшЎрь, ёюуырёэю
\eqref{eq:zgrowth}, яюыєўшь юЎхэъє $\left\|\left(\M^{-1}\right)^*\right\|_\C\le e^{a}$, ¤ътштрыхэЄэє■ шёъюьющ.
\end{proof}
╬яЁхфхышь ьрЄЁшЎ√
$$
\Q(x)=\M^{-1}(x)(\A(x)-\D(x))\M(x),\qquad \R(x,\la)=\M^{-1}(x)\C(x,\la)\M(x),
$$
¤ыхьхэЄ√ ъюЄюЁ√ї ь√ сєфхь юсючэрўрЄ№ $q_{jk}(x)$ ш $r_{jk}(x,\la)$. ╤Ёрчє цх чрьхЄшь, ўЄю
\begin{equation}\label{eq:qest}
\|\Q(x)\|_{L_1}\le ae^{2a}\qquad \text{ш} \qquad \|\R(x,\la)\|_{L_1}\le \gamma(\la)e^{2a}.
\end{equation}
\begin{Lemma}\label{lem:3}
┼ёыш фы  ъръющ--ышсю ярЁ√ шэфхъёют $b_k=b_l$, Єю $q_{kl}(x)=q_{lk}(x)=0$.
\end{Lemma}
\begin{proof}
╧Ёхцфх тёхую, чрьхЄшь, ўЄю ьрЄЁшЎр $\M^{-1}(x)$ шьххЄ Єръє■ цх ёЄЁєъЄєЁє, ўЄю ш $\M(x)$, яЁшўхь Ёютэю ё Єръшьш цх яю
ЁрчьхЁє ш яюыюцхэш■ сыюърьш, Ёртэ√ьш $\M_j^{-1}(x)$. ╥юуфр, Ёрчсштр  ьрЄЁшЎє $\A(x)$ эр ёююЄтхЄёЄтє■∙шх сыюъш, шьххь
\begin{gather*}
\Q(x)=\M^{-1}(x)(\A(x)-\D(x))\M(x)=\begin{pmatrix}\M^{-1}_1(x)&\O&\O&\hdotsfor{2}&\O\\\O&\M^{-1}_2(x)&\O&\hdotsfor{2}&\O\\
\hdotsfor{6}\\\O&\hdotsfor{2}&\O&\O&\M^{-1}_\nu(x)\end{pmatrix}\times\\
\times\begin{pmatrix}\O&\A_{12}(x)&\A_{13}(x)&\hdotsfor{1}&\A_{1\nu}(x)\\
\A_{21}(x)&\O&\A_{23}(x)&\hdotsfor{1}&\A_{2\nu}(x)\\
\hdotsfor{5}\\\A_{\nu1}(x)&\hdotsfor{2}&\A_{\nu\,\nu-1}(x)&\O\end{pmatrix}\times\begin{pmatrix}\M_1(x)&\O&\O&\hdotsfor{2}&\O\\\O&\M_2(x)&\O&\hdotsfor{2}&\O\\
\hdotsfor{6}\\\O&\hdotsfor{2}&\O&\O&\M_\nu(x)\end{pmatrix},
\end{gather*}
ўЄю ш тыхўхЄ яюёых яхЁхьэюцхэш  ьрЄЁшЎ єЄтхЁцфхэшх ыхьь√.
\end{proof}
╧юыюцшь $\rv_{jl}(x,\la)=q_{jl}(x)+r_{jl}(x,\la)$. ╤шьтюыюь $(\pm)_{jk}$ сєфхь юсючэрўрЄ№ тхышўшэє, Ёртэє■ $-1$ яЁш
$j<k$ ш $1$ яЁш $j\ge k$. ─ы  чряшёш юёЄрЄъют юяЁхфхышь ЇєэъЎшш
\begin{gather}\label{eq:defups}
\upsilon_{jkl}(s,x,\la)=(\pm)_{jk}(\pm)_{lk}\int
q_{jl}(t)e^{(b_l-b_k)\la(p(t)-p(s))+(b_j-b_k)\la(p(x)-p(t))}\,dt,\\\notag
\varrho_{jkl}(s,x,\la)=(\pm)_{jk}(\pm)_{lk}\int r_{jl}(t,\la)e^{(b_l-b_k)\la(p(t)-p(s))+(b_j-b_k)\la(p(x)-p(t))}\,dt,
\end{gather}
уфх яЁхфхы√ шэЄхуЁшЁютрэш  Ёртэ√ (ь√ ёўшЄрхь, ўЄю шэЄхуЁры Ёртхэ эєы■, хёыш эшцэшщ яЁхфхы сюы№°х тхЁїэхую)
\begin{equation}\label{eq:lim}
\begin{cases}\text{юЄ}\ x\ \text{фю}\ s,\qquad
&\text{яЁш}\ j,\,l<k;\\\text{юЄ}\ \max\{x,s\}\ \text{фю}\ 1,\qquad &\text{яЁш}\ j<k\le l;\\\text{юЄ}\ 0\ \text{фю}\ \min\{x,s\},\qquad &\text{яЁш}\ l<k\le j;\\
\text{юЄ}\ s\ \text{фю}\ x,\qquad &\text{яЁш}\ k\le j,\,l.
\end{cases}
\end{equation}
─ы  юЎхэъш юёЄрЄъют т яЁхфёЄртыхэшш \eqref{eq:matras} сєфхь шёяюы№чютрЄ№ ЇєэъЎшш
\begin{gather}\notag
\Upsilon(\la)=\Upsilon_\infty(\la)=\max_{j,\,k,\,l,\,s,\,x}|\upsilon_{jkl}(s,x,\la)|, \qquad
\Upsilon(x,\la)=\max_{j,\,k}|\upsilon_{jkk}(0,x,\la)|,\\
\Upsilon_\mu(\la)=\max_{j,\,k,\,l}\left(\int_0^1\int_0^1|\upsilon_{jkl}(s,x,\la)|^\mu\,ds\,dx\right)^{1/\mu}+
\max_{j,\,k}\left(\int_0^1|\upsilon_{jkk}(0,x,\la)|^\mu dx\right)^{1/\mu},\label{eq:defUps}
\end{gather}
уфх $\mu\in[1,\infty)$. ╦хуъю тшфхЄ№, ўЄю $\Upsilon_\mu(\la)\le\Upsilon(\la)$ ш $\Upsilon(x,\la)\le\Upsilon(\la)$ фы 
ы■сюую $x\in[0,1]$. ═рь яюЄЁхсєхЄё  ёыхфє■∙шщ рэрыюу ыхьь√ ╨шьрэр--╦хсхур.
\begin{Lemma}\label{lem:RL}
$\Upsilon(\la)\to0$ яЁш $\Dom\ni\la\to\infty$. ╩Ёюьх Єюую, т юсырёЄш $\Dom$ т√яюыэхэр юЎхэър
$\max_{j,k,l,s,x}|\varrho_{jkl}(s,x,\la)|\le e^{2hp}\gamma(\la)$.
\end{Lemma}
\begin{proof}
╟рЇшъёшЁєхь яЁюшчтюы№э√х шэфхъё√ $j$, $k$, $l$ ш Єюўъш $s,\,x\in[0,1]$. ╧Ёхцфх тёхую чрьхЄшь, ўЄю яЁш $l<k$ яЁхфхы√
шэЄхуЁшЁютрэш  т \eqref{eq:defups} ЁрёёЄртыхэ√ Єръ, ўЄю $t\le s$, р чэрўшЄ шч \eqref{eq:d} ш єёыютш  (iii) ёыхфєхЄ
$$
\Re(b_l-b_k)\la(p(t)-p(s))<h(p(s)-p(t)).
$$
╧Ёш $l\ge k$ шьххь $t\ge s$ ш
$$
\Re(b_l-b_k)\la(p(t)-p(s))<h(p(t)-p(s)),
$$
Єръ ўЄю т ы■сюь ёыєўрх
$$\Re(b_l-b_k)\la(p(t)-p(s))<h|p(s)-p(t)|\le hp.
$$
(чфхё№ ш фрыхх юсючэрўрхь $p=p(1)$.) └эрыюушўэю яюърч√трхь, ўЄю
$$
\Re(b_j-b_k)\la(p(x)-p(t))<h|p(x)-p(t)|\le hp.
$$
╬Єё■фр ёЁрчє ёыхфєхЄ юЎхэър $\max_{j,k,l,s,x}|\varrho_{jkl}(s,x,\la)|\le e^{2hp}\gamma(\la)$. ╬сЁрЄшьё  ъ ЇєэъЎшш
$\Upsilon(\la)$. ┬ ёшыє ыхьь√ \ref{lem:3}, $\upsilon_{jkl}\equiv0$, хёыш $b_j=b_l$. ┼ёыш цх $b_j\ne b_l$, Єю
\begin{equation*}%\label{eq:ups1}
|\upsilon_{jkl}(s,x,\la)|\le \left|\int q_{jl}(t)\exp\{(b_l-b_k)\la(p(t)-p(s))+(b_j-b_k)\la(p(x)-p(t))\}\,dt\right|,
\end{equation*}
уфх яЁхфхы√ т яюёыхфэхь шэЄхуЁрых ЁрёёЄртыхэ√ ёюуырёэю \eqref{eq:lim}. ╤фхырхь чрьхэє $\xi=p(t)$, $\xi\in[0,p]$, яюёых
ъюЄюЁющ шэЄхуЁры яЁшьхЄ тшф
\begin{equation}\label{eq:ups2}
\left|\int \frac{q_{jl}(t(\xi))}{\rho(t(\xi))}e^{(b_l-b_k)\la(\xi-p(s))+(b_j-b_k)\la(p(x)-\xi)}\,d\xi\right|.
\end{equation}
╧Ёхфхы√ шэЄхуЁшЁютрэш  чфхё№ ЁрёёЄртыхэ√, ёюуырёэю яЁртшыє \eqref{eq:lim}, ё єўхЄюь Єюую ўЄю эхЁртхэёЄтю $a\le t\le b$
¤ътштрыхэЄэю $p(a)\le\xi\le p(b)$. ╟рьхЄшь, ўЄю ЇєэъЎш  $\frac{q_{jl}(t(\xi))}{\rho(t(\xi))}$ ёєььшЁєхьр эр юЄЁхчъх
$[0,p]$, яюёъюы№ъє
$$
\int_0^{p}\frac{|q_{jl}(t(\xi))|}{\rho(t(\xi))}\,d\xi=\int_0^1|q_{jl}(t)|\,dt.
$$
─ы  фрээюую $\eps>0$ яюфсхЁхь эхяЁхЁ√тэю фшЇЇхЁхэЎшЁєхьє■ ЇєэъЎш■ $\wt q_{jl}(\xi)$ Єръ, ўЄю
$$
\int_0^{p}\left|\frac{q_{jl}(t(\xi))}{\rho(t(\xi))}-\wt q_{jl}(\xi)\right|\,d\xi<\eps.
$$
╧юёъюы№ъє тх∙хёЄтхээр  ўрёЄ№ т√Ёрцхэш  т яюърчрЄхых ¤ъёяюэхэЄ√ т шэЄхуЁрых \eqref{eq:ups2} яю-яЁхцэхьє эх яЁхтюёїюфшЄ
$2hp$, Єю шэЄхуЁры \eqref{eq:ups2} фюяєёърхЄ юЎхэъє
\begin{equation}\label{eq:ups3}
\le\eps e^{2hp}+\left|\int\wt q_{jl}(\xi)e^{(b_l-b_k)\la(\xi-p(s))+(b_j-b_k)\la(p(x)-\xi)}\,d\xi\right|.
\end{equation}
╧юёых шэЄхуЁшЁютрэш  яю ўрёЄ ь т юёЄрт°хьё  шэЄхуЁрых яЁшфхь ъ юЎхэъх
$$
\le\eps e^{2hp}+\frac{e^{2hp}}{|b_l-b_j|\cdot|\la|}\left(2\max_{0\le\xi\le p}|\wt q_{jl}(\xi)|+\int_0^{p}|\wt
q_{jl}'(\xi)|\,d\xi\right)\le2\eps e^{2hp}
$$
яЁш фюёЄрЄюўэю сюы№°юь $|\la|$.
\end{proof}

\section{╬ёэютэющ Ёхчєы№ЄрЄ фы  ёшёЄхь}\setcounter{equation}{0}

╧ЁшёЄєяшь ъ  т√тюфє рёшьяЄюЄшўхёъшї ЇюЁьєы фы  ЇєэфрьхэЄры№эющ ьрЄЁшЎ√ ёшёЄхь√ \eqref{eq:main}.
\begin{Theorem}\label{tm:main}
╨рёёьюЄЁшь ёшёЄхьє \eqref{eq:main}, ъю¤ЇЇшЎшхэЄ√ ъюЄюЁющ єфютыхЄтюЁ ■Є єёыютш ь (i) --- (iii). ╧єёЄ№ $\wt \Gamma_\kappa$
--- юфшэ шч ёхъЄюЁют, юяЁхфхыхээ√ї т√°х. ╥юуфр т юсырёЄш
$\Dom=\{\la\in\wt\Gamma_\kappa:\,|\la|>\la_0\}$ юяЁхфхыхэр ьрЄЁшЎр $\Y(x,\la)$ ЇєэфрьхэЄры№эющ ёшёЄхь√ Ёх°хэшщ
єЁртэхэш  \eqref{eq:main}, єфютыхЄтюЁ ■∙р  рёшьяЄюЄшўхёъюьє єёыютш■
\begin{gather}\label{eq:matras}
\Y(x,\la)=\Y^0(x,\la)+\cA(x,\la)\cE(x,\la),\quad \text{уфх}\\
\cA(x,\la)=(\al_{jk}(x,\la))_{j,\,k=1}^n,\qquad\max_{j,\,k,\,x}|\alpha_{jk}(x,\la)|\le
C(\Upsilon(\la)+\gamma(\la)),\notag
\end{gather}
р ЇєэъЎшш $\Y^0(x,\la)$ ш $\cE(x,\la)$ юяЁхфхыхэ√ т \eqref{eq:Y0}. ┼ёыш ¤ыхьхэЄ√ ьрЄЁшЎ√ $\A(x)-\D(x)$ ёєььшЁєхь√ т
ёЄхяхэш $\mu'\in(1,\infty)$, Єю, яюыюцшт $1/\mu+1/\mu'=1$, ьюцхь юЎхэшЄ№
\begin{gather}\label{eq:matras1}
\max_{j,\,k}|\alpha_{jk}(x,\la)|\le \Upsilon(x,\la)+C(\Upsilon_{\mu}(\la)+\gamma(\la)),\\
\max_{j,\,k}\left(\int_0^1|\alpha_{jk}(x,\la)|^{\mu}\,dx\right)^{1/\mu}\le
C(\Upsilon_{\mu}(\la)+\gamma(\la)).\label{eq:matras2}
\end{gather}
\end{Theorem}
\begin{Note}
╤ єўхЄюь ыхьь√ \ref{lem:RL},
$$
\Y(x,\la)=\left(\M(x,\la)+o(1)\right)\cE(x,\la),\qquad\wt\Gamma_\kappa\ni\la\to\infty,
$$
ЁртэюьхЁэю яю $x\in[0,1]$.
\end{Note}
\begin{proof}[─юърчрЄхы№ёЄтю ЄхюЁхь√ \ref{tm:main}.]\textit{╪ру 1. ╧хЁхїюф ъ ёшёЄхьх шэЄхуЁры№э√ї єЁртэхэшщ.}\qquad
═ряюьэшь, ўЄю єЁртэхэш  ёшёЄхь√ \eqref{eq:main} чрэєьхЁютрэ√ Єръ, ўЄюс√ т юсырёЄш $\Dom$ т√яюыэ ышё№ эхЁртхэёЄтр
\eqref{eq:d}. ┴єфхь шёърЄ№ Ёх°хэшх \eqref{eq:matrmain} т тшфх $\Y(x,\la)=\M(x)\Z(x,\la)\cE(x,\la)$. ╥юуфр, ё єўхЄюь
ЁртхэёЄт $\cE'=\la\rho\B\cE$ ш $\M'=\D\M$,
\begin{gather*}
\M'\Z\cE+\M\Z'\cE+\M\Z\cE'=\la\rho\B\M\Z\cE+\D\M\Z\cE+(\A-\D+\C)\M\Z\cE\quad\Longleftrightarrow\\
\M\Z'\cE+\la\rho\M\Z\B\cE=\la\rho\B\M\Z\cE+(\A-\D+\C)\M\Z\cE.
\end{gather*}
─юьэюцшь ¤Єю ЁртхэёЄтю эр $\M^{-1}$ ёяЁртр ш эр $\cE^{-1}$ ёыхтр ш тюёяюы№чєхьё  ЁртхэёЄтюь $\M^{-1}\B=\B\M^{-1}$. ┬
Ёхчєы№ЄрЄх эр°х єЁртэхэшх яЁшьхЄ тшф
$$
\Z'=\la\rho(\B\Z-\Z\B)+(\Q+\R)\Z
$$
(эряюьэшь, ўЄю $\Q(x):=\M^{-1}(x)(\A(x)-\D(x))\M(x)$, $\R(x,\la):=\M^{-1}(x)\C(x,\la)\M(x)$). ┬ ъююЁфшэрЄэющ чряшёш
яюыєўхээюх єЁртэхэшх шьххЄ тшф
\begin{equation}\label{eq:mainzdif}
z_{jk}'(x,\la)=\la(b_j-b_k)\rho(x) z_{jk}(x,\la)+\sum_{l=1}^n\rv_{jl}(x,\la)z_{lk}(x,\la).
\end{equation}
╤ўшЄр  фрыхх шэфхъё $k$ ЇшъёшЁютрээ√ь, яЁютхфхь шэЄхуЁшЁютрэшх т \eqref{eq:mainzdif}, т√сшЁр  эрўры№э√х єёыютш 
$z_{jk}(1,\la)=0$ фы  $j<k$, $z_{jk}(0,\la)=0$ фы  $j>k$ ш $z_{kk}(0,\la)=1$:
\begin{align}
&z_{jk}(x,\la)=-\sum_{l}\int_x^1\rv_{jl}(t,\la)e^{(b_j-b_k)\la(p(x)-p(t))}z_{lk}(t,\la)\,dt,\qquad j<k,\notag\\
&z_{kk}(x,\la)-1=\sum_{l}\int_0^x\rv_{kl}(t,\la)z_{lk}(t,\la)\,dt, \label{eq:mainz}\\
&z_{jk}(x,\la)=\sum_{l}\int_0^x\rv_{jl}(t,\la)e^{(b_j-b_k)\la(p(x)-p(t))}z_{lk}(t,\la)\,dt,\qquad j>k.\notag
\end{align}
╬сючэрўр  ўхЁхч $\V_k(\la)$ шэЄхуЁры№э√щ юяхЁрЄюЁ т яЁртющ ўрёЄш ¤Єющ ёшёЄхь√, яюыєўшь $\z_k=\z^0_k+\V_k\z_k$, уфх
$\z^0_k=(\delta_j^k)_{j=1}^n$ (ўхЁхч $\z_k$ ь√ юсючэрўрхь $k$--Є√щ ёЄюысхЎ ьрЄЁшЎ√ $\Z$). ┬тхфхь Єръцх юяхЁрЄюЁ√
$\Q_k(\la)$ ш $\R_k(\la)$, юяЁхфхыхээ√х яЁртющ ўрёЄ№■ \eqref{eq:mainz} ё чрьхэющ ЇєэъЎшщ $\rv_{jk}$ эр $q_{jk}$ ш
$r_{jk}$ ёююЄтхЄёЄтхээю. ╧юёъюы№ъє $\rv_{jk}(x,\la)=q_{jk}(x)+r_{jk}(x,\la)$, Єю $\V_k=\Q_k+\R_k$.

\noindent\textit{╪ру 2. ╬яхЁрЄюЁ $(V_k(\la))^2$ --- ёцрЄшх.}\qquad ╨х°хэшх ёшёЄхь√ \eqref{eq:mainz} ь√ сєфхь шёърЄ№ т
тшфх Ё фр $\z_k=\sum_{\nu=0}(\V_k(\la))^\nu\z^0_k$, ёїюф ∙хуюё  т яЁюёЄЁрэёЄтх $(\AC[0,1])^n$. ─ы  фюърчрЄхы№ёЄтр
ёїюфшьюёЄш Ё фр ш юЎхэъш юёЄрЄъют т \eqref{eq:matras} эрь эєцэю юЎхэшЄ№ эюЁь√ юяхЁрЄюЁют $\Q_k(\la)$ ш $\R_k(\la)$. ╧Ёш
¤Єюь, яюёъюы№ъє ь√ їюЄшь яюыєўшЄ№ Єръцх эхЁртхэёЄтр \eqref{eq:matras1}, эр°ш юЎхэъш фюыцэ√ єўшЄ√трЄ№ чэрўхэшх
$\mu'\in[1,\infty]$. ╬ЄьхЄшь, ўЄю юЎхэъш фы  $\|\Q_k(\la)\|$ ш фы  $\|\R_k(\la)\|$ шьх■Є Ёрчышўэє■ яЁшЁюфє. ╧хЁт√х
ёт чрэ√ ё ыхььющ ╨шьрэр--╦хсхур, р тЄюЁ√х --- ё єс√трэшхь эюЁь√ $\|\C(\cdot,\la)\|_{L_1[0,1]}=\gamma(\la)$ яЁш
$\la\to\infty$. ═рўэхь ё юЎхэюъ
\begin{gather}
\|\R_k(\la)\|_{L_\infty\to\AC}\le2e^{hp+2a}\gamma(\la),\qquad\|\V_k(\la)\|_{L_\infty\to\AC}\le2e^{hp+2a}(a+c),\notag\\
\|\Q_k(\la)\|_{L_\mu\to\AC}\le 2e^{hp+2a}\|\A(x)-\D(x)\|_{L_{\mu'}}.\label{eq:RQ1}
\end{gather}
╟рьхЄшь, ўЄю яЁхфхы√ шэЄхуЁшЁютрэш  т \eqref{eq:mainz} ЁрёёЄртыхэ√ Єръ, ўЄю яЁш $\la\in\wt\Gamma_\kappa$
$$
\Re(b_j-b_k)\la(p(x)-p(t))<h|p(x)-p(t)|\le hp.
$$
╧юыюцшь $\f\in L_\infty[0,1]$, $\g_k=\R_k(\la)\f\in\AC[0,1]$. ╚ч \eqref{eq:qest} ш \eqref{eq:mainz} ёыхфєхЄ, ўЄю
\begin{multline*}
\|\g_k(x,\la)\|_{\AC}=\sum_{j=1}^n\int_0^1|g_{jk}(x,\la)|+|g'_{jk}(x,\la)|\,dx\le\\\le
2e^{hp}\sum_{j=1}^{n}\sum_{l=1}^n \int_0^1|r_{jl}(t,\la)||f_{l}(t)|\,dt\le 2e^{hp+2a}\gamma(\la)\|\f\|_{L_\infty}
\end{multline*}
ш яхЁтюх эхЁртхэёЄтю т \eqref{eq:RQ1} фюърчрэю. ─юърчрЄхы№ёЄтю тЄюЁюую яюыэюёЄ№■ рэрыюушўэю. ╬ёЄрхЄё  юЎхэшЄ№ эюЁьє
юяхЁрЄюЁр $\Q_k(\la)$. ╧юыюцшь $\f\in L_\mu[0,1]$, $\g_k=\Q_k(\la)\f\in\AC[0,1]$ ш, шёяюы№чє  \eqref{eq:qest} ш
\eqref{eq:mainz}, яюыєўшь
\begin{multline*}
\|\g_k(x,\la)\|_{\AC}=\sum_{j=1}^n\int_0^1|g_{jk}(x,\la)|+|g'_{jk}(x,\la)|\,dx\le\\\le
2e^{hp}\sum_{j=1}^{n}\sum_{l=1}^n \int_0^1|q_{jl}(t)||f_{l}(t)|\,dt\le 2e^{hp}\|\Q(x)\|_{L_{\mu'}}\|\f\|_{L_\mu}.
\end{multline*}
╧юърцхь, ўЄю яЁш фюёЄрЄюўэю сюы№°юь $|\la|$, $\la\in\Dom$, юяхЁрЄюЁ $(\V_k(\la))^2$  ты хЄё  ёцрЄшхь т $L_\infty[0,1]$.
┼ёыш $\f\in L_\infty$, $\g_k=\V_k^2(\la)\f$, Єю
$$
g_{jk}(x,\la)=\!\!\!\sum_{l,\,m=1}^n(\pm)_{jk}(\pm)_{lk}\!\int\!\!\left(\rv_{jl}(t,\la)e^{(b_j-b_k)\la(p(x)-p(t))}\!\int
\rv_{lm}(s,\la)e^{(b_l-b_k)\la(p(t)-p(s))}f_m(s)\,ds\right)dt.
$$
┬эх°эшх шэЄхуЁры√ т√ўшёы ■Єё  яю юЄЁхчъє $[x,1]$ яЁш $j<k$, р яЁш $j\ge k$
--- яю юЄЁхчъє $[0,x]$. ┬эєЄЁхээшх шэЄхуЁры√ --- яю $[t,1]$ яЁш $l<k$ ш яю $[0,t]$ яЁш $l\ge k$. ╧хЁхёЄртшт шэЄхуЁры√,
яюыєўшь
\begin{equation}\label{eq:Vsquare}
g_{jk}(x,\la)=\sum_{m=1}^n\int_0^1\left(\sum_{l=1}^{n}\rv_{lm}(s,\la)(\upsilon_{jkl}(s,x,\la)+\varrho_{jkl}(s,x,\la))\right)f_m(s)\,ds.
\end{equation}
╤юуырёэю эхЁртхэёЄтє \eqref{eq:qest} ш ыхььх \ref{lem:RL},
\begin{equation*}%\label{eq:Tsqest}
\sup_{0\le
x\le1}\int_0^1\left|\rv_{lm}(s,\la)(\upsilon_{jkl}(s,x,\la)+\varrho_{jkl}(s,x,\la))\right|\,ds\le(a+\gamma(\la))e^{2a}(\Upsilon(\la)+\gamma(\la))\to0
\end{equation*}
яЁш $\Dom\ni\la\to\infty$. ╥юуфр
\begin{equation*}
\|(\V_k(\la))^2\f\|_{L_\infty}=\max_{j,x}|g_{jk}(x,\la)|\le
n^2(a+\gamma(\la))e^{2a}(\Upsilon(\la)+\gamma(\la))\max_{m,s}|f_m(s)|.%,\qquad\text{Є.х.}\ \
%\|(\V_k(\la))^2\|_{L_\infty}\le(a+\gamma(\la))e^{2a}\Upsilon(\la).
\end{equation*}
╬сючэрўшт $C_0=n^2(a+c)e^{2a}$, яюыєўшь
\begin{equation}\label{eq:Vsqnorm}
\|(\V_k(\la))^2\|_{L_\infty}\le C_0(\Upsilon(\la)+\gamma(\la))\to0\quad \text{яЁш}\ |\la|\to\infty,
\end{equation}
Є.х. юяхЁрЄюЁ $(\V_k(\la))^2$ т яЁюёЄЁрэёЄтх $L_\infty$ --- ёцрЄшх яЁш сюы№°шї $|\la|$.

\noindent\textit{╪ру 3. ╚Єюуют√х юЎхэъш.}\qquad ┬хЁэхьё  ъ ёшёЄхьх \eqref{eq:mainz}. ╧ЁхфёЄртшь хх Ёх°хэшх т тшфх (яюър
ЇюЁьры№эюую) Ё фр
$$
\z_k=\z_k^0+\V_k(\la)\sum_{\nu=0}^\infty(V_k(\la))^\nu\z_k^0.
$$
╚ч ¤Єющ чряшёш ш шч юуЁрэшўхээюёЄш юяхЁрЄюЁр $V_k(\la)$, фхщёЄтє■∙хую шч $L_\infty$ т $\AC$ (юуЁрэшўхээюёЄ№ юяхЁрЄюЁр
фюърчрэр т \eqref{eq:RQ1}) тшфэю, ўЄю фюёЄрЄюўэю фюърчрЄ№ ёїюфшьюёЄ№ Ё фр яю эюЁьх $L_\infty[0,1]$. ─ы  ¤Єюую яхЁхяш°хь
хую т тшфх
$$
\z_k=\z^0_k+\V_k\z^0_k+\sum_{\nu=0}^\infty(\V_k(\la))^{2\nu}\left((\V_k(\la))^2\z_k^0+(\V_k(\la))^3\z_k^0\right).
$$
╙ўшЄ√тр  \eqref{eq:Vsqnorm}, єтхышўшь ўшёыю $\la_0$ Єръ, ўЄюс√ яЁш тёхї $\la\in \Dom$ с√ыю т√яюыэхэю
$\|(\V_k(\la))^2\|_{L_\infty}<1/2$. ▌Єю урЁрэЄшЁєхЄ ёїюфшьюёЄ№ Ё фр т эюЁьх $L_\infty$ ш фрхЄ юЎхэъє
\begin{equation}\label{eq:estz1}
\|\z_k-\z_k^0-\V_k(\la)\z_k^0\|_{L_\infty}\le 2\left\|(\V_k(\la))^2\z_k^0+(\V_k(\la))^3\z_k^0\right\|_{L_\infty}.
\end{equation}
┬√ўшёышь ЇєэъЎш■ $\wt\z^1_k=(\wt z^1_{jk}(x,\la))_{j=1}^{n}=\Q_k(\la)\z^0_k$ :
\begin{align}
&\wt z_{jk}^1(x,\la)=-\int_x^1q_{jk}(t,\la)e^{(b_j-b_k)\la(p(x)-p(t))}\,dt,\qquad j<k,\label{eq:estz2}\\
&\wt z_{jk}^1(x,\la)=\int_0^xq_{jk}(t,\la)e^{(b_j-b_k)\la(p(x)-p(t))}\,dt,\qquad j\ge k.\notag
\end{align}
╥ръшь юсЁрчюь, шьххь $\wt z^1_{jk}(x,\la)=\upsilon_{jkk}(0,x,\la)$ ш
\begin{equation}\label{eq:z1}
\max_{1\le j,\,k\le n}|\wt z_{jk}^1(x,\la)|\le\Upsilon(x,\la).
\end{equation}
╬Єё■фр ёыхфєхЄ юЎхэъш $\|\wt\z^1_k\|_{L_\mu}\le \Upsilon_\mu(\la)$ ш
\begin{equation}\label{eq:estz3}
\|\V_k\z^0_k\|_{L_\mu}\le\|\wt \z^1_k\|_{L_\mu}+\|\R_k(\la)\z_k^0\|_{L_\infty}\le
C_1(\Upsilon_\mu(\la)+\gamma(\la)),\quad C=2e^{hp+2a}.
\end{equation}
%C_1=2e^{hp+2a}
╙ёыютшьё  фры№°х юсючэрўрЄ№ $\z^\nu_k=(V_k)^\nu\z_k^0$. ╥хяхЁ№ чрьхЄшь, ўЄю
\begin{gather}\notag
\|\z_k^2\|_{\AC}\le\|\Q_k\|_{L_\mu\to
\AC}\|\z^1_k\|_{L_\mu}+\|\R_k\|_{L_\infty\to\AC}\|\z^1_k\|_{L_\infty} \le\\
\le C_1\|\Q_k\|_{L_\mu\to\AC}(\Upsilon_\mu(\la)+\gamma(\la))+2e^{hp+2a}\gamma(\la)\|\V_k\|_{L_\infty}\le
C_2(\Upsilon_\mu(\la)+\gamma(\la)),\label{eq:estz5}
\end{gather}
уфх $C_2=4e^{hp+4a}(a+c+m)$, $m=\|\A-\D\|_{L_{\mu'}}$. ═ръюэхЎ,
\begin{equation}\label{eq:estz6}
\|\z_k^3\|_{L_\infty}\le\|\V_k(\la)\|_{L_\infty}C_2(\Upsilon_\mu(\la)+\gamma(\la))\le
C_3(\Upsilon_\mu(\la)+\gamma(\la)),
\end{equation}
$C_3=8e^{3hp+6a}(a+c)(a+c+m)$. ╧юфёЄрты   яюыєўхээ√х юЎхэъш т \eqref{eq:estz1}, яюыєўрхь
\begin{equation}\label{eq:estz4}
\|\z_k-\z_k^0-\z_k^1\|_{L_\infty}\le C_4(\Upsilon_\mu(\la)+\gamma(\la)),
\end{equation}
$C_4=2C_2+2C_3$, р єўшЄ√тр  \eqref{eq:estz3},
\begin{equation}\label{eq:estz7}
\|\z_k-\z_k^0\|_{L_\mu}\le C_5(\Upsilon_\mu(\la)+\gamma(\la)),
\end{equation}
$C_5=C_4+2e^{hp+2a}$. ╬Ўхэъш \eqref{eq:matras2} ёыхфє■Є юЄё■фр, яюёых фюьэюцхэш  ьрЄЁшЎ√ $\Z$ ёыхтр эр
$\M(x)\in\AC[0,1]$ (яЁш ¤Єюь ъюэёЄрэЄр т \eqref{eq:matras2} Ёртэр $C=C_5e^a$). ╬Ўхэъш \eqref{eq:matras} --- ўрёЄэ√щ
ёыєўрщ \eqref{eq:matras2} яЁш $\mu=\infty$. ─ы  фюърчрЄхы№ёЄтр \eqref{eq:matras1} чрьхЄшь, ўЄю
$$
|\z_k^1(x,\la)|\le|\wt\z^1_k(x,\la)|+\|\R_k(\la)\z_k^0\|_{L_\infty}\le \Upsilon(x,\la)+2e^{hp+2a}\gamma(\la).
$$
╬Єё■фр ш шч \eqref{eq:estz4} яюыєўрхь \eqref{eq:matras1}.
\end{proof}
%\begin{Note}
%╫шёыю $\la_0$ ш ъюэёЄрэЄ√ т юЎхэърї \eqref{eq:matras}--\eqref{eq:matras2} ьюуєЄ с√Є№ яЁхф· тыхэ√. └эрышчшЁє 
%фюърчрЄхы№ёЄтю ЄхюЁхь√ \ref{tm:main}, яюыєўрхь, ўЄю $\la_0$ т√сЁрэю Єръ, ўЄюс√ яЁш тёхї $\la\in\wt\Gamma_\kappa$,
%$|\la|>\la_0$ т√яюыэ ышё№ эхЁртхэёЄтр \eqref{eq:d} ш юЎхэъш $n^2(a+c)e^{2a}(\Upsilon(\la)+\gamma(\la))\le1/2$,
%$\gamma(\la)\le1$. ╩юэёЄрэЄ√ т \eqref{eq:matras}--\eqref{eq:matras2} ьюцэю яюыюцшЄ№ Ёртэ√ьш
%$C=3(8m^2+20m+9)e^{4hp+6a}$, уфх $m=\|\A(x)-\D(x)\|_{L_{\mu'}}$.
%\end{Note}
\begin{Note}
┼ёыш ¤ыхьхэЄ√ ьрЄЁшЎ√ $\C$ --- ЇєэъЎшш $c_{jk}(x,\la)$ уюыюьюЁЇэ√ яю ярЁрьхЄЁє $\la$ т юсырёЄш $\Dom_{\kappa,\la_0}$,
Єю ЇєэъЎш  $\Y(x,\la)$ Єръцх уюыюьюЁЇэр т $\Dom_{\kappa,\la_0}$ (ъръ ¤ыхьхэЄ яЁюёЄЁрэёЄтр $(\AC[0,1])^{n^2}$).
\end{Note}
─юърчрЄхы№ёЄтю яюыэюёЄ№■ рэрыюушўэю фюърчрЄхы№ёЄтє чрьхўрэш  \ref{note:2}.
\begin{Note}
╬яЁхфхышЄхы№ ьрЄЁшЎ√ $\Y(x,\la)$ шьххЄ тшф
$$
\det \Y(x,\la)=\det\Y^0(x,\la)(\det\M(x)+o(1))\ne0
$$
яЁш $\la\in\Dom_{\kappa,\la_0}$ яЁш фюёЄрЄюўэю сюы№°юь $\la_0$, Є.х. эрщфхээр  эрьш $\Y(x,\la)$ фхщёЄтшЄхы№эю  ты хЄё 
ЇєэфрьхэЄры№эющ ьрЄЁшЎхщ Ёх°хэшщ ёшёЄхь√ \eqref{eq:main} т $\Dom$.
\end{Note}
┬ ЄхюЁхьх \ref{tm:main} ь√, ЇръЄшўхёъш, эр°ыш ыш°№ ёЄрЁ°шщ ўыхэ т рёшьяЄюЄшўхёъюь яЁхфёЄртыхэшш ьрЄЁшЎ√ $\Y(x,\la)$.
─рфшь эхёъюы№ъю єЄюўэхэшщ рёшьяЄюЄшўхёъшї ЇюЁьєы \eqref{eq:matras}--\eqref{eq:matras2}.
\begin{Theorem}\label{tm:app}
┬ єёыютш ї ЄхюЁхь√ \ref{tm:main}, т юс∙хщ ёшЄєрЎшш, ъюуфр $\A(x)\in L_1[0,1]$,
\begin{gather}\label{eq:matras3}
\Y(x,\la)=\M(x,\la)\Z(x,\la)\cE(x,\la),\quad \text{уфх}\quad \Z(x,\la)=I+\Z^1(x,\la)+\Z^2(x,\la)+\EuScript A(x,\la),\\
\EuScript A(x,\la)=(\al_{jk}(x,\la))_{j,\,k=1}^n,\qquad\max_{j,\,k,\,x}|\alpha_{jk}(x,\la)|\le
C(\Upsilon(\la)+\gamma(\la))^2,\notag\\
\label{eq:matras35} Z(x,\la)=I+\Z^1(x,\la)+\EuScript
A(x,\la),\qquad\text{уфх}\quad\max_{j,\,k}\|\alpha_{jk}(x,\la)\|_{\AC}\le C(\Upsilon+\gamma(\la)).
\end{gather}
┼ёыш ¤ыхьхэЄ√ ьрЄЁшЎ√ $\A(x)-\D(x)$ ёєььшЁєхь√ т ёЄхяхэш $\mu'\in(1,\infty)$, Єю
\begin{gather}\label{eq:matras4}
\Z(x,\la)=I+\Z^1(x,\la)+\EuScript A(x,\la))),\qquad \text{уфх}\quad\max_{j,\,k,\,x}|\alpha_{jk}(x,\la)|\le
C(\Upsilon_\mu(\la)+\gamma(\la)),\\ \notag
\Z(x,\la)=I+\Z^1(x,\la)+\Z^2(x,\la)+\Z^3(x,\la)+\EuScript A(x,\la))),\\
\text{уфх}\quad\max_{j,\,k}\|\alpha_{jk}(x,\la)\|_{L_\mu}\le C(\Upsilon_\mu(\la)+\gamma(\la))^2.\label{eq:matras5}
\\ \notag
\Z(x,\la)=I+\Z^1(x,\la)+\Z^2(x,\la)+\Z^3(x,\la)+\Z^4(x,\la)+\EuScript A(x,\la))),\\
\text{уфх}\quad\max_{j,\,k}\|\alpha_{jk}(x,\la)\|_{\AC}\le C(\Upsilon_\mu(\la)+\gamma(\la))^2.\label{eq:matras6}
\end{gather}
╟фхё№ ўхЁхч $\Z^\nu(x,\la)$ ь√ юсючэрўрхь ьрЄЁшЎ√, ёюёЄртыхээ√х шч ёЄюысЎют $\z_k^\nu(x,\la)$.
\end{Theorem}
\begin{proof}
╧хЁт√х ЄЁш єЄтхЁцфхэш  ёыхфє■Є шч, яюыєўхээ√ї яЁш фюърчрЄхы№ёЄтх ЄхюЁхь√ \ref{tm:main}, юЎхэюъ. ─ы  т√тюфр
\eqref{eq:matras3} фюёЄрЄюўэю эхЁртхэёЄт \eqref{eq:estz3} ш \eqref{eq:Vsqnorm}, \eqref{eq:matras4} ёыхфєхЄ шч
\eqref{eq:estz4}. ╟ряшёрт ЁртхэёЄтю
$$
\z_k-\z_k^0-V_k\z_k^0=\V_k(\z_k-\z_k^0)
$$
ш тюёяюы№чютрт°шё№ \eqref{eq:estz7} фы  $\mu=\infty$ ш \eqref{eq:RQ1}, яЁшфхь ъ \eqref{eq:matras35}. ─ы  фюърчрЄхы№ёЄтр
яюёыхфэшї фтєї єЄтхЁцфхэшщ ЄЁхсє■Єё  фюяюыэшЄхы№э√х Ёрёёєцфхэш .

╧юърцхь, ўЄю яЁш $\mu'>1$ ёцрЄшхь  ты хЄё  ш юяхЁрЄюЁ $(\Q_k(\la))^2$, фхщёЄтє■∙шщ шч $L_\infty$ єцх т $L_\mu$. ┬юч№ьхь
яЁюшчтюы№эє■ ЇєэъЎш■ $\f\in L_\infty$, яюыюцшь $\g_k=(\Q_k(\la))^2\f$ ш, Єръ цх, ъръ ш фы  юяхЁрЄюЁр $(\V_k(\la))^2$,
яЁшфхь ъ ЁртхэёЄтрь \eqref{eq:Vsquare} (ё чрьхэющ $\rv_{lm}(s,\la)$ эр $q_{lm}(s)$). ╥юуфр
\begin{multline*}
\left(\int_0^1\left|\int_0^1q_{lm}(s)\upsilon_{jkl}(s,x,\la)\,ds\right|^\mu dx\right)^{\frac1{\mu}}\le\\
\le\left(\int_0^1\left(\int_0^1|q_{lm}(s)|^{\mu'}ds\right)^{\frac{\mu}{\mu'}}\int_0^1|\upsilon_{jkl}(s,x,\la)|^{\mu}ds\right)^\frac1{\mu}
\le\|q_{lm}\|_{L_{\mu'}}\Upsilon_\mu(\la),
\end{multline*}
\begin{multline*}
\int_0^1|g_{jk}(x,\la)|^\mu\,dx\le
n^{2/\mu'}\|\f\|^\mu_{L_\infty}\sum_{l,\,m=1}^n\int_0^1\left|\int_0^1q_{lm}(s)\upsilon_{jkl}(s,x,\la)\,ds\right|^\mu
dx\le\\\le n^{2/\mu'}\|\f\|^\mu_{L_\infty}\left(\|\Q(x)\|_{L_{\mu'}}\Upsilon_\mu(\la)\right)^\mu,
\end{multline*}
юЄъєфр $\|\g_k\|_{L_\mu}\le n^2\|\Q(x)\|_{L_{\mu'}}\Upsilon_\mu(\la)\|\f\|_{L_\infty}$, Є.х.
\begin{equation}\label{eq:Qsqnorm}
\|(\Q_k(\la))^2\|_{L_\infty\to L_\mu}\le n^2e^{2a}\|\A(x)-\D(x)\|_{L_{\mu'}}\Upsilon_\mu(\la).
\end{equation}
╧юыюцшт $C_6=n^2e^{2a}\max\{a+c,\|\A(x)-\D(x)\|_{L_{\mu'}}\}$, шьххь $\|(\Q_k(\la))^2\|_{L_\infty\to L_\mu}\le
C_6\Upsilon_\mu(\la)$. ┬ ърўхёЄтх ёыхфёЄтш  яюыєўрхь юЎхэъє
\begin{gather*}
\|(\V_k(\la))^2\|_{L_\infty\to L_\mu}\le\|(\Q_k(\la))^2\|_{L_\infty\to L_\mu}+\|\Q_k(\la)\R_k(\la)\|_{L_\infty\to
L_\mu}+\|\R_k(\la)\Q_k(\la)\|_{L_\infty\to L_\mu}+\\+\|(\R_k(\la))^2\|_{L_\infty\to L_\mu} \le
C_6\Upsilon_\mu(\la)+8e^{2hp+2a}\|\A(x)-\D(x)\|_{L_\mu'}\gamma(\la)+4e^{2hp}\gamma^2(\la),
\end{gather*}
Є.х. яЁш фюёЄрЄюўэю сюы№°шї $|\la|>\la_0$, ъюуфр $\gamma(\la)<1$,
\begin{equation}\label{eq:Vsqnorm1}
\|(\V_k(\la))^2\|_{L_\infty\to L_\mu}\le C_7(\Upsilon_{\mu}(\la)+\gamma(\la)).
\end{equation}
┬ \eqref{eq:estz5} ш \eqref{eq:estz6} ь√ юЎхэшыш $\|\z_k^2\|_{L_\infty}$ ш $\|\z_k^3\|_{L_\infty}$ тхышўшэющ
$C(\Upsilon_\mu(\la)+\gamma(\la))$. ╬ёЄрхЄё  єўхёЄ№, ўЄю $\|(\V_k(\la)^2\|_{L_\infty}\le1/2$, р чэрўшЄ
$$
\left\|\sum_{\nu=0}^\infty(\V_k(\la))^{2\nu}(\z_k^2+\z_k^3)\right\|_{L_\infty}\le 2\|\z_k^2+\z_k^3\|_{L_\infty}\le
C_8(\Upsilon_\mu(\la)+\gamma(\la)).
$$
╧Ёшьхэ   юяхЁрЄюЁ $(\V_k(\la))^2$ ш єўшЄ√тр  \eqref{eq:Vsqnorm1}, яюыєўшь
$$
\left\|\sum_{\nu=1}^\infty(\V_k(\la))^{2\nu}(\z_k^2+\z_k^3)\right\|_{L_\mu}\le C_9(\Upsilon_\mu(\la)+\gamma(\la))^2
$$
ш \eqref{eq:matras5} фюърчрэю. ─ы  фюърчрЄхы№ёЄтр \eqref{eq:matras6} юЎхэшь тэрўрых $\|\z_k^5\|_{L_\infty}$.
╨рёъырф√тр  юяхЁрЄюЁ $\V_k$ т ёєььє $\V_k=\Q_k+\R_k$, яЁхфёЄртшь $(\V_k)^5$ т тшфх ёєьь√ ЄЁшфЎрЄш фтєї ёырурхь√ї,
ърцфюх шч ъюЄюЁ√ї хёЄ№ яЁюшчтхфхэшх ўхЄ√Ёхї юяхЁрЄюЁют. ╤ырурхь√х, т ъюЄюЁ√х ьэюцшЄхы№ $\R_k$ тїюфшЄ фтр ш сюыхх Ёрч
юЎхэштр■Єё  тхышўшэющ $C\gamma^2(\la)$ (ёь. \eqref{eq:RQ1}). ╬ёЄры№э√х ёырурхь√х юЎхэшь, шёяюы№чє  \eqref{eq:RQ1},
\eqref{eq:Qsqnorm} ш \eqref{eq:z1}:
\begin{gather*}
\|(\Q_k)^5\z_k^0\|_{\AC}\le\|\Q_k\|_{L_\mu\to \AC}\|(\Q_k)^2\|_{L_\infty\to
L_\mu}\|\Q_k\|_{L_\mu\to L_\infty}\|\wt\z_k^1\|_{L_\mu}\le C\Upsilon^2_\mu(\la),\\
\|(\Q_k)^4\R_k\z_k^0\|_{\AC}\le\|(\Q_k)^2\|_{L_\mu\to \AC}\|(\Q_k)^2\|_{L_\infty\to L_\mu}\|\R_k\|_{L_\infty}\le
C\Upsilon_\mu(\la)\gamma(\la),\\
\|\R_k\Q_k\z_k^0\|_{\AC}\le\|\R_k\|_{L_\infty\to\AC}\|\wt\z_k^1\|_{L_\infty}\le C\Upsilon_\mu(\la)\gamma(\la).
\end{gather*}
╥ръшь юсЁрчюь, $\|\z_k^5\|_{\AC}\le C(\Upsilon_\mu(\la)+\gamma(\la))^2$. ╧Ёшьхэ   юяхЁрЄюЁ $\V_k$ тшфшь, ўЄю Єрър  цх
юЎхэър тхЁэр ш фы  $\z_k^6$ ш $\z_k^7$. ╧юы№чє ё№ \eqref{eq:RQ1} ш \eqref{eq:Vsqnorm}, єтхышўшь $\la_0$ эрёЄюы№ъю, ўЄю
$$
\|(\V_k(\la))^3\|_{\AC}\le\|\V_k\|_{L_\infty\to\AC}\|(\V_k)^2\|_{L_\infty}\le\frac12.
$$
╥хяхЁ№ \eqref{eq:matras6} ёыхфєхЄ шч юЎхэъш
$$
\left\|\sum_{\nu=0}^\infty(\V_k(\la))^{3\nu}(\z_k^5+\z_k^6+\z_k^7)\right\|_{\AC}\le 2\|\z_k^5+\z_k^6+\z_k^7\|_{\AC}\le
C(\Upsilon_\mu(\la)+\gamma(\la))^2.
$$
\end{proof}
╥хюЁхь√ \ref{tm:main} ш \ref{tm:app} фр■Є юЎхэъш юёЄрЄюўэ√ї ўыхэют т рёшьяЄюЄшўхёъюь яЁхфёЄртыхэшш ьрЄЁшЎ√ $\Y(x,\la)$
т ЄхЁьшэрї ЇєэъЎшщ $\Upsilon(\la)$, $\Upsilon_\mu(\la)$, $\Upsilon(x,\la)$ ш $\gamma(\la)$. ┬ юс∙хщ ёшЄєрЎшш, ъюуфр
т√яюыэхэ√ єёыютш  (i)--(iv) эр ъю¤ЇЇшЎшхэЄ√ ёшёЄхь√ \eqref{eq:main} ¤Єш ЇєэъЎшш єс√тр■Є ъ эєы■ яЁш
$\Dom\ni\la\to\infty$ (ёь. ыхььє \ref{lem:RL}). ╒рЁръЄхЁ ёЄЁхьыхэш  ъ эєы■ ЇєэъЎшш $\gamma(\la)$ Ўхышъюь юяЁхфхы хЄё 
ёЄхяхэ№■ єс√трэш  ЇєэъЎшщ $c_{ij}(\cdot,\la)$. ╩ръ тшфэю шч ЄхюЁхь√ \ref{tm:ho}, фы  ёыєўр  юяхЁрЄюЁют т√ёюъюую яюЁ фър
шьххь $\gamma(\la)=O(|\la|^{-1})$. ╙с√трэшх ъ эєы■ ЇєэъЎшщ $\Upsilon(\la)$, $\Upsilon_\mu(\la)$ ш $\Upsilon(x,\la)$
юяЁхфхы хЄё  "<Ёхуєы ЁэюёЄ№■"> ЇєэъЎшш $\rho(x)$ ш ¤ыхьхэЄют ьрЄЁшЎ√ $\A(x)-\D(x)$. ╧Ёш яют√°хэшш ёєььшЁєхьюёЄш
$\mu'\in(1,2]$ ¤Єшї ЇєэъЎшщ, ьюцэю юЎхэштрЄ№ эюЁьє $\Upsilon_\mu(\la)$ ш $\Upsilon(x,\la)$ т яЁюёЄЁрэёЄтрї ╒рЁфш
$H^{\mu}$. ╧Ёш яют√°хэшш уырфъюёЄш, т√Ёрцхээющ т юЎхэърї эр ьюфєы№ эхяЁхЁ√тэюёЄш шыш эр эюЁьє т яЁюёЄЁрэёЄтрї ┴хёютр
ЇєэъЎшш $\rho(x)$ ш ¤ыхьхэЄют ьрЄЁшЎ√ $\A(x)-\D(x)$, ьюцэю фртрЄ№ юЎхэъш эр ёъюЁюёЄ№ єс√трэш  ЇєэъЎшщ $\Upsilon(\la)$,
$\Upsilon_\mu(\la)$ ш $\Upsilon(x,\la)$. ╤ююЄтхЄёЄтє■∙шх ЄхюЁхь√ сєфєЄ яЁштхфхэ√ т ёыхфє■∙хщ ёЄрЄ№х ртЄюЁют.

\section{└ёшьяЄюЄшъш фы  юс√ъэютхээ√ї фшЇЇхЁхэЎшры№э√ї єЁртэхэшщ ё  ъю¤ЇЇшЎшхэЄрьш-ЁрёяЁхфхыхэш ьш}
\setcounter{equation}{0}

┬ ¤Єюь ярЁруЁрЇх ь√ яюыєўшь эєцэ√х рёьяЄюЄшъш фы  ЇєэфрьхэЄры№эющ ёшёЄэхь√ Ёх°хэшщ єЁртэхэш  \eqref{1.1}.  ╩ы■ўхтє■ Ёюы№ яЁш ¤Єюь сєфєЄ шуЁрЄ№
 Ёхчєы№ЄрЄ ЁрсюЄ√ \cite{MirzSh} ю Ёхуєы ЁшчрЎшш єЁртэхэш  \eqref{1.1}  ш Ёхчєы№ЄрЄ√, яюыєўхээ√х эрьш т яЁхф√фє∙хь ярЁруЁрЇх.
┬√яш°хь юЄфхы№эю эршсюыхх юс∙шщ тшф (эхёрьюёюяЁ цхээюую) фшЇЇхЁхэЎшры№эюую т√Ёрцхэш  ўхЄэюую яюЁ фър ё ъюЄюЁ√ь сєфхь ЁрсюЄрЄ№:
\begin{equation}\label{eq:deg}
\tau(y)=\sum_{k,\,s=0}^m(\tau_{k,\,s}(x)y^{(m-k)}(x))^{m-s}
\end{equation}
ш сєфхь яюырурЄ№, ўЄю ъюьяыхъёэючэрўэ√х ЇєэъЎшш   $\tau_{k,\,s}$
яюфўшэхэ√  єёыютш ь \eqref{eq:cond2},  уфх  $\tau_0:=\tau_{0,\,0}$.
%\begin{equation}\label{eq:cond1}
%\frac1{\sqrt{|\tau_0|}},\ \frac1{\sqrt{|\tau_0|}}\tau_{k,\,s}^{(-l)}\in L^2[0,1], \quad 0\le k,\,s\le m,
%\end{equation}
%уфх $l=l_{ks}=\min\{k,\,s\}$, р яюф $f^{(-l)}$ ь√ яюэшьрхь $l$ Ёрч яЁюшэЄхуЁшЁютрээє■ ЇєэъЎш■ $f$.
╬ЄьхЄшь ёЁрчє, ўЄю
т√сюЁ ъюэёЄрэЄ шэЄхуЁшЁютрэш  т \eqref{eq:cond2} эх тыш хЄ эр т√яюыэхэшх ¤Єшї  єёыютшщ.

 ┴єфхь уютюЁшЄ№, ўЄю фшЇЇхЁхэЎшры№эюх
т√Ёрцхэшх \textit{яЁштхфхэю ъ эюЁьры№эющ ЇюЁьх}, хёыш юэю шьххЄ тшф
\begin{multline}\label{eq:den}
l(y)=\sum_{k=0}^m(-1)^{m-k}(\tau_{k}(x)y^{(m-k)}(x))^{(m-k)}+\\
+\i\sum_{k=0}^{m-1}(-1)^{m-k-1}\left[(\sigma_k(x)y^{(m-k-1)})^{(m-k)}+(\sigma_k(x)y^{(m-k)})^{(m-k-1)}\right].
\end{multline}
┬шфэю, ўЄю \eqref{eq:deg} ётюфшЄё  ъ \eqref{eq:den}, хёыш тёх $\tau_{k,\,s}\equiv0$ яЁш $|k-s|\ge2$, р
$\tau_{k,\,k+1}\equiv\tau_{k+1,\,k}$. ╙ёыютш  \eqref{eq:cond2} фы  фшЇЇхЁхэЎшры№э√ї т√Ёрцхэшщ тшфр \eqref{eq:den}
яЁшэшьр■Є тшф
\begin{equation}\label{eq:cond1}
\frac1{\sqrt{|\tau_0|}},\ \frac1{\sqrt{|\tau_0|}}\tau_{k}^{(-k)},\ \frac1{\sqrt{|\tau_0|}}\sigma_{k}^{(-k)}\in
L^2[0,1].
\end{equation}
╬ЄьхЄшь, ўЄю т ЁрсюЄх \cite{MirzSh} єёыютш  \eqref{eq:cond2} ш \eqref{eq:cond1} яЁхфяюырурышё№ т√яюыэхээ√ьш ыш°№
ыюъры№эю эр шэЄхЁтрых $(0,1)$. ▌Єю яючтюы хЄ ъюЁЁхъЄэю юяЁхфхы Є№ фшЇЇхЁхэЎшры№э√х т√Ёрцхэш , ёЄЁюшЄ№ ьшэшьры№э√щ ш
ьръёшьры№э√щ юяхЁрЄюЁ√ ш Є.ф., эю эх фюёЄрЄюўэю фы  эр°шї Ўхыхщ. ╤ыхфє■∙хх єЄтхЁцфхэшх яютЄюЁ хЄ ЄхюЁхьє 1 ЁрсюЄ√
\cite{MirzSh}. ╠√ яЁштхфхь хую ё  яюфЁюсэ√ь фюърчрЄхы№ёЄтюь.
\begin{Proposition}
─ы  ы■сюую фшЇЇхЁхэЎшры№эюую т√Ёрцхэш  тшфр \eqref{eq:deg} эрщфхЄё  ёютярфр■∙хх ё эшь т√Ёрцхэшх тшфр \eqref{eq:den}
(яюф ёютярфхэшхь ь√ яюэшьрхь ЁртхэёЄтю т яЁюёЄЁрэёЄтх юсюс∙хээ√ї ЇєэъЎшщ $l(y)=\tau(y)$ фы  ы■сющ схёъюэхўэю уырфъющ
ЇєэъЎшш $y$ ё эюёшЄхыхь т $(0,1)$).
\end{Proposition}
\begin{proof}
─юуютюЁшьё  эрч√трЄ№ шэфхъёюь ёырурхьюую $(\tau_{k,\,s}y^{(m-k)})^{(m-s)}$ ўшёыю $k-s$. ╨рёёьюЄЁшь тэрўрых
(хфшэёЄтхээюх) ёырурхьюх $(\tau_{m,\,0}y)^{(n)}$ ё шэфхъёюь $m$. ┬ ёшыє \eqref{eq:cond2}
$|\tau_0|^{-1/2}\tau_{m,\,0}\in L_2[0,1]$. ╟рьхэшь т фшЇЇхЁхэЎшры№эюь т√Ёрцхэшш $\tau(y)$ ¤Єю ёырурхьюх ёєььющ
$$
(\tau_{m,\,0}y)^{(n)}=(\tau_{m,\,0}'y)^{(m-1)}+(\tau_{m,\,0}y')^{(m-1)}
$$
фтєї ёырурхь√ї ё шэфхъёрьш $m-1$ ш $m-2$. ╧Ёш ¤Єюь єёыютш  \eqref{eq:cond2} т яЁхюсЁрчютрээюь фшЇЇхЁхэЎшры№эюь
т√Ёрцхэшш Єръцх т√яюыэхэ√, яюёъюы№ъє фы  ¤Єшї ёырурхь√ї $l=1$, р $|\tau_0|^{-1/2}(\tau_{m,\,0}')^{(-1)}\in L_2[0,1]$,
$|\tau_0|^{-1/2}(\tau_{m,\,0})^{(-1)}\in L_2[0,1]$ (ёЄхяхэ№ $(-1)$ ючэрўрхЄ чфхё№ тч Єшх яхЁтююсЁрчэющ). ┬ Ёхчєы№ЄрЄх
яЁшфхь ъ фшЇЇхЁхэЎшры№эюьє т√Ёрцхэш■, т ъюЄюЁюь ёырурхь√х ё шэфхъёюь $m$ юЄёєЄёЄтє■Є, р єёыютш  \eqref{eq:cond2}
т√яюыэхэ√. ╧юы№чє ё№ яЁхюсЁрчютрэшхь
$$
(\tau_{k,\,s}y^{(m-k)})^{(m-s)}=(\tau_{k,\,s}'y^{(m-k)})^{(m-s-1)}+(\tau_{k,\,s}y^{(m-k+1)})^{(m-s-1)}
$$
Єюўэю Єръцх яюёЄєяшь ё фтєь  ёырурхь√ьш шэфхъёр $m-1$ ш Є.ф., яюър эх єэшўЄюцшь тёх т√Ёрцхэш  ё шэфхъёрьш $k-s\ge2$.
╦хуъю тшфхЄ№, ўЄю єёыютш  \eqref{eq:cond1} яЁш ¤Єюь ёюїЁрэ Єё . ╟рЄхь, яюы№чє ё№ ЄюцфхёЄтюь
$$
(\tau_{k,\,s}y^{(m-k)})^{(m-s)}=(\tau_{k,\,s}y^{(m-k-1)})^{(m-s+1)}-(\tau'_{k,\,s}y^{(m-k-1)})^{(m-s)}
$$
(яЁш ¤Єюь яЁхюсЁрчютрэшш шэфхъё т√Ёрцхэш  єтхышўштрхЄё ), ышътшфшЁєхь тёх ёырурхь√х ё шэфхъёрьш $k-s\le-2$. ┬
Ёхчєы№ЄрЄх юёЄрэєЄё  т√Ёрцхэш  ё шэфхъёрьш $-1$, $0$ ш $1$. ═ръюэхЎ, ЄюцфхёЄтю
\begin{multline*}
(\tau_{k+1,\,k}y^{(m-k-1)})^{(m-k)}+(\tau_{k,\,k+1}y^{(m-k)})^{(n-k-1)}=
\left(\frac12(\tau_{k+1,\,k}+\tau_{k,\,k+1})y^{(m-k-1)}\right)^{(m-k)}+\\+
\left(\frac12(\tau_{k+1,\,k}+\tau_{k,\,k+1})y^{(m-k)}\right)^{(m-k-1)}-\left(\frac12(\tau_{k+1,\,k}+\tau_{k,\,k+1})'y^{(m-k-1)}\right)^{(m-k-1)},
\end{multline*}
уфх $k=0,\dots,\,m-1$, яЁштюфшЄ фшЇЇхЁхэЎшры№эюх т√Ёрцхэшх ъ тшфє \eqref{eq:den}.
\end{proof}
\begin{Note}
═хёыюцэю тшфхЄ№, ўЄю ёфхырээ√х яЁхюсЁрчютрэш  эх ьхэ ■Є ёЄрЁ°хую ёырурхьюую $(\tau_0(x)y^{(m)})^{(m)}$.
\end{Note}
╚Єръ, фрыхх сєфхь ЁрсюЄрЄ№ ё т√Ёрцхэш ьш тшфр \eqref{eq:den}. ╧юърцхь, ъръ юЄ єЁртэхэш  $l(y)=\la^{n}\varrho y$ яхЁхщЄш
ъ ёшёЄхьх тшфр \eqref{eq:gen} (ь√ яЁхфяюырурхь $\varrho\in L_1[0,1]$). ╬яЁхфхышь ЇєэъЎшш $\cT_k=\tau_k^{(-k)}$ ш
$\cS_k=\sigma_k^{(-k)}$ (т√сюЁ ъюэёЄрэЄ шэЄхуЁшЁютрэш  яЁюшчтюыхэ), $\varphi_k=\cT_k+\i\cS_{k-1}$,
$\psi_k=\cT_k-\i\cS_{k-1}$. ╥хяхЁ№ ттхфхь ьрЄЁшЎє $\F(x)=(f_{j,\,k}(x))_{j,\,k=1}^{n}$. ╧Ёхцфх тёхую, яюыюцшь
\begin{gather*}
f_{j,\,k}(x)=0,\quad\text{хёыш}\ \ \begin{cases}1\le j\le m-1\\k\ne j+1\end{cases};\quad\begin{cases}j=m\\k\ge
m+1\end{cases};\quad\begin{cases}m+1\le j\le n\\ k\ge m+2\\k\ne j+1,\end{cases}
\end{gather*}
\begin{gather*}
f_{j,\,j+1}=\begin{cases}1\ &\text{яЁш}\ j\ne m,\\
\tau_0^{-1}\ &\text{яЁш}\ j=m,\end{cases}
\end{gather*}
Єръ ўЄю ьрЄЁшЎр $\F$ яЁшюсЁхЄхЄ тшф
$$
\F=\begin{pmatrix}0&1&0&0&\ldots&0&0&\ldots&0&0\\
0&0&1&0&\ldots&0&0&\ldots&0&0\\\hdotsfor{10}\\f_{m,\,1}&f_{m,\,2}&\hdotsfor{2}&f_{m,\,m}&\tau_0^{-1}&0&\ldots&0&0\\
f_{m+1,\,1}&f_{m+1,\,2}&\hdotsfor{3}&f_{m+1,\,m+1}&1&\ldots&0&0\\
\hdotsfor{10}\\f_{n-1,\,1}&f_{n-1,\,2}&\hdotsfor{3}&f_{n-1,\,m+1}&0&\ldots&0&1\\f_{n,\,1}&f_{n,\,2}&f_{n,\,3}&\hdotsfor{2}&f_{n,\,m+1}&0&\ldots&0&0
\end{pmatrix}
$$
╬яЁхфхышь ёЄЁюъє $f_{m,\,k}=(-1)^{m-k}\varphi_{m+1-k}\tau_0^{-1}$, $1\le k\le m$, ш ёЄюысхЎ
$f_{j,\,m+1}=-\psi_{j-m}\tau_0^{-1}$, $m+1\le j\le n$. ╬ёЄрт°шхё  ¤ыхьхэЄ√ ьрЄЁшЎ√ $F$ юяЁхфхышь ЇюЁьєыющ
$$
f_{m+k,\,m-j}=(-1)^{j+1}\varphi_{j+1}\psi_k\tau_0^{-1}+\chi_{j+k<m}(-1)^jC_{j+k+1}^k\left[\cT_{j+k+1}+\i\frac{j-k+1}{j+k+1}\cS_{j+k}\right],
$$
уфх $1\le k\le m$, $0\le j\le m-1$, $C_{j+k+1}^k$ --- сшэюьшры№э√щ ъю¤ЇЇшЎшхэЄ, р ўшёыю $\chi_{j+k<m}$ Ёртэю $1$ яЁш
$j+k<m$ ш $0$ шэрўх. ╬ЄьхЄшь трцэ√щ ЇръЄ, ъюЄюЁ√щ т√ЄхърхЄ эхяюёЁхфёЄтхээю шч єёыютшщ \eqref{eq:cond1}: \textit{ЇєэъЎшш
$f_{j,\,k}(x)$ ёєььшЁєхь√ эр $[0,1]$.}

╥хяхЁ№ ттхфхь ътрчшяЁюшчтюфэ√х\footnote{┬ Ёрчышўэ√ї шёЄюўэшърї (ёь., эряЁшьхЁ, \cite{AEZ}, \cite{E}), т Єюь ўшёых т
ЁрсюЄрї ртЄюЁют \cite{SavSh99}, \cite{SavSh03}, шёяюы№чє■Єё  ш фЁєушх юяЁхфхыхэш  ътрчшяЁюшчтюфэ√ї. ═ряЁшьхЁ, фы 
ёыєўр  юяхЁрЄюЁр ╪ЄєЁьр--╦шєтшыы , ъюуфр $n=2$ ш эхюсїюфшью юяЁхфхышЄ№ Єюы№ъю яхЁтє■ ътрчшяЁюшчтюфэє■, єфюсэю юЄсЁюёшЄ№
т юяЁхфхыхэшш ьэюцшЄхы№ $1/f_{1,\,2}$. ╥ръшх Ёрчышўш  ёыхфєхЄ єўшЄ√трЄ№ яЁш ёЁртэхэшш Ёхчєы№ЄрЄют юс рёшьяЄюЄшъх
Ёх°хэшщ.} ЇєэъЎшш $y$, ъюЄюЁ√х ь√ сєфхь юсючэрўрЄ№ $y^{[j]}(x)$
\begin{equation}\label{eq:quasider}
y^{[0]}=y,\quad y^{[k]}=\frac1{f_{k,\,k+1}(x)}\left[(y^{[k-1]})'-\sum_{j=1}^{k}f_{k,\,j}(x)y^{[j-1]}\right],\quad 1\le
k\le n-1.
\end{equation}
╤ЇюЁьєышЁєхь ЄхюЁхьє 2 ЁрсюЄ√ \cite{MirzSh} (хх фюърчрЄхы№ёЄтю ЄЁхсєхЄ ёхЁ№хчэющ  Єхїэшўхёъющ ЁрсюЄ√ ш чфхё№ эх
яЁштюфшЄё ).
\begin{Proposition}
┼ёыш $l(y)$ --- фшЇЇхЁхэЎшры№эюх т√Ёрцхэшх тшфр \eqref{eq:den}, Єю т эют√ї юсючэрўхэш ї юэю шьххЄ тшф
\begin{equation}\label{eq:dequasi}
l(y)=(y^{[n-1]})'-\sum_{j=1}^{n}f_{n,\, j}(x)y^{[j-1]}.
\end{equation}
\end{Proposition}
═ряюьэшь, ўЄю эр°р Ўхы№ ёюёЄюшЄ т шёёыхфютрэшш єЁртэхэш  $l(y)=\la^{n}\varrho(x)y$, уфх $\la^{n}\in\bC$ ---
ёяхъЄЁры№э√щ ярЁрьхЄЁ. ╠√ сєфхь яЁхфяюырурЄ№ $\varrho(x)\in L_1[0,1]$ (яючфэхх ь√ эрыюцшь фюяюыэшЄхы№э√х юуЁрэшўхэш ).
╟рьхЄшь, ўЄю ЁртхэёЄтр \eqref{eq:quasider} ш \eqref{eq:dequasi} єцх яючтюы ■Є чряшёрЄ№ єЁртэхэшх
$l(y)=\la^{n}\varrho(x)y$ т тшфх ёшёЄхь√ $n$ фшЇЇхЁхэЎшры№э√ї єЁртэхэшщ яхЁтюую яюЁ фър юЄэюёшЄхы№эю тхъЄюЁ--ЇєэъЎшш,
ёюёЄртыхээющ шч ътрчшяЁюшчтюфэ√ї $y^{[j]}(x)$, $0\le j\le n-1$. ╠рЄЁшЎхщ ¤Єющ ёшёЄхь√, юўхтшфэю,  ты хЄё  ьрЄЁшЎр
$\F(x)$ ё хфшэёЄтхээ√ь шчьхэхэшхь --- ъ хх ¤ыхьхэЄє $f_{n,\,1}(x)$ фюсрты хЄё  ёырурхьюх $(-1)^m\la^{n}\varrho(x)$.
╙фюсэхх, юфэръю, юёє∙хёЄтшЄ№ яхЁхїюф ъ ёшёЄхьх, яюыюцшт $\u(x)=(u_j(x))_{j=1}^{n}$, уфх
$u_j(x):=\la^{1-j}y^{[j-1]}(x)$.
\begin{Proposition}
┼ёыш $l(y)$ --- фшЇЇхЁхэЎшры№эюх т√Ёрцхэшх тшфр \eqref{eq:den}, Єю єЁртэхэшх $l(y)=\la^{n}\varrho(x)y$ т эр°шї
юсючэрўхэш ї яЁшэшьрхЄ тшф
\begin{equation}\label{eq:systF1}
\u'=\F(x,\la)\u,\qquad\text{уфх}\ \
\F(x,\la)=\Big(f_{j,\,k}(x)\la^{k-j}+(-1)^m\la\varrho(x)\delta_{j}^{n}\delta_k^1\Big)_{j,\,k=1}^{n}.
\end{equation}
\end{Proposition}
\begin{proof}
╧юыюцшь $\v(x)=(y^{[j-1]}(x))_{j=1}^{n}$. ╥юуфр $\v'=(\F(x)+(-1)^m\la^{n}\E_{n,\,1})\v$, уфх ўхЁхч $\E_{\al,\,\beta}$
ь√ юсючэрўрхь ьрЄЁшЎє $E_{\al,\,\beta}=\big(\delta_{j}^{\al}\delta_k^\beta\big)_{j,\,k=1}^{n}$, $\delta^k_j$ --- ёшьтюы
╩ЁюэхъхЁр. ╤юуырёэю эр°шь юсючэрўхэш ь, $\v=\diag\{1,\la,\dots,\la^{n-1}\}\u$, р чэрўшЄ
$$
\u'=\diag\{1,\la^{-1},\dots,\la^{1-n}\}(\F(x)+(-1)^m\la^{n}\E_{n,\,1})\diag\{1,\la,\dots,\la^{n-1}\}\u.
$$
╧хЁхьэюцшт ¤Єш ЄЁш ьрЄЁшЎ√, яЁшїюфшь ъ \eqref{eq:systF1}.
\end{proof}
╥хяхЁ№ чрьхЄшь, ўЄю ьрЄЁшЎр $\F(x,\la)$ т \eqref{eq:systF1} Ёрёъырф√трхЄё  т ёєььє
\begin{gather*}
\F(x,\la)=\la\F_1(x)+\F_0(x)+\la^{-1}\F_{-1}(x)+\dots+\la^{1-n}\F_{1-n}(x),\quad\text{уфх}\\
\F_1(x)=\sum_{j=1}^{m-1}\E_{j,\,j+1}+\tau_0^{-1}(x)\E_{m,\,m+1}+\sum_{j=m+1}^{n-1}\E_{j,\,j+1}+(-1)^m\varrho(x)\E_{n,\,1},\\
\F_0(x)=f_{m,\,m}(x)\E_{m,\,m}+f_{m+1,\,m+1}(x)\E_{m+1,\,m+1},\\
\F_{-k}(x)=\sum_{j=k+1}^{n}f_{j,\,j-k}\E_{j,\,j-k},\quad 1\le k\le n-1.
\end{gather*}
╥ръшь юсЁрчюь, ёшёЄхьр \eqref{eq:systF1} шьххЄ тшф \eqref{eq:gen} ш, ёыхфє  єЄтхЁцфхэш■ \ref{pr:gen}, ь√ ётхфхь хх ъ
ЇюЁьх \eqref{eq:main}. ╧Ёш ¤Єюь юфэръю, эрь яЁшфхЄё  эрыюцшЄ№ фюяюыэшЄхы№э√х єёыютш  эр ЇєэъЎшш $\tau_0$ ш $\varrho$.
─рыхх ўхЁхч $\omega_k$, $0\le k\le n-1$, юсючэрўрхь ъюЁэш ёЄхяхэш $n$ шч $(-1)^m$, чрэєьхЁютрээ√х т яюЁ фъх
$$
\omega_k=\begin{cases}\epsilon_l^{\frac{2k+1}{2}},\quad&\text{хёыш ўшёыю }m\text{ эхўхЄэю,}\\
\epsilon_l^k,\quad&\text{хёыш }m\text{ ўхЄэю,}\end{cases}\qquad\epsilon_l=e^{\frac{i\pi l}m}.
$$
╬сючэрўшь Єръцх $\rho(x)=\varrho^{\frac1n}(x)\tau_0^{-\frac1n}(x)$.
\begin{Theorem}\label{tm:ho}
╧єёЄ№ фшЇЇхЁхэЎшры№эюх т√Ёрцхэшх $\tau$ шьххЄ тшф \eqref{eq:deg}, т√яюыэхэ√ єёыютш  \eqref{eq:cond1}, р ЇєэъЎш 
$\tau_0$ рсёюы■Єэю эхяЁхЁ√тэр ш яюыюцшЄхы№эр эр $[0,1]$. ╧єёЄ№ ЇєэъЎш  $\varrho$ Єръцх рсёюы■Єэю эхяЁхЁ√тэр ш
яюыюцшЄхы№эр. ╥юуфр єЁртэхэшх $\tau(y)=\la^{n}\varrho(x)y$ ьюцэю ётхёЄш ъ ёшёЄхьх тшфр \eqref{eq:main}, ё ьрЄЁшЎхщ
$\B=\diag\{\omega_0,\,\omega_1,\dots,\,\omega_{n-1}\}$. ╧Ёш ¤Єюь тёх ¤ыхьхэЄ√ ьрЄЁшЎ $\A(x)$ ш $\C(x,\la)$ ёєььшЁєхь√
эр $[0,1]$, р $\|\C(\cdot,\la)\|_{L_1}=O(|\la|^{-1})$, $\la\to\infty$.
\end{Theorem}
\begin{proof}
┴юы№°р  ўрёЄ№ ЁрсюЄ√ єцх яЁюфхырэр --- ь√ яюърчрыш, ъръ ётхёЄш єЁртэхэшх $\tau(y)=\la^{n}\varrho(x)y$ ъ ёшёЄхьх
\eqref{eq:systF1}. ╬ёЄрхЄё  фшруюэрышчютрЄ№ ьрЄЁшЎє $\F_1(x)$. ╦хуъю тшфхЄ№, ўЄю їрЁръЄхЁшёЄшўхёъшщ ьэюуюўыхэ ¤Єющ
ьрЄЁшЎ√ Ёртхэ $\chi(s)=s^{n}+(-1)^{m+1}\varrho(x)\tau^{-1}(x)$, р чэрўшЄ хую ъюЁэ ьш  ты ■Єё  ўшёыр
$s_k=\omega_k\rho(x)$. ╟рьхЄшь, ўЄю ЇєэъЎш  $\rho(x)$ яюыюцшЄхы№эр, р чэрўшЄ фы  ы■сюую $x\in[0,1]$ ўшёыр $s_k$ яюярЁэю
Ёрчышўэ√. ▌Єю ючэрўрхЄ, ўЄю эрщфхЄё  ьрЄЁшЎр яхЁхїюфр $\W(x)$ Єрър , ўЄю
$\W^{-1}(x)\F_1(x)\W(x)=\rho(x)\diag\{\omega_0,\,\omega_1,\dots,\,\omega_{n-1}\}$. ╬яєёър  т√ўшёыхэш , эрщфхь
\begin{gather}\label{eq:matrW}
\W(x)=(w_{jk}(x))_{j,\,k=1}^{n},\qquad
w_{jk}(x)=(\omega_{k-1}\rho(x))^{j-1}\cdot\begin{cases}1,\quad&j\le m\\
\tau_0,\quad&j>m,\end{cases}\\\notag \W^{-1}(x)=(\wt w_{jk}(x))_{j,\,k=1}^{n},\qquad \wt
w_{jk}(x)=\frac1{n}(\omega_{j-1}\rho(x))^{1-k}\cdot\begin{cases}1,\quad&k\le m\\
\tau^{-1}_0,\quad&k>m.\end{cases}
\end{gather}
╥ръшь юсЁрчюь, чрьхэр $\y=\W^{-1}\u$ яЁштюфшЄ ёшёЄхьє \eqref{eq:systF1} ъ тшфє
\begin{equation}\label{eq:hosyst}
\y'=\la\rho(x)\B\y+\A(x)\y+\C(x,\la)\y,\qquad\text{уфх}\ \ \B=\diag\{\omega_0,\,\omega_1,\dots,\,\omega_{n-1}\}.
\end{equation}
╧Ёш ¤Єюь,
$$
\A(x)=-\W^{-1}(x)\W'(x)+\W^{-1}(x)\F_0(x)\W(x),\qquad \text{р} \ \
\C(x,\la)=\W^{-1}(x)\sum_{k=1}^{n-1}\la^{-k}\F_{-k}(x)\W(x).
$$
╠√ эх сєфхь яЁштюфшЄ№ чфхё№ ЇюЁьєы√, ёт ч√тр■∙шх ¤ыхьхэЄ√ ьрЄЁшЎ√ $\C$ ё ЇєэъЎш ьш $\tau_k$ ш $\sigma_k$, ттшфє шї
уЁюьючфъюёЄш. ╬ЄьхЄшь ыш°№, ўЄю тёх ¤ыхьхэЄ√ ¤Єющ ьрЄЁшЎ√ яЁшэрфыхцрЄ яЁюёЄЁрэёЄтє $L_1[0,1]$ (¤Єю ёЁрчє цх ёыхфєхЄ шч
рсёюы■Єэющ эхяЁхЁ√тэюёЄш ЇєэъЎшщ $w(x,\la)$ ш $\wt w(x,\la)$ ш ёєььшЁєхьюёЄш ЇєэъЎшщ $f_{j,\,k}(x)$). ╩Ёюьх Єюую,
$\C(x,\la)=O(|\la|^{-1})$ яЁш $\la\to\infty$. ╠рЄЁшЎє $\A(x)=(a_{j,\,k}(x))_{j,\,k=1}^{n}$ ь√, эряЁюЄшт, т√яш°хь  тэю
\begin{equation}\label{eq:matrA}
a_{j,\,k}(x)=\frac1{n}\begin{cases}\dfrac{1}{1-\epsilon_{k-j}}\left(\dfrac{\varrho'(x)}{\varrho(x)}-\epsilon_{k-j}^m\dfrac{\tau_0'(x)}{\tau_0(x)}\right)
+\dfrac{\epsilon_{k-j}^m}{\tau_0(x)}\left(
\varphi_1(x)\epsilon_{j-k}-\psi_1(x)\right),\ &j\ne k,\\\,\\
\dfrac{1-n}2(\ln\varrho(x))'-\dfrac12(\ln\tau_0(x))'+\dfrac{2\i\sigma_0(x)}{\tau_0(x)},\ &j=k.\end{cases}
\end{equation}
\end{proof}

╥хяхЁ№ ь√ яюыєўшь рёшьяЄюЄшўхёъшх ЇюЁьєы√ фы  Ёх°хэшщ фшЇЇхЁхэЎшры№эюую єЁртэхэш  $l(y)=\la^{n}\varrho(x)y$ яюЁ фър
$n=2m$ тшфр ёю ёяхъЄЁры№э√ь ярЁрьхЄЁюь т яЁртющ ўрёЄш. ╘єэъЎш■ $\varrho(x)$ сєфхь яЁхфяюырурЄ№ рсёюы■Єэю эхяЁхЁ√тэющ ш
яюыюцшЄхы№эющ эр $[0,1]$. ┴єфхь ёўшЄрЄ№, ўЄю фшЇЇхЁхэЎшры№эюх т√Ёрцхэшх $l(y)$ єцх яЁштхфхэю ъ эюЁьры№эюьє тшфє
\eqref{eq:den}. ─рыхх схч эряюьшэрэшщ яюы№чєхьё  юсючэрўхэш ьш, ттхфхээ√ьш т√°х. ╤юуырёэю ЄхюЁхьх \ref{tm:ho}, ьрЄЁшЎр
$\B$ шЄюуютющ ёшёЄхь√ Ёртэр $\B=\diag\{\omega_0,\dots,\,\omega_{n-1}\}$. ┬эрўрых юяш°хь ёшёЄхьє ёхъЄюЁют $\Gamma_k$,
уЁрэшЎ√ ъюЄюЁ√ї юяЁхфхы ■Єё  єЁртэхэш ьш $\Re(\omega_j\la)=\Re(\omega_k\la)$, $0\le j\ne k\le n-1$. ╙ўшЄ√тр   тэ√щ тшф
ўшёхы $\omega_j$, ¤Єш єЁртэхэш  ыхуъю Ёх°р■Єё  (Ёх°хэшх ЄЁшуюэюьхЄЁшўхёъшї єЁртэхэш  ь√ чфхё№ юяєёърхь). ╩юэхўэ√щ
Ёхчєы№ЄрЄ ёыхфє■∙шщ: яЁш $n=2$ яыюёъюёЄ№ ЁрчсштрхЄё  эр фтр ёхъЄюЁр --- тхЁїэ■■ ш эшцэ■■ яюыєяыюёъюёЄ№; яЁш $n>2$
яыюёъюёЄ№ ЁрчсштрхЄё  эр $2n$ ёхъЄюЁют тшфр
$$
\Gamma_k=\left\{\la:\frac{\pi(k-1)}{n}\le\arg\la\le\frac{\pi k}{n}\right\},\quad k=1,\dots,\,2n.
$$
╩ръ ш Ёрэхх, Ёрё°шЁшь ърцф√щ ёхъЄюЁ $\Gamma_k$, ёфтшэєт хую эр $r$ тфюы№ сшёёхъЄЁшё√. ╧юыєўхээ√х ёхъЄюЁ√ юсючэрўрхь
$\wt\Gamma_k$.
\begin{Theorem}\label{tm:final}
╧єёЄ№ $l(y)$ --- фшЇЇхЁхэЎшры№эюх т√Ёрцхэшх тшфр \eqref{eq:den}, фы  ъю¤ЇЇшЎшхэЄют ъюЄюЁюую т√яюыэхэ√ єёыютш 
\eqref{eq:cond1}. ╧єёЄ№ ъ Єюьє цх ЇєэъЎшш $\tau_0(x)$ ш $\varrho(x)$ рсёюы■Єэю эхяЁхЁ√тэ√ ш яюыюцшЄхы№э√ эр $[0,1]$.
╧єёЄ№ $\wt\Gamma_\kappa$ --- юфшэ шч ёхъЄюЁют, юяЁхфхыхээ√ї т√°х. ╥юуфр т ¤Єюь ёхъЄюЁх єЁртэхэшх
$l(y)=\la^{n}\varrho(x)y$ шьххЄ ёшёЄхьє ЇєэфрьхэЄры№э√ї Ёх°хэшщ $y_k(x,\la)$, $k=1,\dots,\,n$, фы  ърцфюую шч ъюЄюЁ√ї
ёяЁртхфыштю рёшьяЄюЄшўхёъюх яЁхфёЄртыхэшх
\begin{equation}\label{eq:final1}
y_k(x)=e^{\omega_{k-1}\la
p(x)}\left[\varrho^{\frac{1-n}{2n}}(x)\tau_0^{-\frac{1}{2n}}(x)\exp\left\{\frac{2\i}n\int_0^x\frac{\sigma_0(t)}{\tau_0(t)}dt\right\}+\zeta_{1\,k}(x,\la)\right].
\end{equation}
╧Ёш ¤Єюь фы  ътрчшяЁюшчтюфэ√ї яюЁ фър $j=1,\dots,\,m-1$ ¤Єшї ЇєэъЎшщ Єръцх ёяЁртхфышт√ рёшьяЄюЄшўхёъшх яЁхфёЄртыхэш 
\begin{equation}\label{eq:final2}
y_k^{[j]}(x)=\la^{j}e^{\omega_{k-1}\la
p(x)}\left[(\omega_{k-1})^j\varrho^{\frac{2j-n+1}{2n}}(x)\tau_0^{-\frac{1+2j}{2n}}(x)\exp\left\{\frac{2\i}n\int_0^x\frac{\sigma_0(t)}{\tau_0(t)}dt\right\}
+\zeta_{j\,k}(x,\la)\right],
\end{equation}
р фы  ътрчшяЁюшчтюфэ√ї яюЁ фър $j=m,\dots,\,n-1$ --- яЁхфёЄртыхэш 
\begin{equation}\label{eq:final3}
y_k^{[j]}(x)=\la^{j}e^{\omega_{k-1}\la
p(x)}\left[(\omega_{k-1})^j\varrho^{\frac{2j-n+1}{2n}}(x)\tau_0^{1-\frac{1+2j}{2n}}(x)\exp\left\{\frac{2\i}n\int_0^x\frac{\sigma_0(t)}{\tau_0(t)}dt\right\}+
\zeta_{j\,k}(x,\la)\right].
\end{equation}
╬ёЄрЄъш т ¤Єшї яЁхфёЄртыхэш ї фюяєёър■Є яЁш $\wt\Gamma_\kappa\ni\la\to\infty$ юЎхэъє
\begin{gather}\label{eq:final4}
\max_{j,\,k}|\zeta_{j,\,k}(x,\la)|\le C(\Upsilon(x,\la)+\Upsilon_\mu(\la)+|\la|^{-1}),\\
\max_{j,\,k}\|\zeta_{j,\,k}(x,\la)\|_{L_\mu}\le C(\Upsilon_\mu(\la)+|\la|^{-1}),\notag
\end{gather}
уфх $1/\mu+1/\mu'=1$, р $\mu'$ --- ьръёшьры№эр  юс∙р  ёЄхяхэ№ ёєььшЁєхьюёЄш ЇєэъЎшщ $\varrho'(x)$, $\tau'_0(x)$,
$\tau_1^{(-1)}(x)$ ш $\sigma_0(x)$. ╘єэъЎшш $\Upsilon_\mu(\la)$ ш $\Upsilon(x,\la)$ юяЁхфхы ■Єё  ЇюЁьєырьш
\eqref{eq:defUps} ш \eqref{eq:defups}, уфх $q_{jl}(x)=a_{jl}(x)$, ьрЄЁшЎр $\A(x)$ юяЁхфхыхэр т \eqref{eq:matrA}.
\end{Theorem}
\begin{proof}
╧юыюцшь $\u(x)=(u_j(x))_{j=1}^{n}$, $u_j(x)=\la^{1-j}y^{[j-1]}(x)$, ш $\y(x)=\W^{-1}(x)\u(x)$, уфх ьрЄЁшЎ√ $\W(x)$ ш
$\W^{-1}(x)$ юяЁхфхыхэ√ т \eqref{eq:matrW}. ╥юуфр фы  тхъЄюЁр $\y(x)$ т√яюыэхэр ёшёЄхьр \eqref{eq:hosyst} ё ьрЄЁшЎхщ
$\A(x)$, юяЁхфхыхээющ т \eqref{eq:matrA}. ╧Ёшьхэшь ЄхюЁхьє \ref{tm:main} фы  яюшёър рёшьяЄюЄшўхёъюую яЁхфёЄртыхэш 
ЇєэфрьхэЄры№эющ ьрЄЁшЎ√ $\Y(x,\la)$ ¤Єющ ёшёЄхь√. ┬эрўрых эрщфхь уыртэ√щ ўыхэ Ёрчыюцхэш  --- ьрЄЁшЎє $\Y^0(x,\la)$.
╟рьхЄшь, ўЄю тёх ўшёыр $b_j=\omega_{j-1}$, $j=1,\dots,\,n$, яюярЁэю Ёрчышўэ√, р чэрўшЄ (ёь. чрьхўрэшх \ref{nt:diagM})
ьрЄЁшЎр $\M(x)$ фшруюэры№эр. ┴юыхх Єюую, шч \eqref{eq:matrA} тшфэю, ўЄю тёх ¤ыхьхэЄ√ $a_{jj}$, $1\le j\le n$, Ёртэ√
ьхцфє ёюсющ ш
$$
\int_0^x
a_{jj}(t)\,dt=\frac{1-n}{2n}\ln\varrho(x)-\frac1{2n}\ln\tau_0(x)+\frac{2\i}{n}\int_0^x\frac{\sigma_0(t)}{\tau_0(t)}\,dt,
$$
р чэрўшЄ
$$
\Y^0(x,\la)=\varrho^{\frac{1-n}{2n}}(x)\tau_0^{-\frac1{2n}}(x)e^{\frac{2\i}{n}\int_0^x\frac{\sigma_0(t)}{\tau_0(t)}\,dt}\cE(x,\la),
\qquad\cE(x,\la=\diag\left\{e^{\omega_0\la p(x)},\ldots,\,e^{\omega_{n-1}\la p(x)}\right\}.
$$
╥хяхЁ№ тюёяюы№чєхьё  ЄхюЁхьющ \ref{tm:main} ш эрщфхь $\Y(x,\la)=\Y^0(x,\la)+\cA(x,\la)\cE(x,\la)$, яЁшўхь фы  ¤ыхьхэЄют
ьрЄЁшЎ√ $\cA$ т√яюыэхэ√ юЎхэъш \eqref{eq:matras1} ш \eqref{eq:matras2}. ╬ёЄрхЄё  ёфхырЄ№ юсЁрЄэє■ чрьхэє, яюыюцшт
$\U(x,\la)=\W(x)\Y(x,\la)$. ╧хЁхьэюцшт ьрЄЁшЎ√, яюыєўшь
\begin{gather*}
\U(x,\la)=\U^0(x,\la)+\W(x)\EuScript A(x,\la),\quad\text{уфх}\\
u^0_{jk}(x,\la)=\varrho^{\frac{1-n}{2n}}(x)\tau_0^{-\frac1{2n}}(x)e^{\omega_{k-1}\la
p(x)+\frac{2\i}{n}\int_0^x\frac{\sigma_0(t)}{\tau_0(t)}\,dt}(\omega_{k-1}\rho(x))^{j-1}\cdot\begin{cases}1,\ \ &j\le m\\
\tau_0,\ \ &j>m\end{cases}.
\end{gather*}
╧юфёЄрэютър $j=1$ фрхЄ \eqref{eq:final1}, р яЁш $2\le j\le n$ яюыєўрхь ЁртхэёЄтр \eqref{eq:final2} ш \eqref{eq:final3},
яюёъюы№ъє $y_k^{[j-1]}(x)=\la^{j-1}u_{jk}(x)$. ═хЁртхэёЄтр \eqref{eq:matras1} ш \eqref{eq:matras2}, ёяЁртхфышт√х фы 
¤ыхьхэЄют ьрЄЁшЎ√ $\cA(x,\la)$, тыхъєЄ юЎхэъш \eqref{eq:final4} фы  ЇєэъЎшщ $\zeta_{jk}(x,\la)$
--- ¤ыхьхэЄют ьрЄЁшЎ√ $\W(x)\cA(x,\la)$, яюёъюы№ъє ьрЄЁшЎр $\W(x)$ рсёюы■Єэю
эхяЁхЁ√тэр ш эх чртшёшЄ юЄ $\la$.
\end{proof}
┬√яш°хь юЄфхы№эю Ёхчєы№ЄрЄ фы  ёыєўр  $n=2$. ▌Єю ьюцэю ёфхырЄ№ эхяюёЁхфёЄтхээю, тч т $n=2$ т єЄтхЁцфхэшш ЄхюЁхь√
\ref{tm:final}, эю фы  эруы фэюёЄш ь√ ъЁрЄъю яютЄюЁшь юёэютэ√х ьюьхэЄ√ хх фюърчрЄхы№ёЄтр фы  ¤Єюую ёыєўр . ╩Ёюьх Єюую,
 тэ√щ тшф ьрЄЁшЎ√ $\A(x)$ яючтюышЄ эрь эхёъюы№ъю юёырсшЄ№ єёыютш  ЄхюЁхь√ \ref{tm:final}. ╚Єръ, яЁш $n=2$ шёёыхфєхьюх
єЁртэхэшх шьххЄ тшф
\begin{equation}\label{eq:de2}
l(y)=-(\tau_0(x)y')'+\i(\sigma_0(x)y)'+\i\sigma_0(x)y'+\tau_1(x)y=\la^2\varrho(x)y,
\end{equation}
уфх
\begin{equation}\label{eq:cond3}
\tau_0(x),\, \varrho(x)\in W_1^1[0,1],\quad \tau_0(x)>0,\ \varrho(x)>0,\quad \tau_1(x)=\cT_1'(x),\quad  \cT_1(x),\,
\sigma_0(x)\in L_2[0,1].
\end{equation}
╧хЁхїюф ъ ёшёЄхьх юёє∙хёЄты хЄё  чрьхэющ
$$
u_1(x)=y(x),\qquad u_2(x)=\la^{-1}y^{[1]}(x)=\la^{-1}\tau_0(x)(y'(x)-f_{1,1}(x)y),
$$
уфх ¤ыхьхэЄ√ ьрЄЁшЎ√ $\F(x)$ Ёртэ√
\begin{gather*}
f_{1,1}(x)=\frac{\cT_1(x)+\i\sigma_0(x)}{\tau_0(x)},\qquad f_{1,2}(x)=\frac1{\tau_0(x)},\\
f_{2,1}(x)=-\frac{\cT_1^2(x)+\sigma_0^2(x)}{\tau_0(x)},\qquad f_{2,2}(x)=\frac{-\cT_1(x)+\i\sigma_0(x)}{\tau_0(x)}.
\end{gather*}
╤рьр ёшёЄхьр юЄэюёшЄхы№эю тхъЄюЁр $\u(x)=(u_1(x),u_2(x))^t$ шьххЄ тшф
$$
\u'(x)=\begin{pmatrix}f_{1,1}(x)&\frac{\la}{\tau_0(x)}\\-\la\varrho(x)+\frac{f_{2,2}(x)}{\la}&f_{2,1}(x)\end{pmatrix}\u(x),
$$
р яюёых чрьхэ√
$$
\y=\frac12\begin{pmatrix}1&-\frac{\i}{\sqrt{\varrho\tau_0}}\\1&\frac{\i}{\sqrt{\varrho\tau_0}}\end{pmatrix}\u,\qquad
\u=\begin{pmatrix}1&1\\ \i\sqrt{\varrho\tau_0}&-\i\sqrt{\varrho\tau_0}\end{pmatrix}\y
$$
яЁшюсЁхЄрхЄ тшф (чфхё№ $\rho(x)=\varrho(x)/\tau_0(x)$)
$$
\y'=\la\sqrt{\rho(x)}\begin{pmatrix}\i&0\\0&-\i\end{pmatrix}\y+
\begin{pmatrix}-\frac{(\varrho\tau_0)'}{4\varrho\tau_0}+\frac{\i\sigma_0}{\tau_0}&\frac{(\varrho\tau_0)'}{4\varrho\tau_0}+\frac{\cT_1}{\tau_0}\\
\frac{(\varrho\tau_0)'}{4\varrho\tau_0}+\frac{\cT_1}{\tau_0}&-\frac{(\varrho\tau_0)'}{4\varrho\tau_0}+\frac{\i\sigma_0}{\tau_0}\end{pmatrix}\y
+\frac1{\la}\cdot\frac{\i(\cT_1^2+\sigma_0^2)}{\tau_0^{3/2}\varrho^{1/2}}\begin{pmatrix}1&1\\-1&-1\end{pmatrix}\y.
$$
\begin{Corollary}
╧єёЄ№ $\wt\Gamma_1=\{\la:\Im\la>-h\}$, $\wt\Gamma_2=\{\la:\Im\la<h\}$, уфх $h>0$ яЁюшчтюы№эю. ╧Ёш єёыютшш
\eqref{eq:cond3} єЁртэхэшх \eqref{eq:de2} шьххЄ т $\wt\Gamma_1$ Єръє■ ярЁє ышэхщэю эхчртшёшь√ї Ёх°хэшщ $y_\pm(x,\la)$,
ўЄю яЁш $\wt\Gamma_1\ni\la\to\infty$ ёяЁртхфышт√ рёшьяЄюЄшўхёъшх яЁхфёЄртыхэш 
\begin{gather*}
y_\pm(x,\la)=e^{\pm \i\la
p(x)}\left[\varrho^{-\frac14}(x)\tau_0^{-\frac14}(x)e^{\i\int_0^x\frac{\sigma_0(t)}{\tau_0(t)}dt}+\zeta_{\pm}(x,\la)\right],
\qquad p(x)=\int_0^x\varrho^{\frac12}(t)\tau_0^{-\frac12}(t)\,dt,\\
y_{\pm}^{[1]}(x,\la)=\pm \i\la e^{\pm \i\la
p(x)}\left[\varrho^{\frac14}(x)\tau_0^{\frac14}(x)e^{\i\int_0^x\frac{\sigma_0(t)}{\tau_0(t)}dt}+\zeta_{\pm,\,1}(x,\la)\right],
\end{gather*}
уфх юёЄрЄъш $\zeta_{\pm}$ ш $\zeta_{\pm,\,1}$ яюфўшэхэ√ юЎхэърь
$$
\|\zeta(\cdot,\la)\|_{L_\infty}\le C\left(\Upsilon(\la)+|\la|^{-1}\right).
$$
┼ёыш ЇєэъЎшш $(\varrho(x)\tau_0(x))'$ ш $\cT_1(x)$ ёєььшЁєхь√ т ёЄхяхэш $\mu'\in[1,\infty]$, Єю\footnote{═ряюьэшь, ўЄю
т ы■сюь ёыєўрх $\cT_1(x)\in L_2[0,1]$, Єръ ўЄю єёыютшх $\cT_1\in L_{\mu'}[0,1]$ ёюфхЁцрЄхы№эю Єюы№ъю яЁш $\mu'>2$.}
$$
\|\zeta(\cdot,\la)\|_{L_\mu}\le C\left(\Upsilon_\mu(\la)+|\la|^{-1}\right),\qquad|\zeta(x,\la)|\le
\Upsilon(x,\la)+C\left(\Upsilon_\mu(\la)+|\la|^{-1}\right).
$$
┬ ЄюўэюёЄш Єръюх цх єЄтхЁцфхэшх (ё фЁєующ ярЁющ ЇєэъЎшщ $y_\pm(x,\la)$) ёяЁртхфыштю т $\wt\Gamma_2$.
\end{Corollary}
\begin{Note}\label{nt:gen2}
═хяюёЁхфёЄтхээ√щ тшф ьрЄЁшЎ√ ёшёЄхь√ яючтюы хЄ эхёъюы№ъю юёырсшЄ№ єёыютш  эр ъю¤ЇЇшЎшхэЄ√ фшЇЇхЁхэЎшры№эюую т√Ёрцхэш 
\eqref{eq:de2}. ┬ьхёЄю ЄЁхсютрэш  $\sigma_0,\,\cT_1\in L_2[0,1]$ фюёЄрЄюўэю т√яюыэхэш  єёыютшщ $\sigma_0,\,\cT_1\in
L_1[0,1]$, $\sigma_0^2+\cT_1^2\in L_1[0,1]$.
\end{Note}
╬ЄьхЄшь, ўЄю яюыєўхээ√х рёшьяЄюЄшўхёъшх яЁхфёЄртыхэш   ёютярфр■Є ё Ёхчєы№ЄрЄрьш ЁрсюЄ√ \cite{VlaSh}. ┬
¤Єющ ЁрсюЄх с√ы ЁрёёьюЄЁхэ ёыєўрщ $\mu'=1$, р юЄэюёшЄхы№эю ъю¤ЇЇшЎшхэЄют ёшёЄхь√ с√ыш ёфхырэ√ ёыхфє■∙шх яЁхфяюыюцхэш 
$\varrho,\,\tau_0\in \AC[0,1]$, $\varrho,\,\tau>0$, $\sigma_0\in L_1[0,1]$, $\i\sigma_0+\cT_1\in L_2[0,1]$,
$\sigma_0(\i\sigma_0+\cT_1)\in L_1[0,1]$, $\rho'(\i\sigma_0+\cT_1)\in L_1[0,1]$, $\tau_0'(\i\sigma_0+\cT_1)\in
L_1[0,1]$. ╦хуъю тшфхЄ№, ўЄю
$$
\begin{cases}\sigma_0\in L_1[0,1]\\ \i\sigma_0+\cT_1\in L_2[0,1]\\ \sigma_0(\i\sigma_0+\cT_1)\in
L_1[0,1]\end{cases}\Longleftrightarrow\begin{cases}\sigma_0\in L_1[0,1]\\ \i\sigma_0+\cT_1\in L_2[0,1]\\
\sigma_0^2+\cT_1^2\in L_1[0,1]\end{cases}\Longrightarrow\begin{cases}\sigma_0\in L_1[0,1]\\ \cT_1\in L_1[0,1]\\
\sigma_0^2+\cT_1^2\in L_1[0,1],\end{cases}
$$
Єръ ўЄю ё єўхЄюь чрьхўрэш  \ref{nt:gen2}, ь√ ЁрёяЁюёЄЁрэ хь Ёхчєы№ЄрЄ  ЁрсюЄ√ \cite{VlaSh}  эр сюыхх °шЁюъшщ ъырёё
фшЇЇхЁхэЎшры№э√ї т√Ёрцхэшщ. ─юяюыэшЄхы№э√х єёыютш , тючэшъ°шх т  \cite{ShVl} юс· ёэ ■Єё  шёяюы№чютрэшхь чрьхэ√
яхЁхьхээющ $\xi=p(x)$ ш яЁшьхэхэшхь фЁєуюую ьхЄюфр фы  фюърчрЄхы№ёЄтр юёэютэющ ЄхюЁхь√. ╬ЄьхЄшь х∙х, ўЄю т ЁрсюЄрї
ртЄюЁют \cite{SavSh03} (фы  ёыєўр  $l(y)=-y''+\tau_1y$) ш \cite{SavSh14} (фы  ёшёЄхь√ ─шЁрър) яЁшьхэ ыё  х∙х юфшэ ьхЄюф
--- ьхЄюф "<єуыр ╧Ё■ЇхЁр">. ╨хчєы№ЄрЄ√ юс рёшьяЄюЄшўхёъюь яютхфхэшш ╘╤╨, яюыєўхээ√х т ¤Єшї ЁрсюЄрї
ёютярфр■Є ё \eqref{eq:matras1} ш \eqref{eq:matras2}. ┬ ёыєўрх юяхЁрЄюЁр ╪ЄєЁьр--╦шєтшыы  ярЁрьхЄЁ $\mu$ т√сшЁрыё 
Ёртэ√ь $2$, р фы  ёшёЄхь√ ─шЁрър с√ыр ЁрёёьюЄЁхэр °ърыр $\mu\in[2,\infty]$). ┬ \cite{SavSh03}, юфэръю, ъЁюьх
"<ъюЁюЄъющ"> рёшьяЄюЄшъш фы  ёюсёЄтхээ√ї чэрўхэшщ ш ёюсёЄтхээ√ї ЇєэъЎшщ, с√ыр (фы  ёыєўр  ъЁрхт√ї єёыютшщ ─шЁшїых)
яюыєўхэр ш "<фышээр "> рёшьяЄюЄшър, юёэютрээр  эр рёшьяЄюЄшўхёъюь яЁхфёЄртыхэшш \eqref{eq:matras6}. ╧Ёш ¤Єюь т шЄюуют√щ
юЄтхЄ эх тю°ыш ёырурхь√х, яюЁюцфрхь√х ЇєэъЎш ьш $\Z^3$ ш $\Z^4$, р юЄ ЇєэъЎшш $\Z^2$ тю°ыю Єюы№ъю юфэю ёырурхьюх.
╧Ёшўшэ√ ¤Єюую  тыхэш  ёяЁ Єрэ√ фютюы№эю уыєсюъю ш ёюёЄю Є т Єюь, ўЄю, їюЄ  ёрьш ЇєэъЎшш $\|\Z^j(\cdot,\la)\|_{L_2}$,
$j=3,\,4$, эх юЎхэштр■Єё  тхышўшэющ $\Upsilon_2^2(\la)$, эю эюЁь√ ¤Єшї ЇєэъЎшщ т яЁюёЄЁрэёЄтх ╒рЁфш $H^1(\wt\Gamma)$
єцх фюяєёър■Є  Єръшх цх яю яюЁ фъє юЎхэъш, ъръ ш $\|\Upsilon_2^2(\la)\|$. ╧юфЁюсэ√щ рэрышч ¤Єюую сєфхЄ яЁютхфхэ т
фЁєующ ЁрсюЄх ртЄюЁют.

\end{document}